\documentclass[a4paper]{article}

\usepackage[english]{babel}
\usepackage[utf8]{inputenc}
\usepackage[T1]{fontenc}

\usepackage[a4paper,top=3cm,bottom=2cm,left=3cm,right=3cm,marginparwidth=1.75cm]{geometry}

    \usepackage{amsthm}
    \usepackage{amssymb}
    \usepackage{amsmath}
    \usepackage{mathtools}
    \usepackage{graphicx}
    \usepackage{caption}
    \usepackage[colorinlistoftodos]{todonotes}
    \usepackage[colorlinks=true, allcolors=blue]{hyperref}
    \usepackage{xcolor}
    \usepackage{pinlabel}
    \usepackage{geometry}
    \usepackage{subcaption}
    \usepackage{hyperref}
    \usepackage{longtable}

\theoremstyle{definition}

\newtheorem*{dimostrazione}{Proof}

\theoremstyle{plain}
\newtheorem{teorema}{Theorem}[section]

\newtheorem{proposizione}{Proposition}[section]

\theoremstyle{remark}
\newtheorem{osservazione}{Remark}

\usepackage[backend=biber,
             style=nature,
             sorting=nty,
             ]{biblatex}
\title{An algorithmic method to compute plat-like Markov moves for genus two \(3\)-manifolds}
\author{Paolo Cavicchioli}
\date{September 2021}

\begin{document}

\maketitle

\begin{abstract}
This article  deals with equivalence of links in \(3\)-manifolds of Heegaard genus 2. Starting from a description of such a manifold introduced in \cite{casali19912},  that uses  \(6\)-tuples of integers  and determines a Heegaard decomposition of the manifold, we construct an algorithm (implemented in c++) which allows to find the words in \(B_{2, 2n}\), the braid group on \(2n\) strands of a surface of genus 2, that realizes the plat-equivalence for links in that manifold. In this way we extend to the case of genus 2 the result obtained in \cite{cattabriga2018markov} for genus 1 manifolds. We describe explicitly the words for a notable group of manifolds. 
\end{abstract}

\section{Introduction}
The representation of links  as plat closure of braids 
dates back to the work of J. Birman \cite{birman1976heegaard} that in 1976 proved that each link in $\mathbb R^3$ could be represented as the plat closure of a braid and described a set of moves that connect braids having isotopic closures. Since then, a lot of work has been done using this representation in the direction of studying the equivalence problem of links, for example by defining and analyzing link invariants  as, for example, the  Jones polynomial (see \cite{birman1988jones, bigelow2002does}). Using Heegaard surfaces, in \cite{doll1993generalization} H. Doll introduced the notion of $(g,b)$-decomposition or generalized bridge decomposition for links in a  closed, connected and orientable \(3\)-manifold, opening the way to the study of links in 3-manifolds via surface braid groups (see for example \cite{bellingeri2012hilden,cattabriga2004all,cattabriga2006alexander, cattabriga2008extending, cristofori2007strongly, goda2005knot}). Recently, the equivalence of  links in 3-manifolds, under isotopy,  has been described in terms of $(g,b)$-decomposition: in \cite{cattabriga2018markov}, the authors find a finite set of moves that connect braids in \(B_{g, 2n}\), the   braid group on \(2n\) strands of a surface of genus $g$,  having isotopic plat closures. In their result some  moves are explicitly described as elements in \(B_{g, 2n}\), while others, called plat slide moves,   depend on a Heegaard surface for the manifold and are explicitly described only for the case of Heegaard genus\footnote{We recall that the Heegaard genus of a closed, connected and orientable 3-manifold the minimum genus of an embedded Heegaard surface.} at most 1, that is for $S^3$, lens spaces and $S^2\times S^1$.   The goal of this article is to determine plat slide moves for links in  genus 2  3-manifolds.  This class of manifolds is quite interesting and vast: for example it includes  remarkable manifolds, as the Poincaré sphere, the Heisenberg manifold or the Weeks manifold and a various class of manifolds as small Seifert manifolds \cite{boileau1991genus}. Moreover, the classification problem for these manifold is still open.  


Our approach uses a representation introduced in \cite{casali19912} and connected with crystallization theory. In the article the authors prove that every 3-manifold of Heegaard genus 2 can be represented using a colored graph depending on a 6-tuple of integer numbers that detrmines the Heegaard splitting; of course different graphs  may represent the same manifolds (the problem of connecting 6-tuples giving homeomorphic manifolds is investigated in \cite{grasselli20002}). Starting from the 6-tuple, we construct an algorithm (also implemented on the computer) that allows to describe the plat slide moves  associated to the Heegaard surface in terms of elements in \(B_{2, 2n}\) and so to determine explicitly the  moves realizing the equivalence of links in the associated 3-manifold. 
This result can be used to define or compute link invariants or, more generally, to classify families of links in 3-manifolds of genus 2. Another interesting direction of development is to explore  the connection of this equivalence moves with those introduced in \cite{diamantis2015braid} using a rational surgery description for the 3-manifold and representing links via mixed braids. \\

In Section \ref{Preliminaries} we will recall some notion  and result about Heegaard splittings, crystallisation  and  generalized plat decomposition. \\
In Section \ref{genus2} we  describe the representation through 6-tuples of 3-manifolds of Heegaard genus 2 and we will obtain, directly from the 6-tuple, an open Heegaard diagram for the 3-manifold. \\
Section \ref{main} contains the proof of the main result of the article, that is the algorithm determining the equivalence moves cited above. \\
Finally, in Section \ref{examples} we will report the  moves for a notable group of 6-tuples analyzed in \cite{bandieri2010census} obtained using  the coded program. 

\section{Preliminaries}\label{Preliminaries}
In this section we recall the notion of crystallization and Heegaard splitting for 3-manifolds and their relationship, and the definition of generalized bridge representation for links in \(3\)-manifolds. \\ 
Manifolds are always assumed to be closed, connected and orientable and links inside them will be considered up to isotopy. \\

\vspace{-8pt}

\begin{figure}[h!]
    \centering
    \includegraphics[width = .6\textwidth]{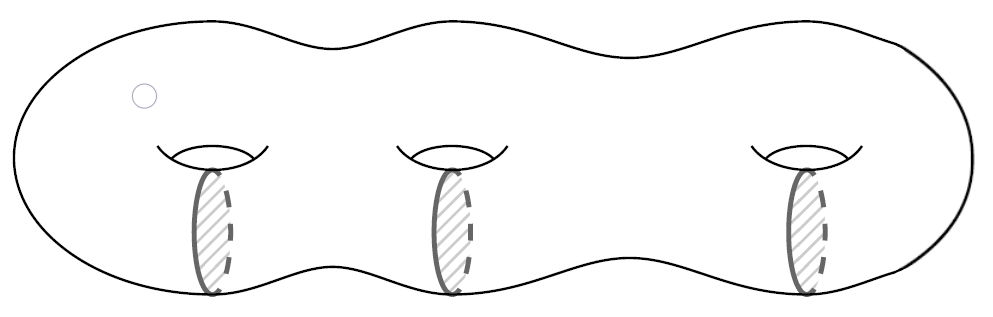}
    \caption{Standard genus \(g\) handlebody.}
    \label{fig:standard}
\end{figure}

A \emph{Heegaard surface} for a 3-manifold \(M\) is a connected closed orientable surface \(\Sigma\) embedded in \(M\) such that \(M\setminus \Sigma\) is the disjoint union of two handlebodies (of the same genus). So we have that \(M\) is homeomorphic to \(H_1 \cup_h H_2\), where \(H_1\) and \(H_2\) are two oriented copies of a standard handlebody in \(\mathbb R^3\) (see Fig. \ref{fig:standard}) and \(h:\partial H_2 \rightarrow \partial H_1\) is an orientation reversing homeomorphism. The triple \((H_1,H_2,h)\) is called \emph{Heegaard splitting} of \(M\) and the Heegaard surface is \(\partial H_1\cup_h\partial H_2\).   
Each \(3\)-manifold admits Heegaard splittings \cite{heegaard1898forstudier}. The \emph{Heegaard genus} of a 3-manifold \(M\) is the minimal genus of a Heegaard surface for \(M\). \\ 
For instance, the \(3\)-sphere \(S^3\) is the only \(3\)-manifold with Heegaard genus \(0\), while the manifolds with Heegaard genus \(1\) are lens spaces (i.e., cyclic quotients of \(S^3\)) and \(S^2\times S^1\). While \(S^3\), as well as lens spaces, have, up to isotopy, only one Heegaard surface of minimal genus (and those of higher genera are stabilizations of that of minimal genus), in general, a manifold may admit non isotopic Heegaard surfaces of the same genus (see \cite{moriah1988heegaard} for the case of Seifert manifolds).\\
Each Heegaard splitting can be represented by means of a Heegaard diagram. Given \((H_1, H_2, h)\) a Heegaard splitting for \(M\) and two sets of meridians\footnote{We recall that a \emph{set of meridians} of \(H_g\), handlebody of genus \(g\), is a collection of closed curves \(u = \{u_1, \dots, u_g\}\) bounding properly embedded disjoint discs \(D_1, \dots, D_g\) such that \(H_g \setminus D_1, \dots, D_g\) is a 3-ball.} \(\mathbf{c} = \{c_1, \dots, c_g\}, \> \mathbf{c}^* = \{c^*_1, \dots, c^*_g\}\) for \(H_1\) and \(H_2\) respectively on \(\Sigma_g = \partial H_1 \cup_h \partial H_2\), we call \((\Sigma_g, \mathbf{c}, \mathbf{c}^*)\) a \emph{closed Heegaard diagram} for \(M\). If we cut \(\Sigma_g\) along \(\mathbf{c}\), we obtain a sphere with \(2g\) holes, say \(D_1, \dots, D_{2g}\), distinct and paired, each pair corresponding to a certain \(\mathbf{c}_i\). Therefore the meridians \(\mathbf{c}^*_i\) will be naturally cut into arcs joining the holes in various ways, giving a graph on the sphere with \(2g\) holes called \emph{open Heegaard diagram} of \((\Sigma_g, \mathbf{c}, \mathbf{c}^*)\). \\

\vspace{-8pt}

Another way to represent 3-manifolds are crystallizations. Let us recall the definition. 
An \emph{edge-coloring on a graph} \(\Gamma = (V(\Gamma), E(\Gamma))\) is a map \(\gamma: E(\Gamma) \rightarrow \Delta_{n+1} = \{0, 1, ..., n\}\) such that \(\gamma(e) \neq \gamma(f)\), for each pair of adjacent edges \(e, f\). If \(v, w\) are the vertices of an edge \(e \in E(\Gamma)\) such that \(\gamma(e)=c\), we say that \(e\) is a \emph{\(c\)-edge} and that \(v, w\) are \emph{\(c\)-adjacent}. Let \(\Gamma\) be a regular graph of degree \(n+1\) (i.e., each vertex has degree \(n+1\)) and \(\gamma: E(\Gamma)\rightarrow \Delta_{n+1}\) be an edge-coloring, then the pair \((\Gamma, \gamma)\) is said to be an \emph{\((n+1)\)-coloured graph}. For each \(\mathcal{B} \subseteq \Delta_{n+1}\), we set \(\Gamma_\mathcal{B}=(V(\Gamma), \gamma^{-1}(\mathcal{B}))\); moreover, each connected component of \(\Gamma_\mathcal{B}\) will be called a \(\mathcal{B}\)-residue. We denote with \(\rho_{\mathcal{B}}\) the number of \(\mathcal{B}\) residues in \((\Gamma, \gamma)\). For each colour \(c \in \Delta_{n+1}\), we set \(\hat{c} = \Delta_{n+1} - \{c\}\). \\
In \cite{ferri1986graph} it is proven that every \((n+1)\)-coloured graph \(\Gamma\) represents an \(n\)-dimensional pseudomanifold \(M(\Gamma)\) which is orientable if and only if \(\Gamma\) is bipartite. A \emph{crystallization} \(\Gamma\) is an \((n+1)\)-coloured graph such that \(\Gamma_{\hat{c}}\) is connected, for every \(c \in \Delta_{n+1}\) and \(M(\Gamma)\) represents an \(n\)-manifold. Every \(n\)-manifold admits a crystallization through a \((n+1)\)-colored graph \cite{pezzana1974sulla}, in particular we can represent every \(3\)-manifold with a \(4\)-colored graph. \\

\vspace{-8pt}

Given a crystallization of a 3-manifold, it is possible to recover a Heegaard for it. Indeed, let \(\Gamma\) be a 4-coloured graph, and let \(c, d\) distincts colors of \(\Delta_4\) and  denote with \(\{c', d'\} = \Delta_4 - \{c, d\}\); in \cite{gagliardi1981regular} it is proven that if \(\Gamma\) is a crystallization of a 3-manifold \(M\), for any choice of a couple of colors \(c, d\), there exists an embedding \(\mathbf{i}\) of \(\Gamma\) into a surface \(F_{cd}\) of genus \(\rho_{c'd'} - \rho_{\hat{c}} - \rho_{\hat{d}} + 1\) that is a Heegaard surface for \(M\). 
Also be \(\xi^0, \dots, \xi^n\) the components of \(\Gamma_{\{c, d\}}\) and \(\xi^{0'}, \dots, \xi^{n'}\) the ones of \(\Gamma_{\{c', d'\}}\). Set \(x^i = \mathbf{i}(\xi^i)\) and \(y^i = \mathbf{i}(\xi^{i'})\) for \(i = 0, \dots, n\) and denote with \(\mathbf{x} = \{x^0, \dots, x^n\}, \> \mathbf{y} = \{y^0, \dots, y^n\}\), \(\mathbf{x}(\hat{j}) = \mathbf{x} \setminus \{x^j\}\), \(\mathbf{y}(\hat{j}) = \mathbf{y} \setminus \{y^j\}\). It is proven in \cite{gagliardi1981extending}, the triple \((F_{\{c, d\}}, \mathbf{x}(\hat{j}), \mathbf{y}(\hat{j}))\) is a Heegaard diagram. \\

\vspace{-8pt}

We end the session by recalling the definition of generalized bridge position for links in \(3\)-manifolds \cite{doll1993generalization}. \\
\vspace{-8pt}

\begin{figure}[h!]
\labellist
\small\hair 2pt
\pinlabel \(1\) at 120 110
\pinlabel \(g\) at 500 110
\pinlabel \(i\) at 310 110
\pinlabel \(1\) at 580 110
\pinlabel \(j\) at 650 110
\pinlabel \(j+1\) at 745 65
\pinlabel \(2n\) at 790 110
\pinlabel \(a_i\) at 420 120
\pinlabel \(b_i\) at 420 10
\pinlabel \(\sigma_j\) at 680 130
\endlabellist
    \centering
    \includegraphics[width = .6\textwidth]{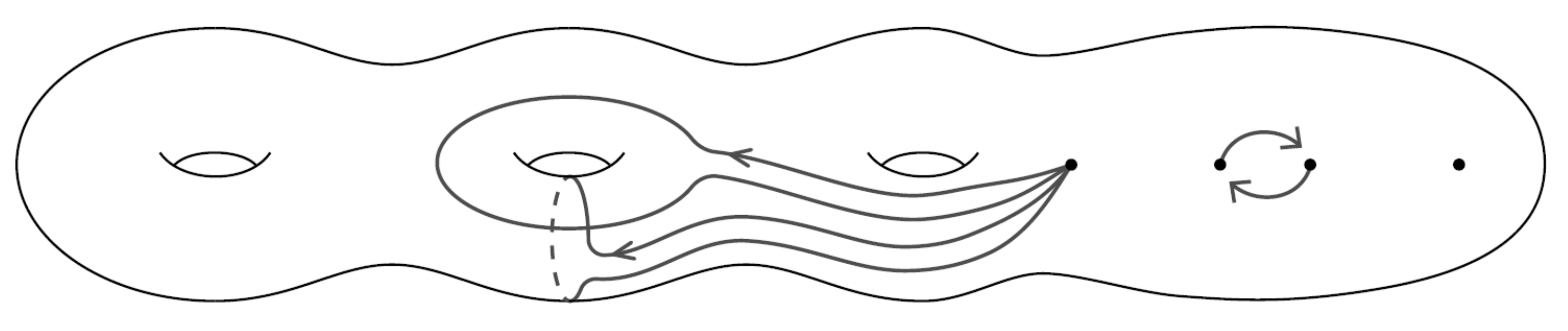}
    \caption{The generators of \(B_{g, 2n}.\)}
    \label{fig:presentation}
\end{figure}
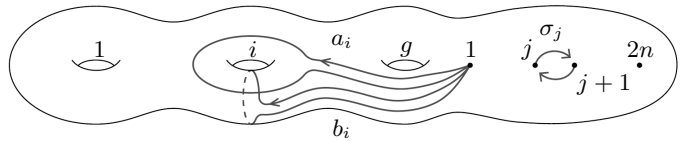

First, we need to say that, in a handlebody \(H\), a set of \(n\) properly embedded disjoint arcs \(\{A_1,\ldots ,A_n\}\) is a \emph{trivial system} of arcs if there exist \(n\) mutually disjoint embedded discs, called \emph{trivializing discs}, \(D_1,\ldots ,D_n\subseteq H\) such that \(A_i\cap D_i= A_i\cap\partial D_i=A_i\), \(A_i\cap D_j=\emptyset\) and \(\partial D_i - A_i\subset\partial H\) for all \(i,j=1,\ldots ,n\) and \(i\neq j\). In this case, we say that the arc \(A_i\) \emph{projects} onto the arc \(\partial D_i-\textup{int}(A_i)\subset\partial H\) via \(D_i\). \\
Let \(\Sigma\) be a Heegaard surface for \(M\). We say that a link \(L\) in \(M\) is in \emph{bridge position} with respect to \(\Sigma\) if:  
\begin{enumerate}
\item \(L\) intersects \(\Sigma\) transversally and
\item the intersection of \(L\) with both handlebodies, obtained  splitting \(M\) by  \(\Sigma\), is a trivial system of arcs.
\end{enumerate}
Such a decomposition for \(L\) is called \((g,n)\)\emph{-decomposition} or \(n\)-\emph{bridge decomposition of genus} \(g\), where  \(g\) is the genus of \(\Sigma\) and \(n\) is the cardinality of the trivial system. \\
The minimal \(n\) such that \(L\) admits a \((g,n)\)-decomposition is called the \emph{genus} \(g\) \emph{bridge number} of \(L\). If \(g=0\), the manifold \(M\) is the 3-sphere and  we get the usual notion of bridge decomposition and bridge number of links in the 3-sphere (or in \(\mathbb R^3\)). \\
Following \cite{cattabriga2018markov} we describe a way to represent links in \(3\)-manifolds via braids. Let \(\mathcal P_{2n}=\{P_1,\ldots,P_{2n}\}\) be a set of \(2n\) distinct points  on \(\Sigma_g\) and denote with  \(B_{g,2n}\), the braid group  on \(2n\)-strands of \(\Sigma_g\), i.e, the fundamental group  of the configuration space of the \(2n\) points in \(\Sigma_g\).  Referring to Fig.\ref{fig:presentation}, the group   \(B_{g,2n}\) is generated  by    \(\sigma_1,\ldots,\sigma_{2n-1}\), the standard braid generators, and \(a_1,\ldots,a_g,b_1,\ldots, b_g\), where \(a_i\) (resp. \(b_i\)) is the braid whose strands are all trivial except the first one which goes once along the \(i\)-th longitude (resp. \(i\)-th meridian) of \(\Sigma_g\) (see  \cite{bellingeri2012hilden}).\\ Fix a set of \(n\) arcs \(\gamma_1,\ldots,\gamma_n\) embedded into \(\Sigma_g\), such that \(\gamma_i\cap \gamma_j=\emptyset\) if \(i\ne j\) and \(\partial \gamma_i=\{P_{2i-1},P_{2i}\}\), for \(i,j=1,\ldots,n\).  Given an element \(\beta\in B_{g,2n}\),  realize it as a geometric braid, that is, as a set of \(2n\) disjoint paths in \(\Sigma_g\times\left[ 0,1\right]\) connecting   \(\mathcal P_{2n}\times \left\{0\right\}\) to \(\mathcal P_{2n}\times \left\{1\right\}\).  The  \emph{plat closure} \(\widehat{\beta}\subseteq M\) of \(\beta\)  is  the link obtained  ``closing'' \(\beta\) by connecting   \(P_{2i-1}\times\{0\}\) with  \(P_{2i}\times\{0\}\) through   \(\gamma_i\times\{0\}\) and \(P_{2i-1}\times\{1\}\) with  \(P_{2i}\times\{1\}\)  through  \(\gamma_i\times\{1\}\), for \(i=1,\ldots,n\). Clearly,  \(\widehat{\beta}\) is in bridge position with respect to \(\Sigma_g\) and so it has genus \(g\) bridge number at most \(n\). 
For each link \(L\) in a \(3\)-manifold \(M\), there exists a braid \(\beta \in \cup_{n \in \mathbb{M}} B_{g, 2n}\) such that \(\hat{\beta}\) is isotopic to \(L\). \\

Moreover, as proved  in  \cite{cattabriga2018markov}, two braids \(\beta_1, \beta_2 \in \bigcup_{n\in\mathbb N} B_{g,2n}\) have isotopic plat closure if and only if \(\beta_1\) and \(\beta_2\) differ by braid isotopy,   a finite sequence of the following moves
\begin{align*}
  (M1)  & \qquad \sigma_1\beta \longleftrightarrow  \beta  \longleftrightarrow \beta \sigma_1\\
  (M2)  & \qquad \sigma_{2i}\sigma_{2i+1}\sigma_{2i-1}\sigma_{2i} \beta \longleftrightarrow \beta \longleftrightarrow \beta \sigma_{2i}\sigma_{2i+1}\sigma_{2i-1}\sigma_{2i} \\
  (M3)  & \qquad \sigma_2\sigma_1^2\sigma_2\beta \longleftrightarrow  \beta  \longleftrightarrow \beta \sigma_2\sigma_1^2\sigma_2\\
   (M4)  & \qquad a_j\sigma_1^{-1}a_j\sigma_1^{-1}\beta \longleftrightarrow \beta\longleftrightarrow \beta a_j\sigma_1^{-1}a_j\sigma_1^{-1}\quad \textup{for } j=1,\ldots, g \\
  (M5) &\qquad b_j\sigma_1^{-1}b_j\sigma_1^{-1}\beta \longleftrightarrow \beta\longleftrightarrow \beta b_j\sigma_1^{-1}b_j\sigma_1^{-1}\quad \textup{for } j=1,\ldots, g \\
  (M6) &\qquad \beta \longleftrightarrow T_k(\beta)\sigma_{2k}\\
  \ &\ \ \quad \textup{ where } T_k: B_{g,2n}\to B_{g, 2n+2} \textup{ is  defined by } 
   T_k(a_i)=a_i, \ T_k(b_i)=b_i \textup{ and }\\
\ & \ \ \quad T_k(\sigma_i)=\left\{\begin{array}{l}\sigma_i \qquad \qquad \qquad \qquad \qquad \qquad\quad\textup{if } i<2k\\
\sigma_{2k}\sigma_{2k+1}\sigma_{2k+2}\sigma_{2k+1}^{-1}\sigma_{2k}^{-1}\ \qquad \quad \textup{if } i=2k\\
\sigma_{i+2} \qquad \qquad \qquad \qquad \qquad \qquad \textup{if } i>2k\\
\end{array}\right..
\end{align*}
and two sets of moves called \emph{plat slide moves} and \emph{dual plat slide moves}  that depend on the Heegaard decomposition of $M$ and are defined as follows
\[psl_i : \beta \longrightarrow \bar{c}_i \# \beta \in B_{g,2n}, \> i = 1, \dots, g\]
\[psl_i^* : \beta \longrightarrow \beta \# \bar{c}_i^* \in B_{g,2n}, \> i = 1, \dots, g\]
where $(\Sigma_g, \mathbf c, \mathbf c^*)$ is a Heegaard diagram for $M$ and for \(\alpha \in B_{g, 2n},\> \beta \in B_{g, 2m}\)  the \emph{plat sum} of \(\alpha\) and \(\beta\) is defined as  \(\alpha \# \beta := \alpha w_{n,m} \beta \in B_{g, 2(n+m-1)}\)
with
\[w_{n,m} = \prod_{i = 0}^{2n-3} \prod_{j = 0}^{2m-3} \sigma_{2n-i+j}.\]

If we assume that \(\mathbf{c}^*\) is the system corresponding to the curves depicted in Fig.\ref{fig:standard}, then we have \(c_i^* = \widehat{b_i}\) and so the dual plat slide move $psl_i^*$ is 
\[psl^*_i : \beta \longrightarrow \beta \# b_i = \beta b_i\]

In \cite{cattabriga2018markov} an explicit formula for plat slide moves in \(L(p,q)\) (i.e., the case of genus one) is
\[psl_{p,q}:\; \beta \longrightarrow \alpha_{p,q} \# \beta = \alpha_{p,q} \, \beta,\]
where
\[\alpha_{p,q} = \prod_{i=1}^{r} b_1^{-1} a_1^{\lceil \frac{p}{q}\rceil} \cdot \prod_{i=1}^{q-r} b_1^{-1} a_1^{\lfloor \frac{p}{q}\rfloor} \in B_{1,2}.\]

\section{Genus-\(2\) manifolds}
\label{genus2}

In this section  we describe, following \cite{casali19912, grasselli20002}, a way of representing all prime \(3\)-manifolds of genus \(g \leq 2\) via crystallizations defined by \(6\)-tuples of integers. \\
Let \(\mathcal{F}\) be the set of the \(6\)-tuples \(f = (h_0, h_1, h_2, q_0, q_1, q_2)\) of integers satisfying the following conditions: 
\begin{enumerate}
    \item \(h_i + h_{i+1} = 2l_{i+1}\), an even integer for \(i = 0, 1, 2\);
    \item \(0\leq q_i < 2l_i\), for \( i = 0, 1, 2\);  
    \item \label{cond_tre} all \(q_i\)'s have the same parity;
    \item \(h_i + q_i\) is odd for \(i = 0, 1, 2\)
    \item there exists at most one \(j\) such that \(h_j = 0\). 
\end{enumerate}
The parameter \(q_i\) will be considered \(\text{mod } 2l_i\) for \(i = 0, 1, 2\). \\
Let \(\mathcal{G} = \{\Gamma(f) | f \in \mathcal{F}\}\) be the class of \(4\)-coloured graphs \(\Gamma(f)\) whose vertices are the elements of the set \[V(f)= \bigcup_{i = 0, 1, 2} \{i\} \times \mathbb{Z}_{2l_i},\] and such that there exists a \(c\)-edge, with \(c = 0, 1, 2, 3\), connecting the vertices \((i,j)\) and \((i', j')\) if and only if \((i', j') = \tau_c (i, j)\) where \(\tau_c\) is the fixed-point-free involution defined by:
\begin{equation}
\begin{array}{lcl}
\tau_0(i, j) & = & (i, j+(-1)^j) \\
\tau_1(i, j) & = & (i, j-(-1)^j) \\
\tau_2(i, j) & = & \left\{\begin{array}{ll}
(i+1, -j-1) & \text{if } j = 0, \dots, h_i - 1 \\
(i-1, 2l_i -j-1) & \text{if } j = h_i, \dots, 2l_i - 1
\end{array} \right., \\
\tau_3(i, j) & = & \rho \tau_2 (i, j) \rho^{-1}
\end{array}
\end{equation}
with \(\rho : V(f) \rightarrow V(f)\)  the bijection defined by \(\rho(i, j) = (i, j+q_i)\). \\
Clearly \(\Gamma(f)\) is regular and is also bipartite because of condition \ref{cond_tre}, so the associated pseudocomplex is indeed a closed orientable \(3\)-pseudomanifold. \\
Note also that \(\Gamma(f)\) has exactly three \(\{0, 1\}\)-residues, that is \(\rho_{\{0, 1\}} = 3\). We denote them with \(C_0, C_1, C_2\), with \(C_i\) having as vertices the elements with \(i\) as first coordinate. In \cite{casali19912}, it is proven that \(\Gamma(f) \in \mathcal{G}\) is a crystallization of a \(3\)-manifold if and only if \(\rho_{\{2, 3\}} = 3\). We call \emph{admissible} a \(6\)-tuple belonging to \(\mathcal{F}\) satisfying this condition and denote with $D_0$, $D_1$ and $D_2$ the three \{2,3\}-residues of $\Gamma(f)$ and with   $M_f$ the 3-manifold associated to  $f$. \\ 

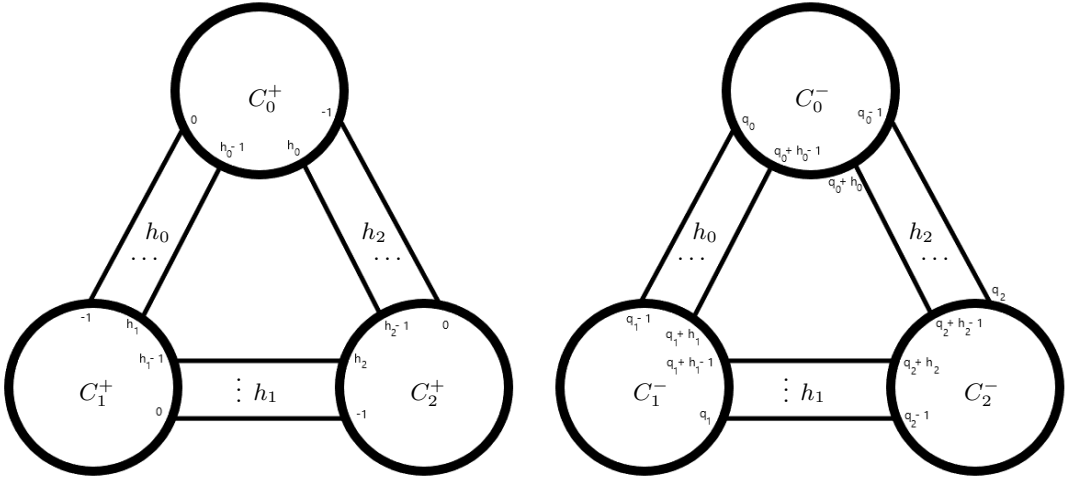
\begin{figure}[h!]
\labellist
\small\hair 2pt
\pinlabel \(C_1^+\) at 100 130
\pinlabel \(C_1^-\) at 510 130
\pinlabel \(C_2^+\) at 345 130
\pinlabel \(C_2^-\) at 755 130
\pinlabel \(C_0^+\) at 225 350
\pinlabel \(C_0^-\) at 630 350
\pinlabel \(h_1\) at 225 130
\pinlabel \(\vdots\) at 205 138
\pinlabel \(h_1\) at 630 130
\pinlabel \(\vdots\) at 610 138
\pinlabel \(h_0\) at 145 250
\pinlabel \(\dots\) at 135 230
\pinlabel \(h_2\) at 305 250
\pinlabel \(\dots\) at 315 230
\pinlabel \(h_0\) at 550 250
\pinlabel \(\dots\) at 540 230
\pinlabel \(h_2\) at 710 250
\pinlabel \(\dots\) at 720 230
\endlabellist
    \centering
    \includegraphics[width = \textwidth]{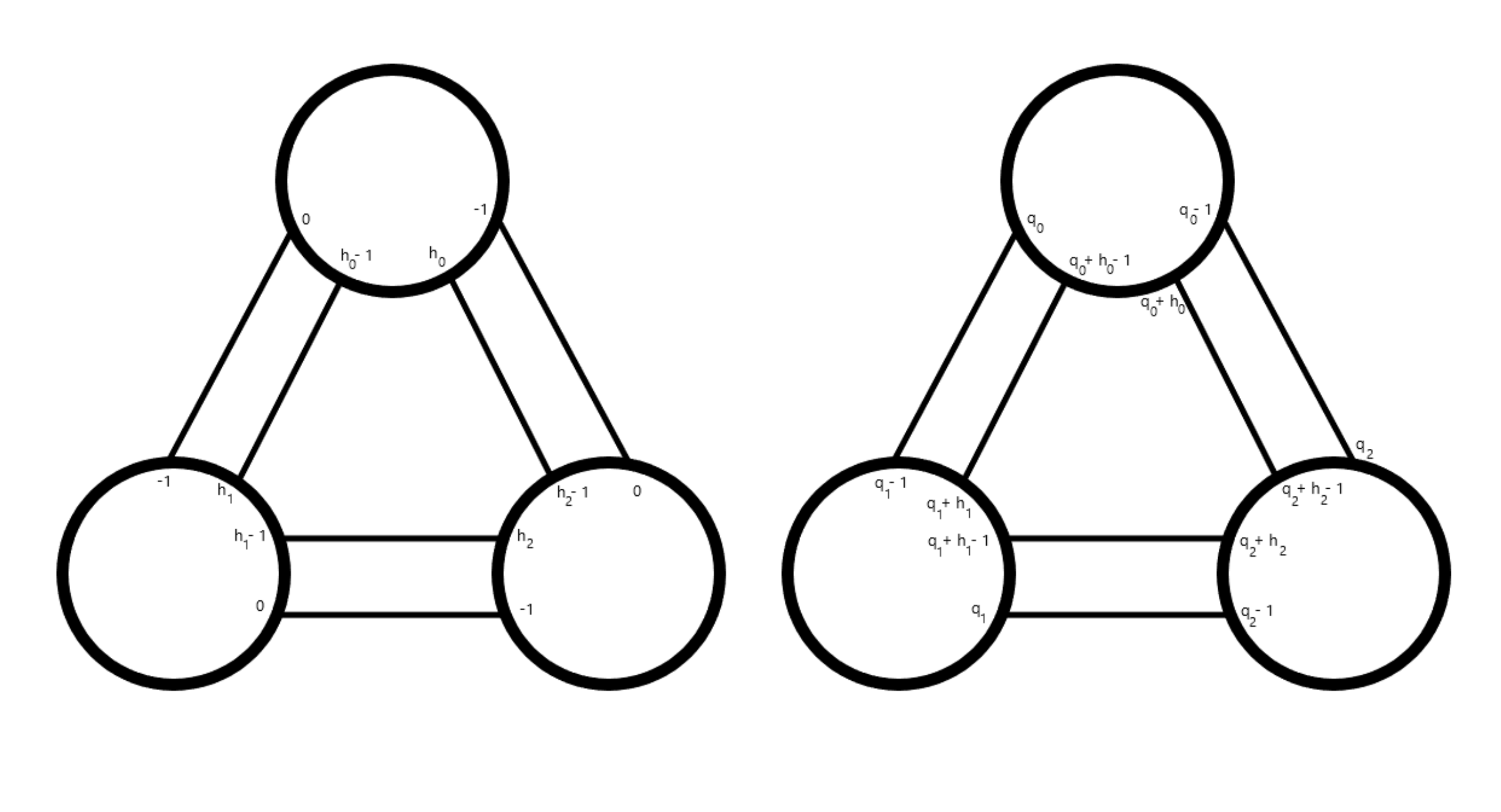}
    \caption{The residues \(\Gamma(f)_{\hat{3}}\), \(\Gamma(f)_{\hat{2}}\). All \(C_i\)'s have a counterclockwise orientation.}
    \label{fig:3-residues}
\end{figure}

As recalled in Section \ref{Preliminaries}, the manifold $M_f$ has genus $\rho_{\{0,1\}} - \rho_{\hat{2}} - \rho_{\hat{3}} + 1 = 3 - 1 - 1 + 1 = 2$. Moreover in order to get a genus $2$ Heegaard diagram of $M_f$ we can proceed as follows (see \cite{gagliardi1981extending}). Consider the planar representation of the two \(3\)-residues \(\Gamma(f)_{\hat{3}}\) and \(\Gamma(f)_{\hat{2}}\) and, referring to  Fig.\ref{fig:3-residues}, denote with \(C_0^+, C_1^+, C_2^+\) (respectively \(C_0^-, C_1^-, C_2^-\)) the copy of $C_0,\> C_1,\> C_2$ belonging to \(\Gamma(f)_{\hat{3}}\) (respectively \(\Gamma(f)_{\hat{2}}\)). Following Fig.\ref{fig:4610111_nedela} (that represents the case of $f=(4,6,10,1,1,1)$) we can embed \(\Gamma(f)_{\hat{3}}\) and \(\Gamma(f)_{\hat{2}}\) in a standard genus 2 surface, that is $F_{\{0,1\}}$, so that: (i) $C_0,C_1,C_2$ are the curves belonging to the intersection with the plane horizontal plane of symmetry for the surface, (ii) \(\Gamma(f)_{\hat{3}}\) is embedded in the upper part of the surface (and represented using violet arcs) and (iii)  \(\Gamma(f)_{\hat{2}}\) is embedded in the lower part (and represented using green arcs). 

\begin{figure}[h!]
\labellist
\small\hair 2pt
\pinlabel \(C_1\) at 390 330
\pinlabel \(C_0\) at 675 330
\pinlabel \(C_2\) at 70 330
\endlabellist
    \centering
    \includegraphics[width = 0.9\textwidth]{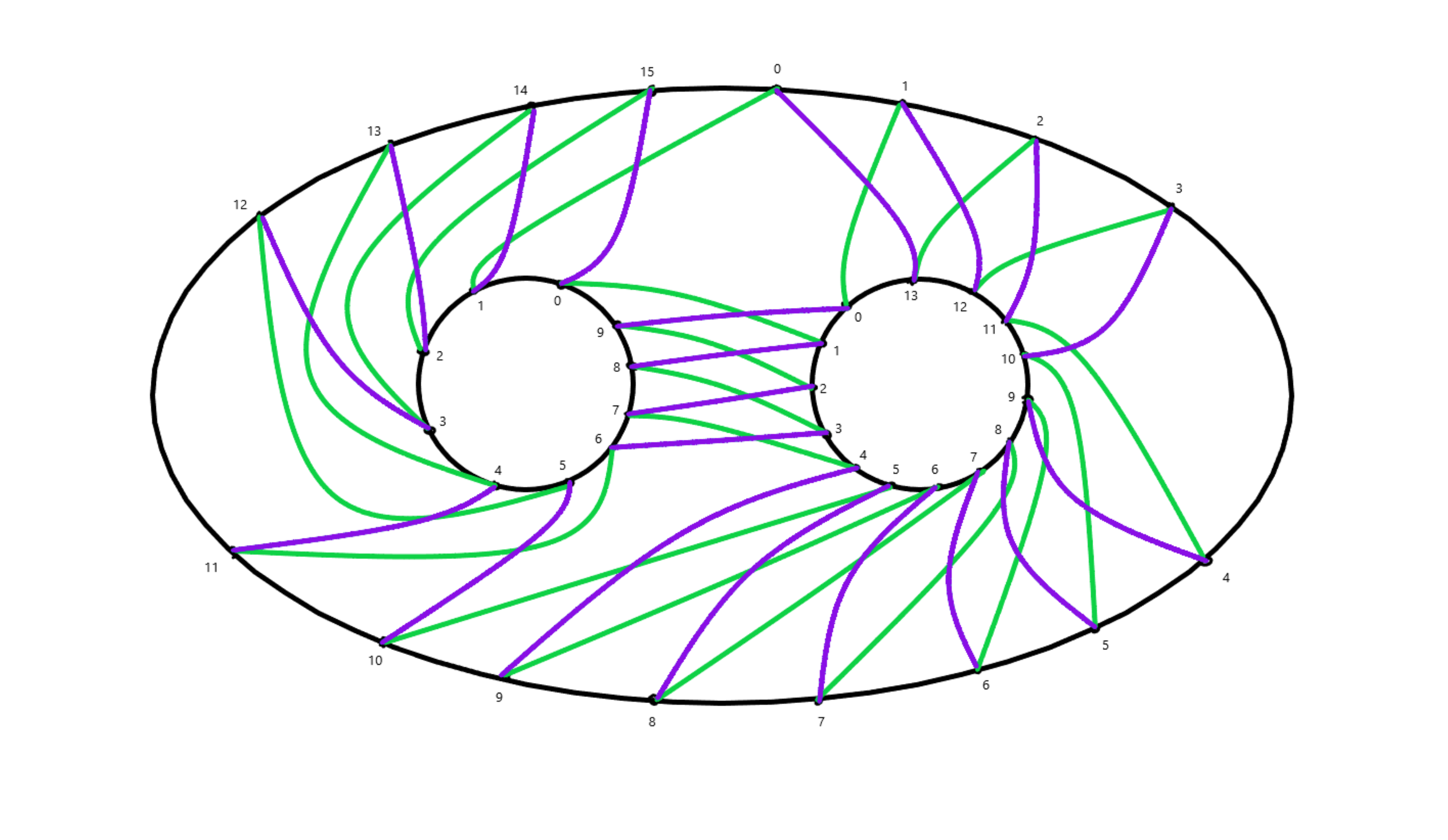}
    \caption{The genus 2 embedding of  \(f = (4, 6, 10, 1, 1, 1)\).}
    \label{fig:4610111_nedela}
\end{figure}

So as recalled in Section \ref{Preliminaries}, the choice of  a couple of curves $C_i$'s and  a couple of curves $D_i$'s gives a closed Heegaard diagram for $M_f$. We call \emph{rich Heegaard diagram}, and denote it with $R_f$ the colored graph obtained by, removing the cycle $C_2$, cutting the genus 2 surface along $C_0$ and $C_1$ and coloring with the same color the arcs belonging to the same $D_i$.  The name is due to the fact that, if we remove the arcs of any of the \(D_i\)'s we get an open Heegaard diagram for $M_f$. In the following proposition we describe explicitly in terms of the parameters of $f$ the rich Heegaard diagram $R_f$. In all figures an arc with label \(k\) denotes \(k\) parallel arcs. 

\begin{figure}[h!]
\begin{subfigure}{.5\textwidth}
\labellist
\small\hair 2pt
\pinlabel \(C_1^-\) at 135 115
\pinlabel \(C_1^+\) at 135 330
\pinlabel \(C_0^-\) at 365 115
\pinlabel \(C_0^+\) at 365 330
\pinlabel \(a\) at 250 350
\pinlabel \(a\) at 250 135
\pinlabel \(b\) at 250 250
\pinlabel \(b\) at 480 100
\pinlabel \(c\) at 385 230
\pinlabel \(d\) at 115 230
\endlabellist
\centering
\includegraphics[width = .8\textwidth]{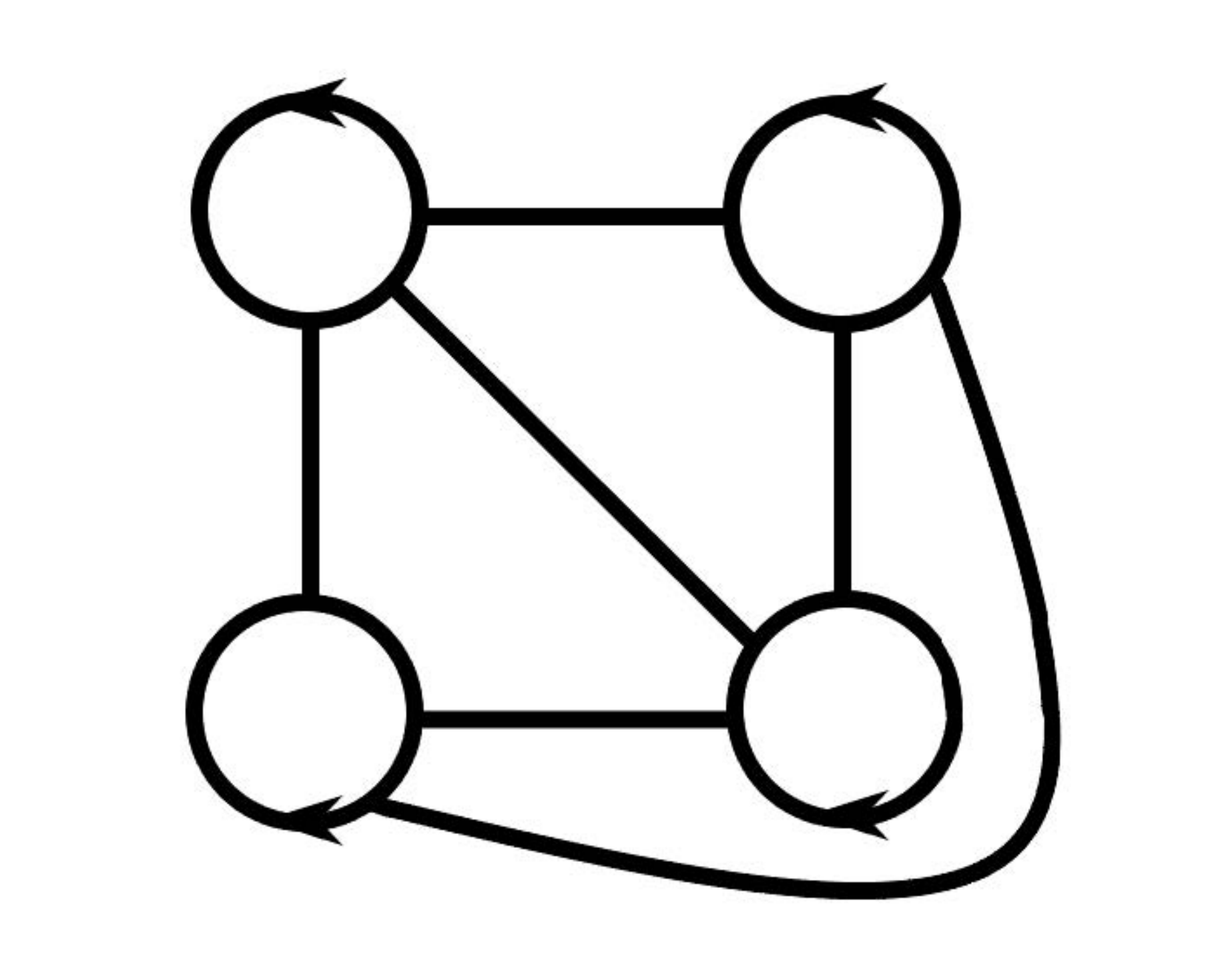}
\caption{}
\label{fig:heeg_1}
\end{subfigure}
\begin{subfigure}{.5\textwidth}
\labellist
\small\hair 2pt
\pinlabel \(C_1^-\) at 245 143
\pinlabel \(C_1^+\) at 245 425
\pinlabel \(C_0^-\) at 545 143
\pinlabel \(C_0^+\) at 545 425
\pinlabel \(a\) at 385 450
\pinlabel \(a\) at 385 165
\pinlabel \(b\) at 350 290
\pinlabel \(b\) at 90 500
\pinlabel \(c\) at 485 330
\pinlabel \(d\) at 215 290
\endlabellist
\centering
\includegraphics[width = .9\textwidth]{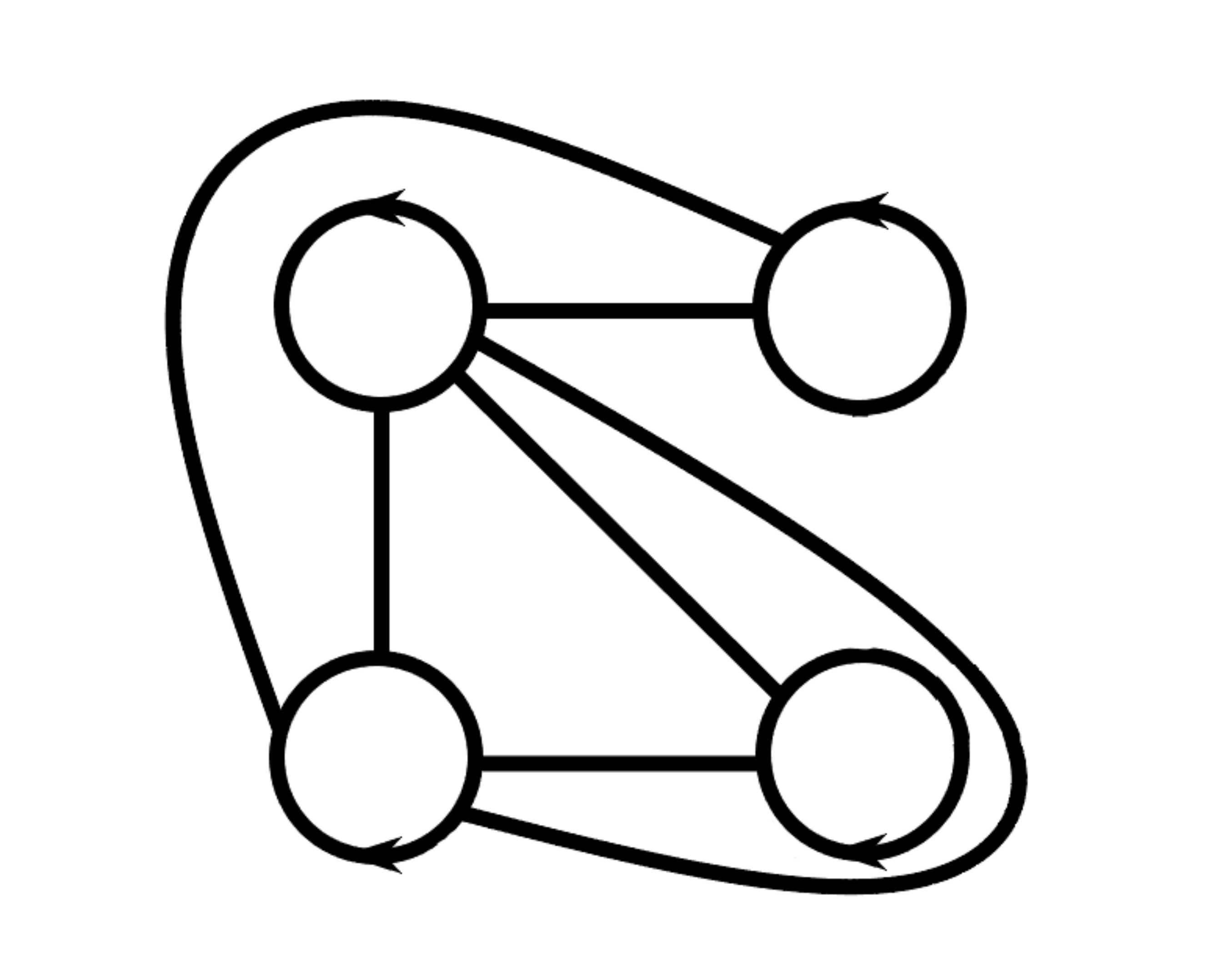}
\caption{}
\label{fig:heeg_2_giusto}
\end{subfigure}

\begin{subfigure}{.5\textwidth}
\labellist
\small\hair 2pt
\pinlabel \(C_0^+\) at 245 200
\pinlabel \(C_0^-\) at 245 448
\pinlabel \(C_1^+\) at 510 200
\pinlabel \(C_1^-\) at 510 448
\pinlabel \(a\) at 370 470
\pinlabel \(a\) at 370 225
\pinlabel \(b\) at 340 320
\pinlabel \(b\) at 110 500
\pinlabel \(c\) at 540 320
\pinlabel \(d\) at 215 320
\endlabellist
\centering
\includegraphics[width = \textwidth]{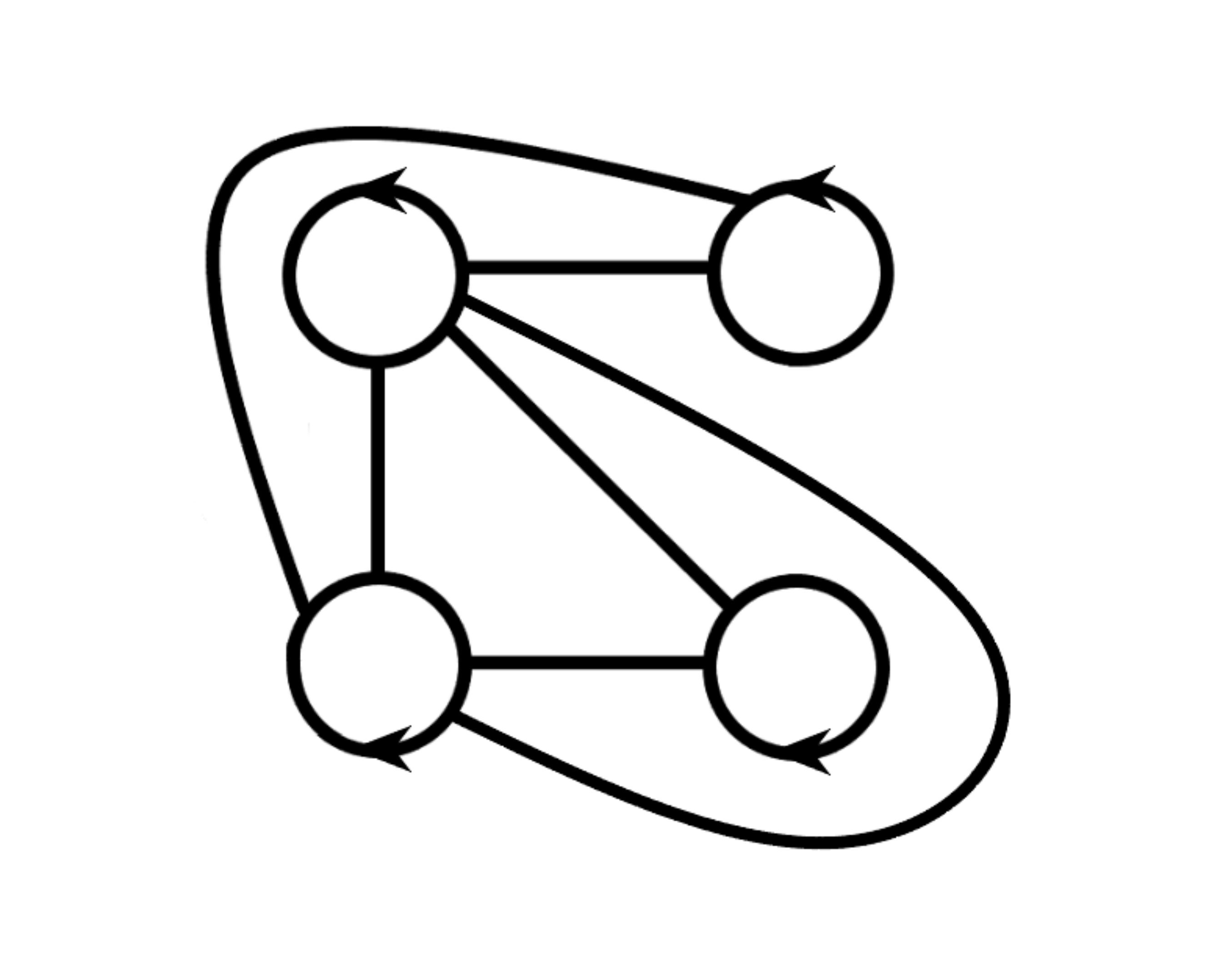}
\caption{}
\label{fig:heeg_3_giusto}
\end{subfigure}
\begin{subfigure}{.5\textwidth}
\labellist
\small\hair 2pt
\pinlabel \(C_1^-\) at 73 96
\pinlabel \(C_1^+\) at 73 226
\pinlabel \(C_0^-\) at 211 96
\pinlabel \(C_0^+\) at 210 226
\pinlabel \(a\) at 140 110
\pinlabel \(a\) at 140 240
\pinlabel \(b\) at 140 180
\pinlabel \(b\) at 280 70
\pinlabel \(d\) at 60 160
\pinlabel \(c\) at 198 160
\endlabellist
\centering
\includegraphics[width = .85\textwidth]{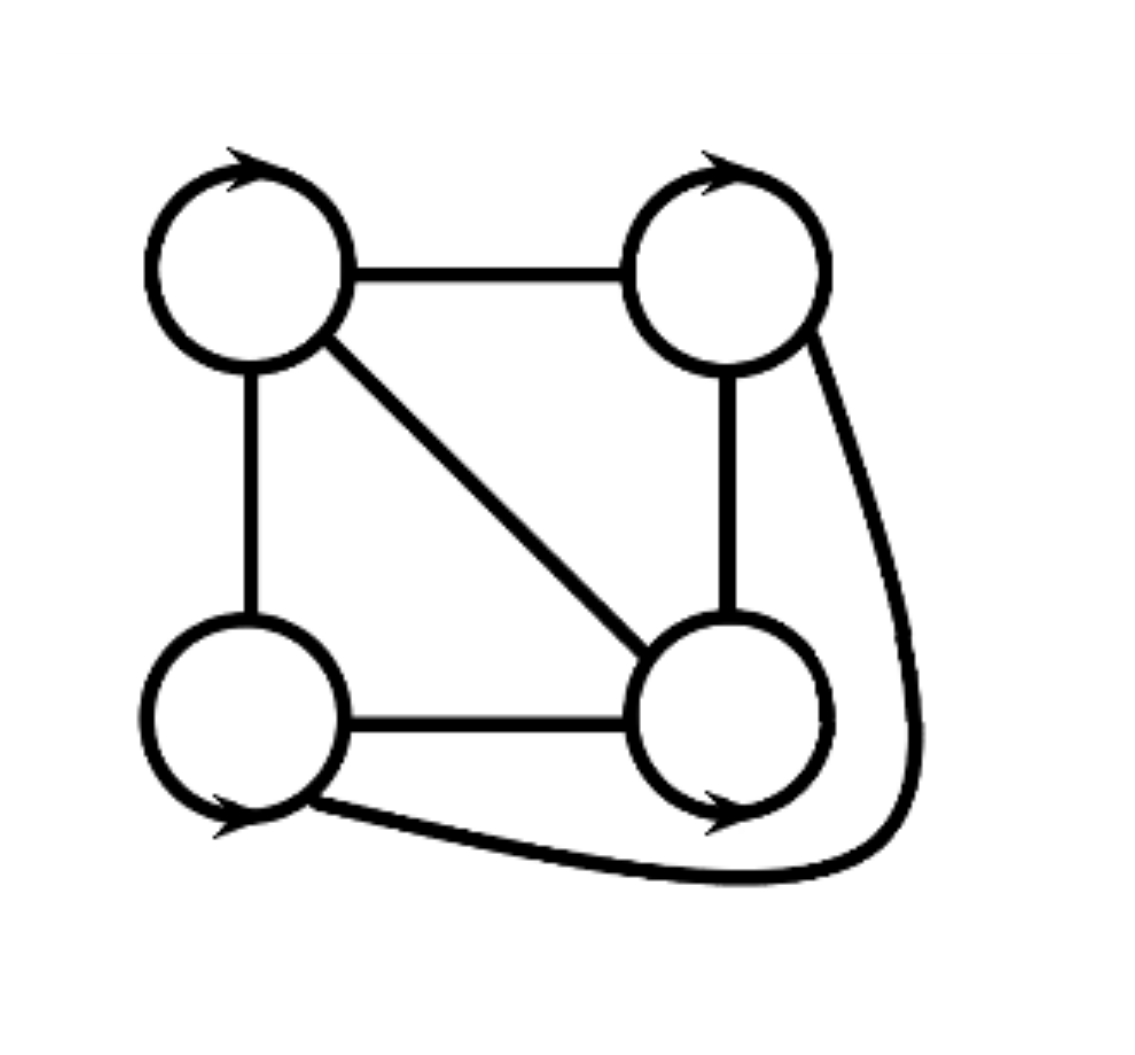}
\caption{}
\label{fig:heeg_4_giusto}
\end{subfigure}
\caption{The four possible cases.}
\label{heeg_generale}
\end{figure}
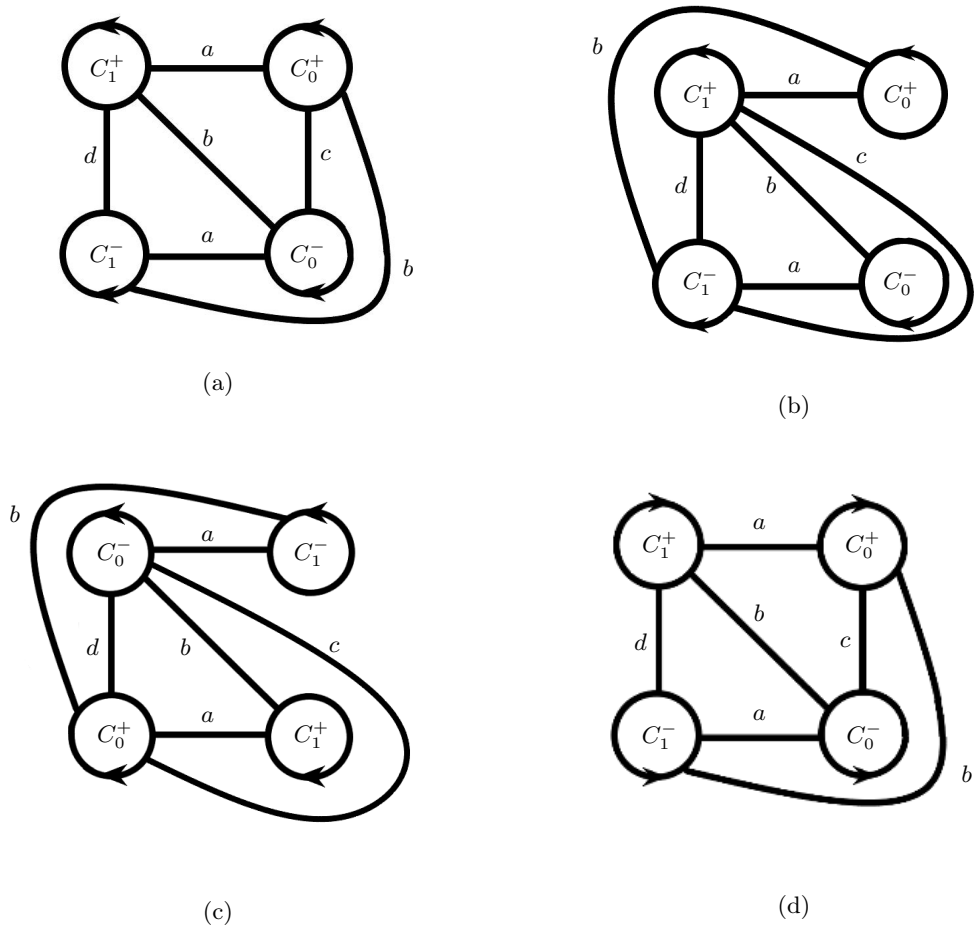

\begin{proposizione}\label{prop:dati}
Let \(f = (h_0, h_1, h_2, q_0, q_1, q_2)\) be an admissible \(6\)-tuple.  
\begin{enumerate}
    \item \label{prop:dati_1} If \(q_2 < h_1, h_2\) then the rich Heegaard diagram \(R_f\) is the one depicted in Fig.\ref{fig:heeg_1}, with \(a = h_0,\> b = q_2,\> c = h_2 - q_2,\> d = h_1 - q_2\); 
    \item \label{prop:dati_2} If \(h_2 \leq q_2 \leq h_1\) then the rich Heegaard diagram \(R_f\) is the one depicted in Fig.\ref{fig:heeg_2_giusto}, with \(a = h_0, b = h_2, c = q_2 - h_2, d = h_1 - q_2\); 
    \item \label{prop:dati_3} If \(h_1 < q_2 < h_2\) then the rich Heegaard diagram \(R_f\) is the one depicted in Fig.\ref{fig:heeg_3_giusto}, with \(a = h_0, b = h_1, d = h_2 - q_2, c = q_2 - h_1\); 
    \item \label{prop:dati_4} If \(q_2 > h_1, h_2\) then the rich Heegaard diagram \(R_f\) is the one depicted in Fig.\ref{fig:heeg_4_giusto}, with \(a = h_0, b = h_1 + h_2 - q_2, c = q_2 - h_1, d = q_2 - h_2\); 
\end{enumerate}
and vertices are labeled as  in Fig.\ref{fig:labelling}. 

\end{proposizione}

\begin{dimostrazione}
By cutting the genus 2 surface $F_{\{0,1\}}$, represented in Fig.\ref{fig:4610111_nedela}, along $C_0$ and $C_1$ we get a diagram on the sphere with four holes   with \(C_2\) as the equator and the two halves of \(C_1\) and \(C_0\) on opposite emispheres (see  Fig.\ref{fig:diag_0}). We give a standard clockwise orientation to \(C_i^-\), a standard counterclockwise orientation to \(C_i^+\), and an orientation from right to left to \(C_2\). \\

\begin{figure}[h!]
\labellist
\small\hair 2pt
\pinlabel \(C_1^-\) at 220 125
\pinlabel \(C_1^+\) at 220 470
\pinlabel \(C_0^-\) at 530 125
\pinlabel \(C_0^+\) at 530 470
\pinlabel \(C_2\) at 730 350
\endlabellist
\begin{center}
\includegraphics[width = .6\textwidth]{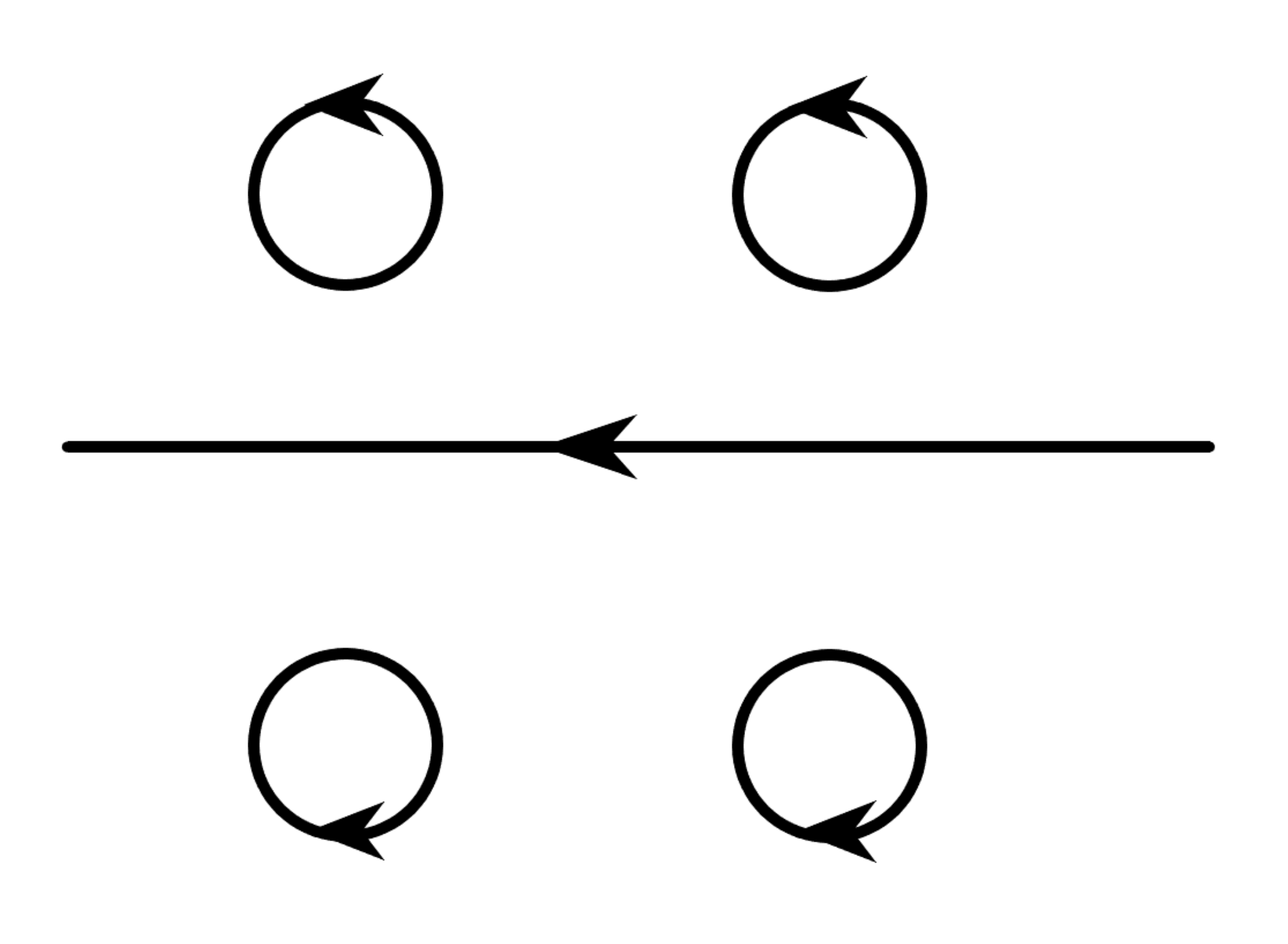}
\caption{The standard orientation.}
\label{fig:diag_0}
\end{center}
\end{figure}
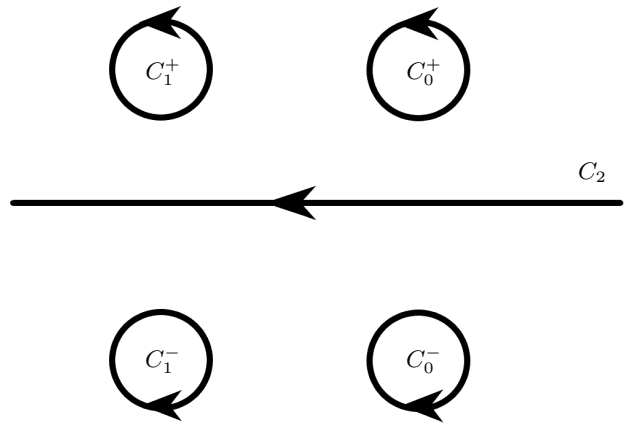

According to the definition of the involutions \(\tau_c\), for \(c = 0, 1, 2, 3\), the labelling of the vertices of \(C_0^\pm, C_1^\pm\) is as follows (see Fig.\ref{fig:labelling}).
In \(C_0^+\) (resp. \(C_0^-\)), the first vertex adjacent to a vertex of \(C_1^+\) (resp. \(C_1^-\)) with respect to the orientation is labelled `\(0\)' (resp. `\(q_0\)'). This vertex is connected to the vertex of \(C_1^+\) (resp. \(C_1^-\)) labelled `\(h_0 + h_1 - 1\)' (resp. `\(q_1 - 1\)'). \\

\begin{figure}[h!]
\labellist
\small\hair 2pt
\pinlabel \(C_1^+\) at 240 450
\pinlabel \(C_0^+\) at 580 450
\pinlabel \(C_1^-\) at 240 150
\pinlabel \(C_0^-\) at 580 150
\pinlabel 0 at 530 495
\pinlabel {\(h_0 + h_1 - 1\)} at 347 495
\pinlabel \(q_0\) at 520 110
\pinlabel \(q_1-1\) at 310 110
\endlabellist
\centering
\includegraphics[width = .6\textwidth]{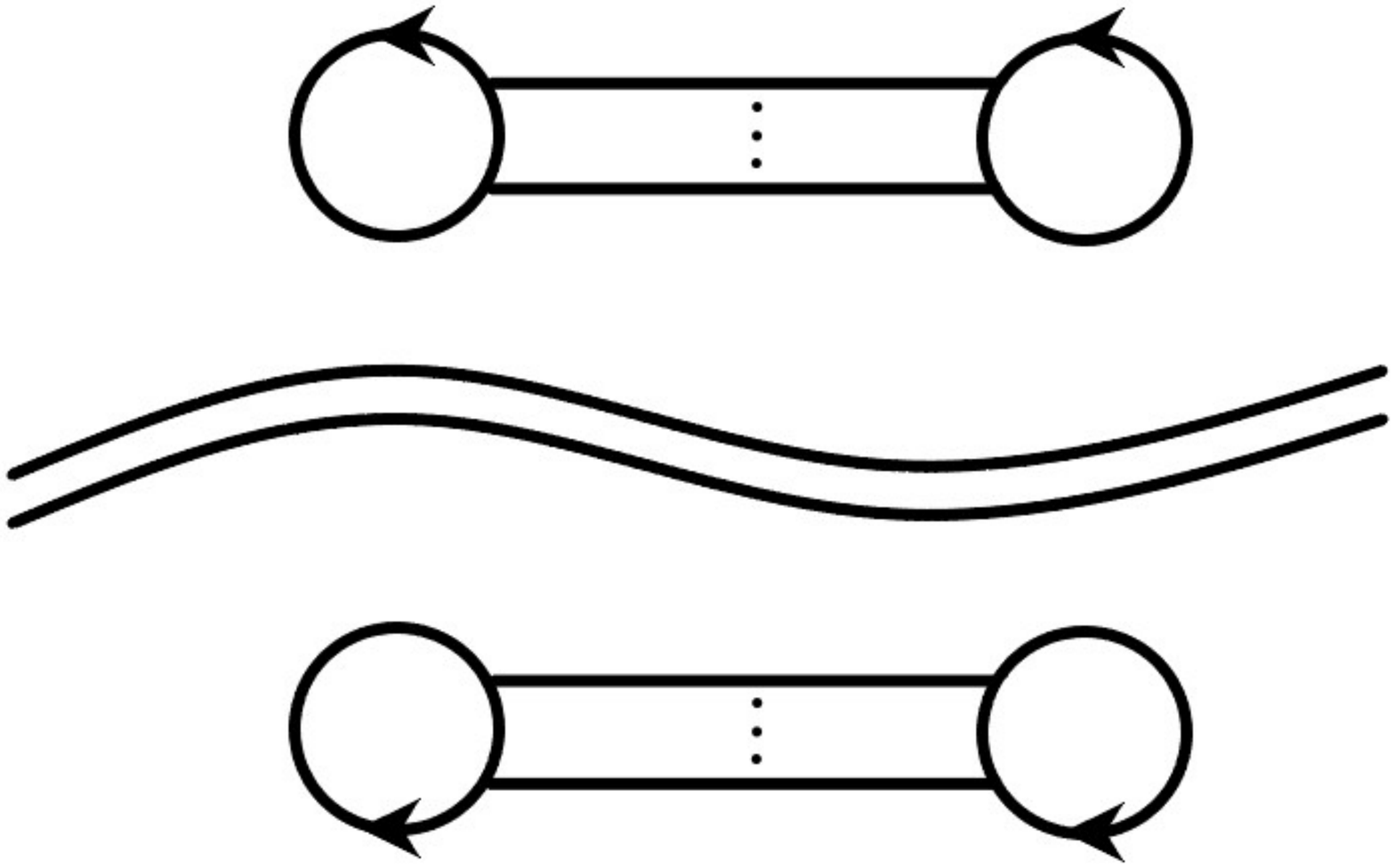}
\caption{The labelling of \(C_0^\pm, \> C_1^\pm\).}
\label{fig:labelling}
\end{figure}

Following the definition of \(\tau_2\) and \(\tau_3\), we observe that (see Fig.\ref{fig:strands_0_a}):
\begin{itemize}
    \item exactly \(h_0\) arcs go from \(C_0^\pm\) to \(C_1^\pm\) and so \(a = h_0\); 
    \item exactly \(h_2\) arcs go from \(C_0^\pm\) to \(C_2\);
    \item exactly \(h_1\) arcs go from \(C_1^\pm\) to \(C_2\).
\end{itemize} 
The involution \(\tau_3\) on the vertices of \(C_2\) differs from \(\tau_2\) by a ``circular shifting'' of \(q_2\) steps 
in the positive direction according to the fixed orientation (as depicted in Fig.\ref{fig:strands_1}). \\
When we remove \(C_2\) in order to obtain the open rich Heegaard diagram, all we have to do is to merge the connecting arcs, obtaining: 
\begin{itemize}
    \item an arc from \(C_0^+\) to \(C_0^-\) if both arcs come from the \(h_2\) arcs; 
    \item an arc from \(C_1^+\) to \(C_1^-\) if both arcs come from the \(h_1\) arcs; 
    \item an arc from \(C_1^+\) to \(C_0^-\) if the upper arc comes from the \(h_1\) arcs and the lower one comes from the \(h_2\) arcs; \item an arc from \(C_0^+\) to \(C_1^-\) if the upper arc comes from the \(h_2\) arcs and the lower one comes from the \(h_1\) arcs;
\end{itemize}
The two last kind of arcs are in the same quantity: indeed, when we shift an \(h_2\) lower arc so that it is connected to an upper \(h_1\) one, we obtain an arc from \(C_1^+\) to \(C_0^-\) and at the same time an \(h_1\) lower arc is connected to an upper \(h_2\) one. \\ 

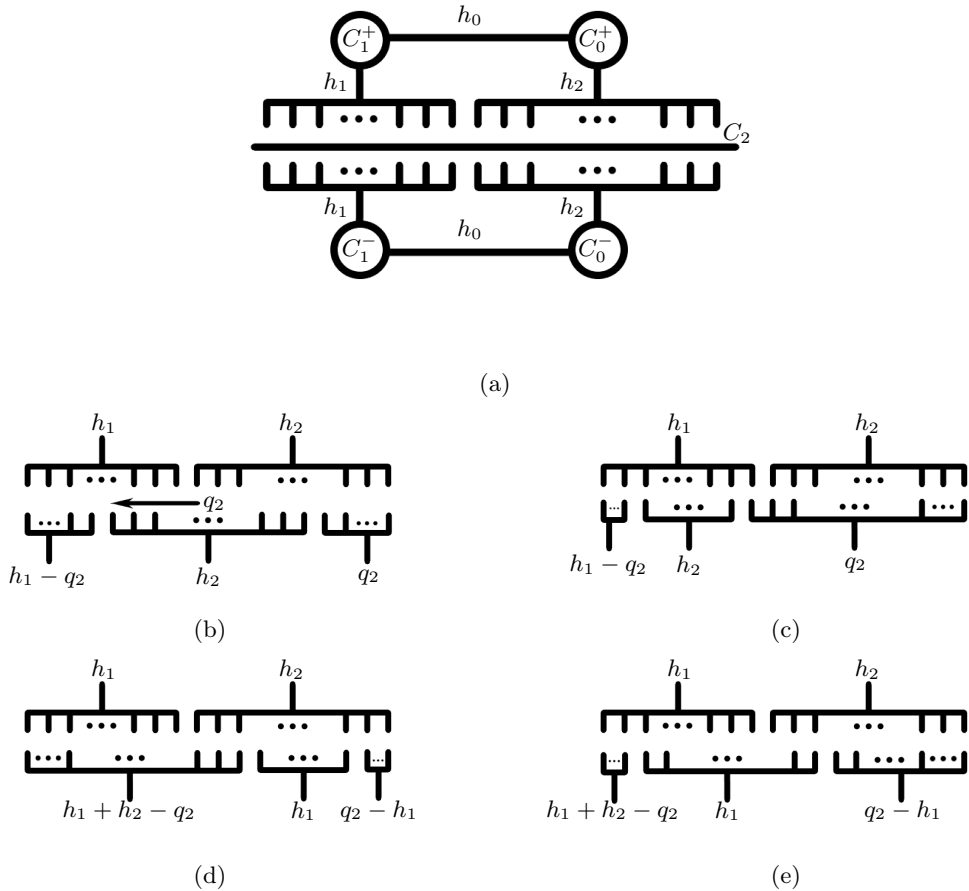
\begin{figure}[h!]
\begin{subfigure}{\textwidth}
\labellist
\small\hair 2pt
\pinlabel \(h_0\) at 350 170
\pinlabel \(h_0\) at 350 460
\pinlabel \(h_1\) at 170 200
\pinlabel \(h_1\) at 170 370
\pinlabel \(h_2\) at 488 200
\pinlabel \(h_2\) at 488 370
\pinlabel \(C_1^-\) at 202 144
\pinlabel \(C_1^+\) at 202 428
\pinlabel \(C_0^-\) at 520 144
\pinlabel \(C_0^+\) at 520 428
\pinlabel \(C_2\) at 710 300
\endlabellist
\centering
\includegraphics[width = .5\textwidth]{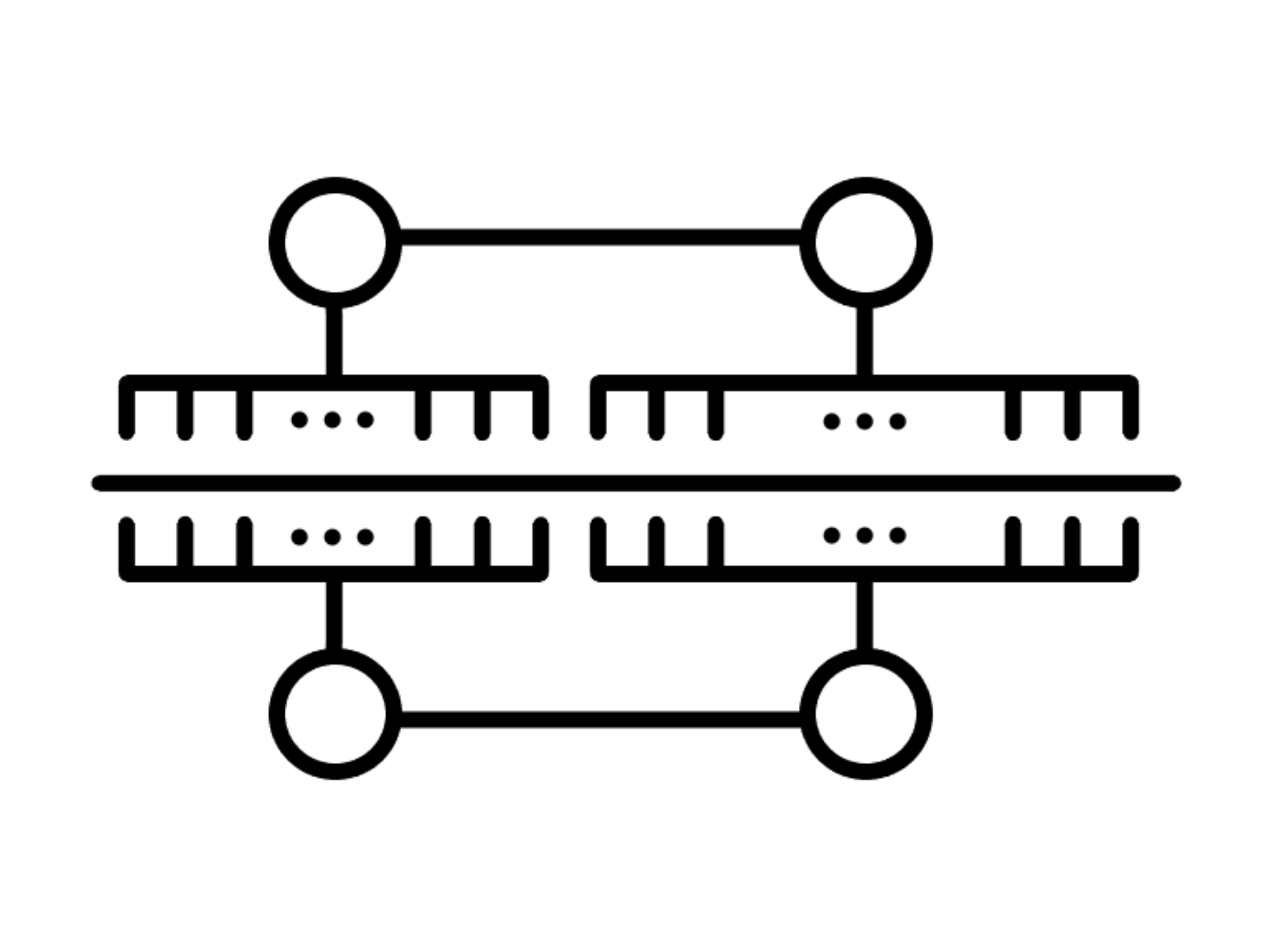}
\caption{}
\label{fig:strands_0_a}
\end{subfigure}

\begin{subfigure}{0.5\textwidth}
\labellist
\small\hair 2pt
\pinlabel \(q_2\) at 390 163
\pinlabel \(h_2\) at 520 300
\pinlabel \(h_2\) at 380 40
\pinlabel \(h_1\) at 205 300
\pinlabel \(h_1-q_2\) at 110 40
\pinlabel \(q_2\) at 650 40
\endlabellist
\centering
\includegraphics[width = .8\textwidth]{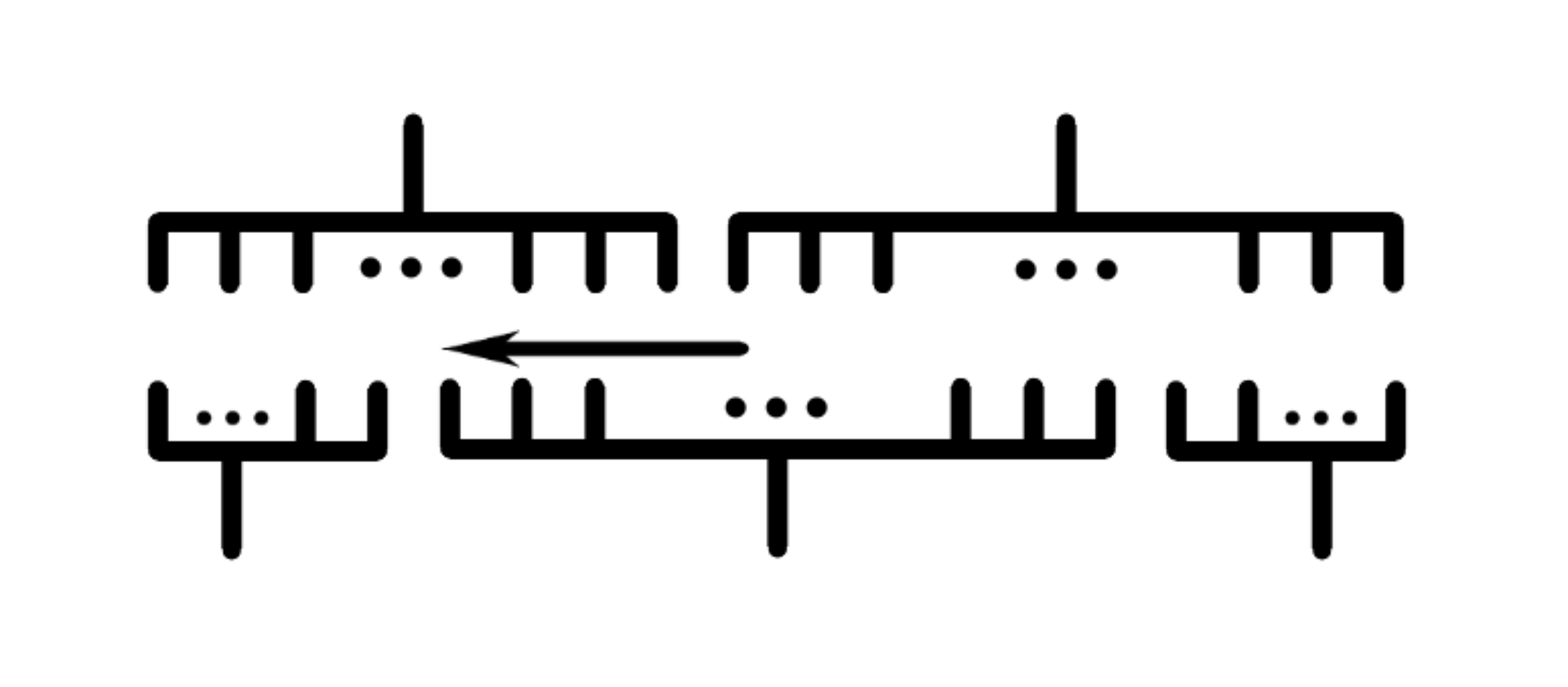}
\caption{}
\label{fig:strands_1}
\end{subfigure}
\begin{subfigure}{0.5\textwidth}
\labellist
\small\hair 2pt
\pinlabel \(h_2\) at 520 300
\pinlabel \(h_2\) at 220 60
\pinlabel \(h_1\) at 205 300
\pinlabel \(h_1-q_2\) at 85 60
\pinlabel \(q_2\) at 500 60
\endlabellist
\centering
\includegraphics[width = .8\textwidth]{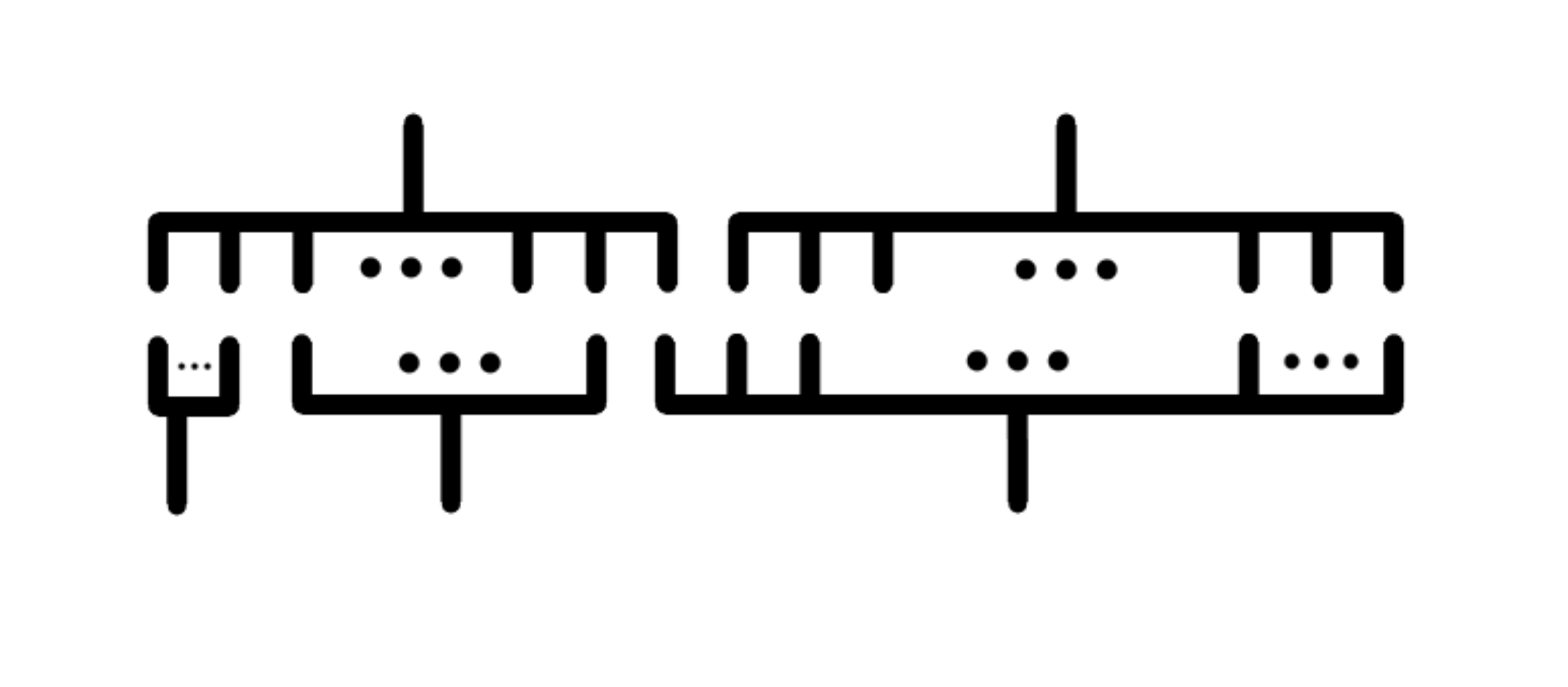}
\caption{}
\label{fig:strands_2}
\end{subfigure}

\begin{subfigure}{0.5\textwidth}
\labellist
\small\hair 2pt
\pinlabel \(h_2\) at 520 300
\pinlabel \(h_1\) at 205 300
\pinlabel \(h_1\) at 540 60
\pinlabel \(h_1+h_2-q_2\) at 247 60
\pinlabel \(q_2-h_1\) at 667 60
\endlabellist
\centering
\includegraphics[width = .8\textwidth]{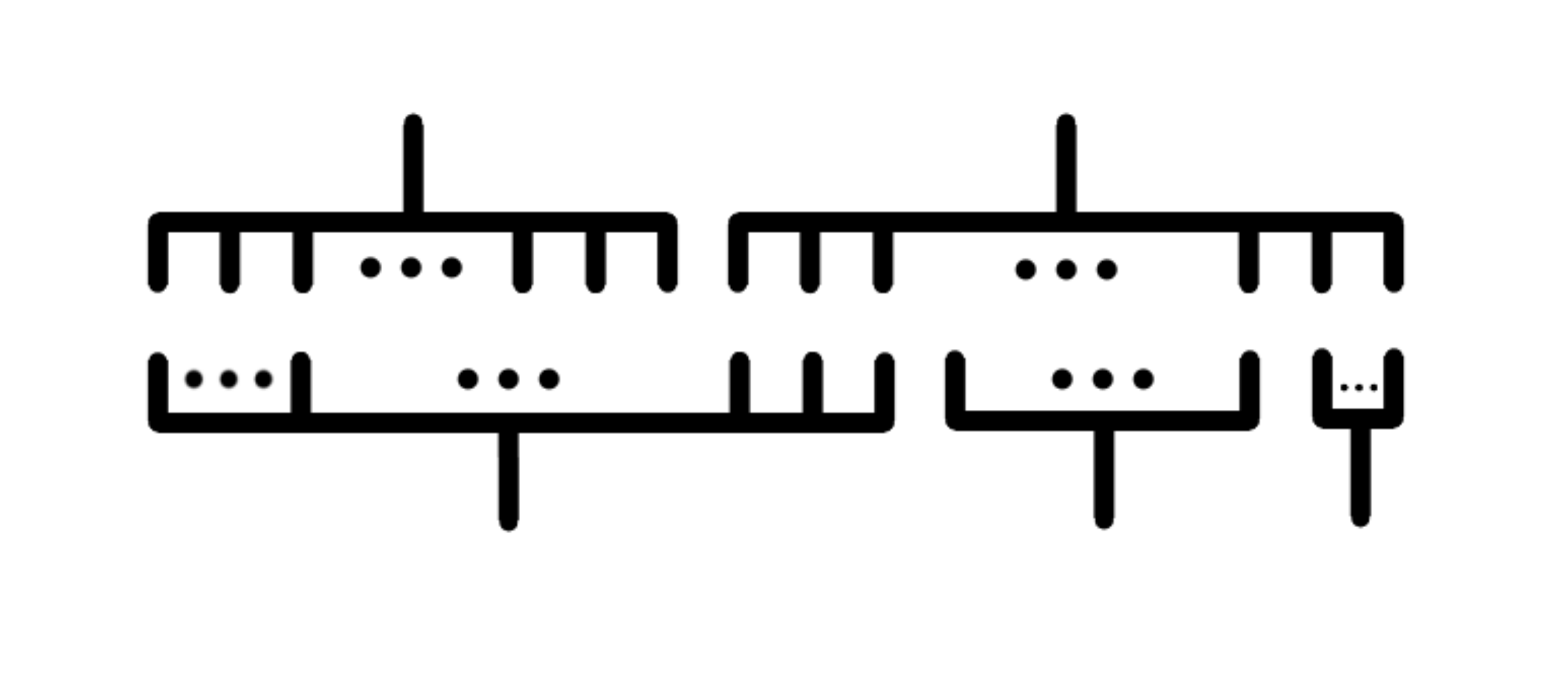}
\caption{}
\label{fig:strands_3}
\end{subfigure}
\begin{subfigure}{0.5\textwidth}
\labellist
\small\hair 2pt
\pinlabel \(h_2\) at 520 300
\pinlabel \(h_1\) at 205 300
\pinlabel \(h_1+h_2-q_2\) at 92 60
\pinlabel \(h_1\) at 285 60
\pinlabel \(q_2-h_1\) at 580 60
\endlabellist
\centering
\includegraphics[width = .8\textwidth]{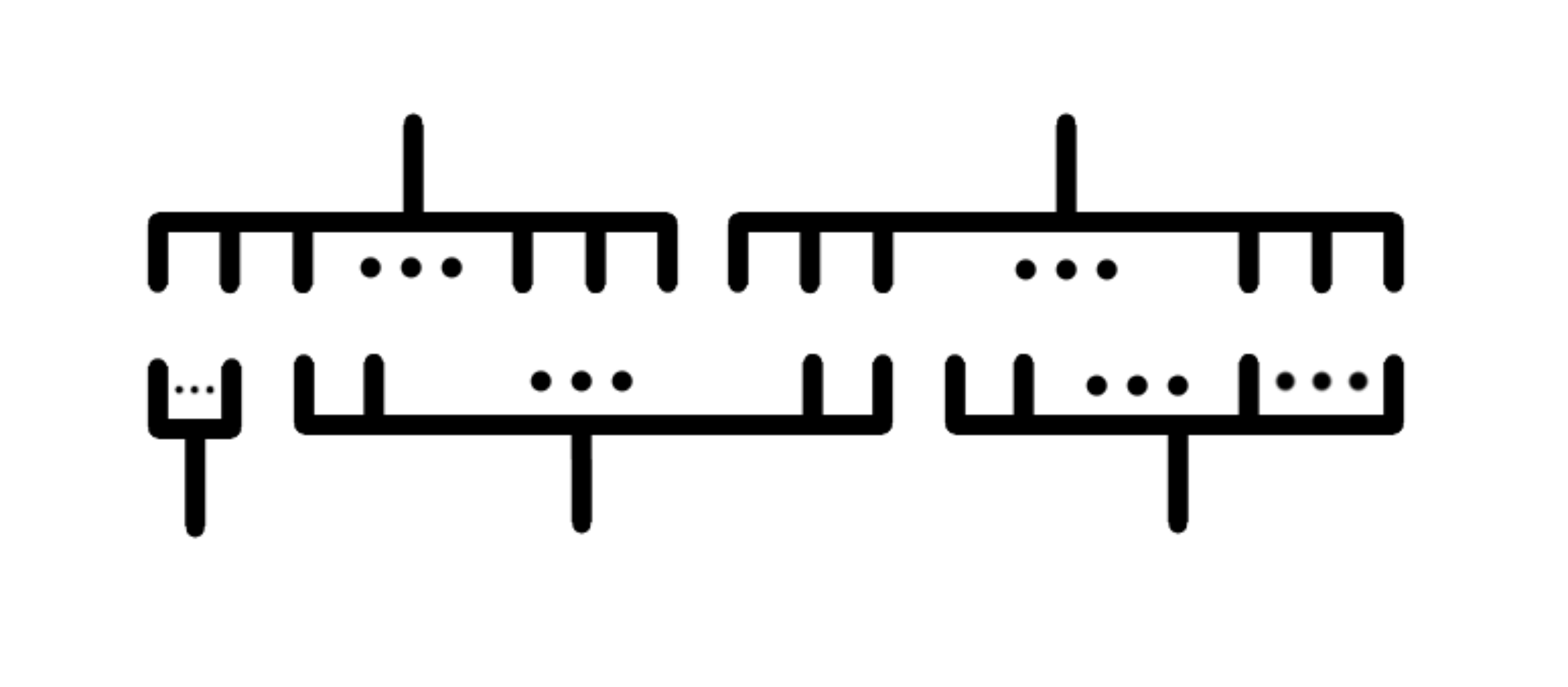}
\caption{}
\label{fig:strands_4}
\end{subfigure}
\caption{The starting position and the four possible arcs positions.}
\end{figure}

\begin{figure}[h!]
\begin{subfigure}{\textwidth}
\labellist
\small\hair 2pt
\pinlabel \(C_1^-\) at 83 135
\pinlabel \(C_1^+\) at 83 350
\pinlabel \(C_0^-\) at 315 135
\pinlabel \(C_0^+\) at 310 350
\pinlabel \(a\) at 200 365
\pinlabel \(a\) at 200 150
\pinlabel \(b\) at 160 240
\pinlabel \(b\) at 460 100
\pinlabel \(c\) at 330 240
\pinlabel \(d\) at 60 240
\endlabellist
\centering
\includegraphics[width = .35\textwidth]{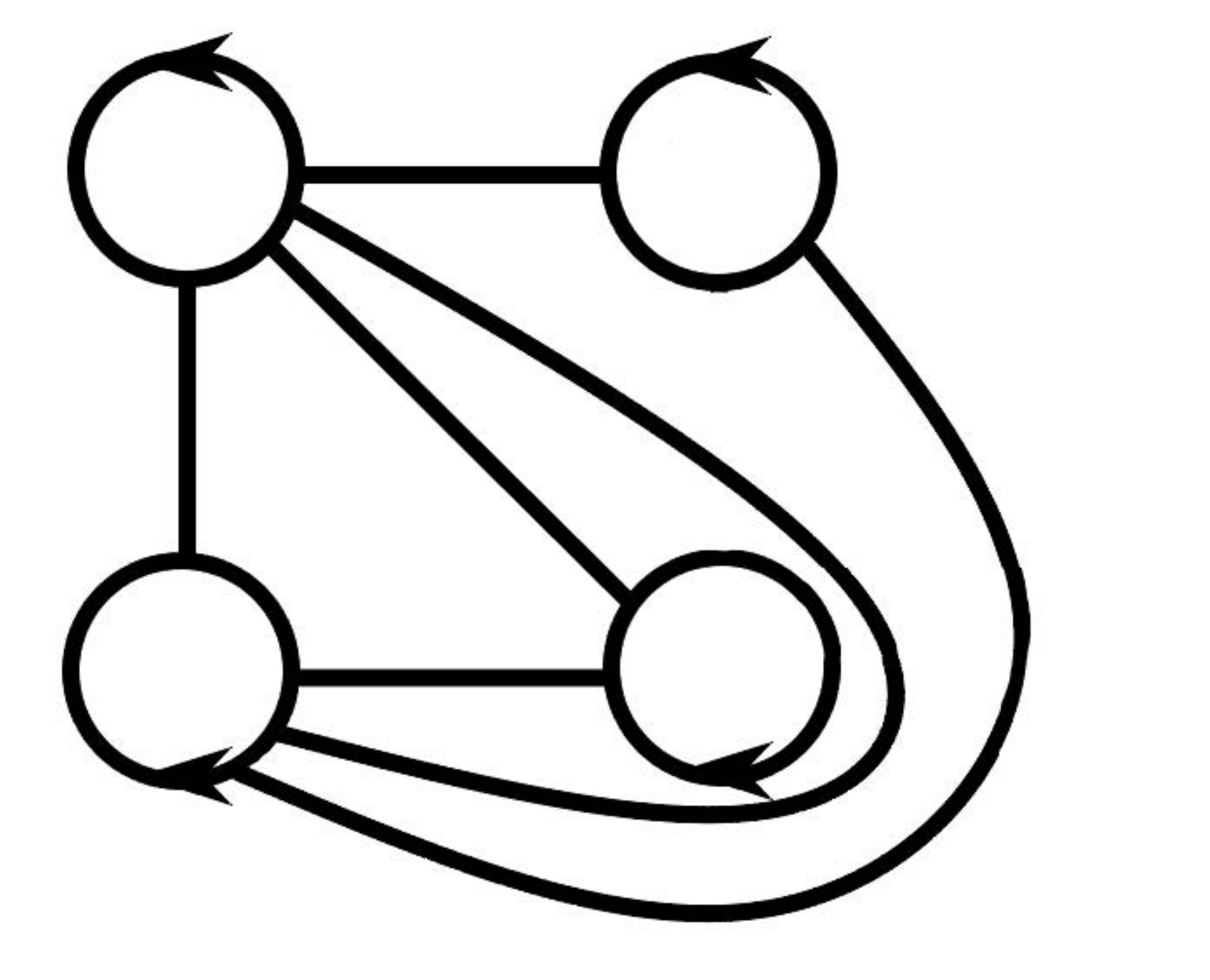}
\caption{}
\label{fig:heeg_2}
\end{subfigure}

\begin{subfigure}{0.5\textwidth}
\labellist
\small\hair 2pt
\pinlabel \(C_1^-\) at 215 242
\pinlabel \(C_1^+\) at 215 495
\pinlabel \(C_0^-\) at 480 242
\pinlabel \(C_0^+\) at 480 495
\pinlabel \(a\) at 345 520
\pinlabel \(a\) at 345 265
\pinlabel \(b\) at 345 400
\pinlabel \(b\) at 585 175
\pinlabel \(d\) at 495 375
\pinlabel \(c\) at 85 300
\endlabellist
\centering
\includegraphics[width = .8\textwidth]{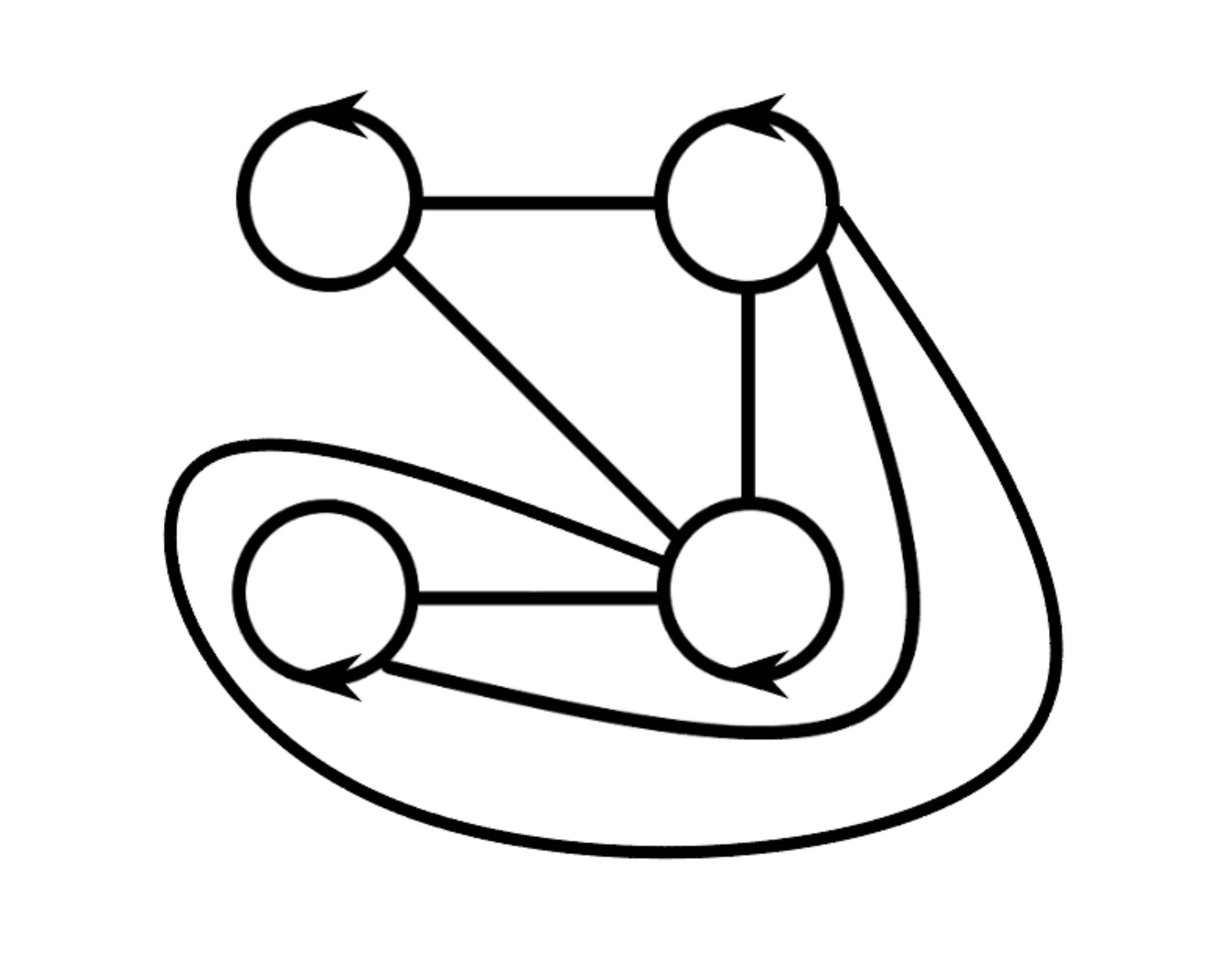}
\caption{}
\label{fig:heeg_3}
\end{subfigure}
\begin{subfigure}{0.5\textwidth}
\labellist
\small\hair 2pt
\pinlabel \(C_1^-\) at 125 165
\pinlabel \(C_1^+\) at 125 345
\pinlabel \(C_0^-\) at 315 165
\pinlabel \(C_0^+\) at 315 345
\pinlabel \(a\) at 220 355
\pinlabel \(a\) at 220 175
\pinlabel \(b\) at 190 255
\pinlabel \(b\) at 420 150
\pinlabel \(c\) at 460 130
\pinlabel \(d\) at 310 270
\endlabellist
\centering
\includegraphics[width = .8\textwidth]{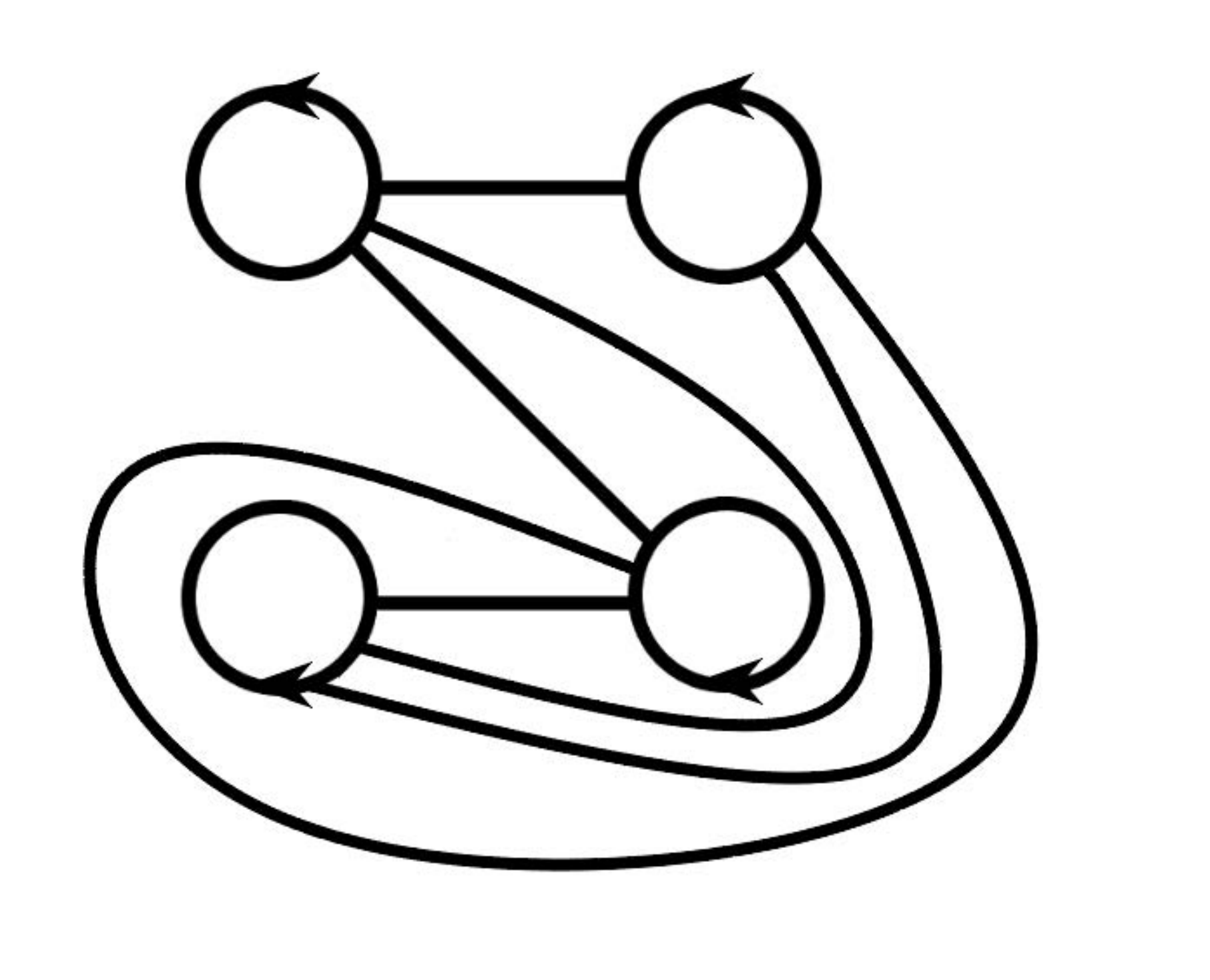}
\caption{}
\label{fig:heeg_4}
\end{subfigure}
\caption{The three unrefined configurations.}
\end{figure}
\vspace{-8pt}

Now, depending on the values of \(h_1, h_2\) and \(q_2\) we can have four different cases: 
\begin{enumerate}
    \item if \(q_2 < h_1, h_2\), we are in the situation depicted in Fig.\ref{fig:strands_1}: here \(q_2\) arcs of the \(h_1\) lower ones are moved to the right. So we obtain a graph like the one in Fig.\ref{fig:heeg_1}. We have exactly \(q_2\) arcs from upper \(h_1\) to lower \(h_2\), so \(b = q_2\). Since the remaining lower arcs from \(h_2\) go into upper \(h_2\) ones, we have \(c = h_2 - q_2\). Lastly, in order to balance the accounts for the lower \(h_1\) arcs, we get \(d = h_1 - q_2\). 
    \item if \(h_2 \leq q_2 \leq h_1\), we are in the situation depicted in Fig.\ref{fig:strands_2}: as before, all the \(h_2\) lower arcs go into upper \(h_1\) ones and \(q_2\) arcs of the \(h_1\) lower ones are moved to the right. So we obtain a graph like the one in Fig.\ref{fig:heeg_2}. We have all the lower \(h_2\) arcs going to upper \(h_1\) ones, so \(b = h_2\). To determine \(c\) and \(d\) we need to count how many arcs are there on the right and on the left. Keeping track of \(q_2\) and of the hypothesis \(h_2 \leq q_2 \leq h_1\), on the right there are \(q_2 - h_2 = c\) arcs while on the left there are \(h_1 - q_2 = d\) arcs. \\
    Using a homeomorphism of the sphere invariant on the set of vertices, we can transform the diagram from of Fig.\ref{fig:heeg_2} into the one of Fig.\ref{fig:heeg_2_giusto}, passing the right \(b\) arcs onto the left. 
    \item if \(h_1 < q_2 < h_2\), we are in the situation depicted in Fig.\ref{fig:strands_3}: here all the \(h_1\) lower arcs and \(q_2 - h_1\) lower \(h_2\) arcs are moved to the right. So, we obtain a graph like the one in Fig.\ref{fig:heeg_3}. We have all the \(h_1\) lower arcs going to upper \(h_2\) ones, so \(b = h_1\). As before, to determine \(c\) and \(d\) we need to count how many arcs are there on the right and on the left. Keeping track of \(q_2\) and of the hypothesis \(h_1 < q_2 < h_2\), on the right there are \(q_2 - h_1 = c\) arcs while on the left there are \(h_2 - q_2 = d\) arcs. \\
    Using a homeomorphism of the sphere invariant on the set of vertices, we can transform the diagram of Fig.\ref{fig:heeg_3} into the one of Fig.\ref{fig:heeg_3_giusto}, passing the right \(d\) arcs onto the left and flipping the diagram horizontally and vertically. 
    \item if \(q_2 > h_1, h_2\), we are in the situation depicted in Fig.\ref{fig:strands_4}: as before, all \(h_1\) lower arcs and \(q_2 - h_1\) lower \(h_2\) arcs are moved to the right. So, we obtain a graph like the one in Fig.\ref{fig:heeg_4}, where \(b = h_1 + h_2 - q_2\), \(d = h_1 - (h_1 + h_2 - q_2) = q_2 - h_2\) and \(c = q_2 - h_1\). \\
    Using a homeomorphism of the sphere invariant on the set of vertices, we can transform the diagram of Fig.\ref{fig:heeg_4} into the one of Fig.\ref{fig:heeg_4_intermedio}, passing the right \(c\) arcs onto the left. Then, following the moves of Fig.\ref{fig:heeg_4_mosse} and flipping horizontally, we obtain the one in Fig.\ref{fig:heeg_4_giusto}. 
\end{enumerate}

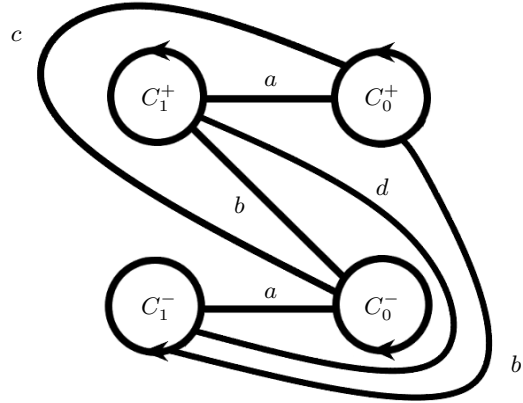
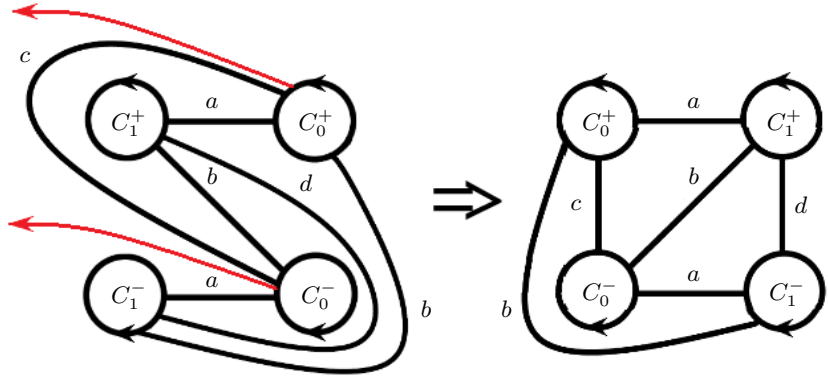
\begin{figure}[h!]
\begin{subfigure}{\textwidth}
\labellist
\small\hair 2pt
\pinlabel \(C_1^-\) at 260 195
\pinlabel \(C_1^+\) at 260 430
\pinlabel \(C_0^-\) at 510 195
\pinlabel \(C_0^+\) at 510 430
\pinlabel \(a\) at 385 450
\pinlabel \(a\) at 385 210
\pinlabel \(b\) at 350 310
\pinlabel \(c\) at 100 500
\pinlabel \(d\) at 510 330
\pinlabel \(b\) at 660 130
\endlabellist
\centering
\includegraphics[width = .6\textwidth]{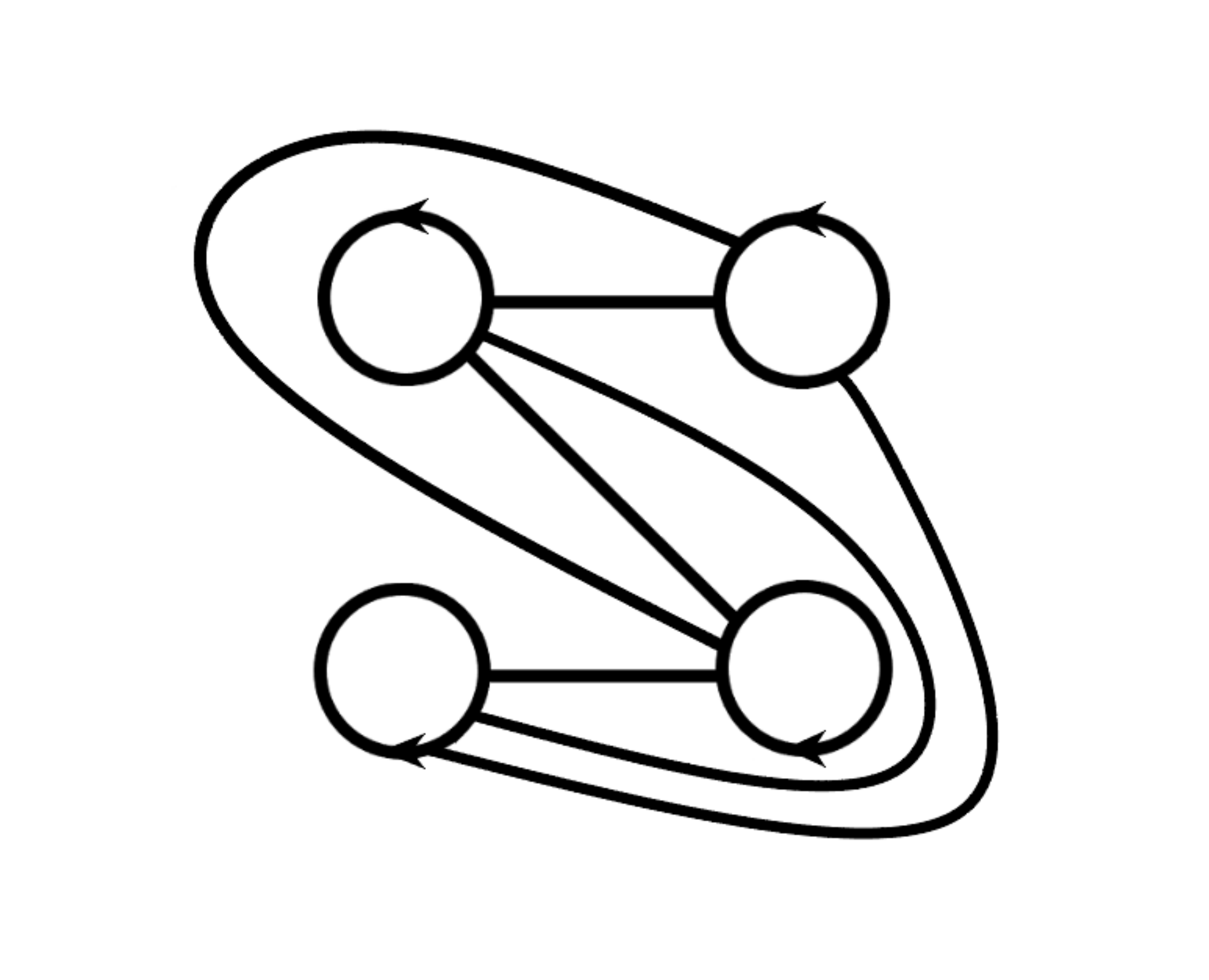}
\caption{}
\label{fig:heeg_4_intermedio}
\end{subfigure}

\begin{subfigure}{\textwidth}
\labellist
\small\hair 2pt
\pinlabel \(C_1^-\) at 178 110
\pinlabel \(C_1^+\) at 178 240
\pinlabel \(C_0^-\) at 320 110
\pinlabel \(C_0^+\) at 318 240
\pinlabel \(a\) at 240 255
\pinlabel \(a\) at 240 120
\pinlabel \(b\) at 240 200
\pinlabel \(d\) at 310 195
\pinlabel \(b\) at 400 100
\pinlabel \(c\) at 100 290
\pinlabel \(C_1^-\) at 670 112
\pinlabel \(C_1^+\) at 668 240
\pinlabel \(C_0^-\) at 530 112
\pinlabel \(C_0^+\) at 530 240
\pinlabel \(a\) at 600 125
\pinlabel \(a\) at 600 255
\pinlabel \(b\) at 460 100
\pinlabel \(b\) at 600 200
\pinlabel \(c\) at 512 178
\pinlabel \(d\) at 680 178
\endlabellist
\centering
\includegraphics[width = 0.9\textwidth]{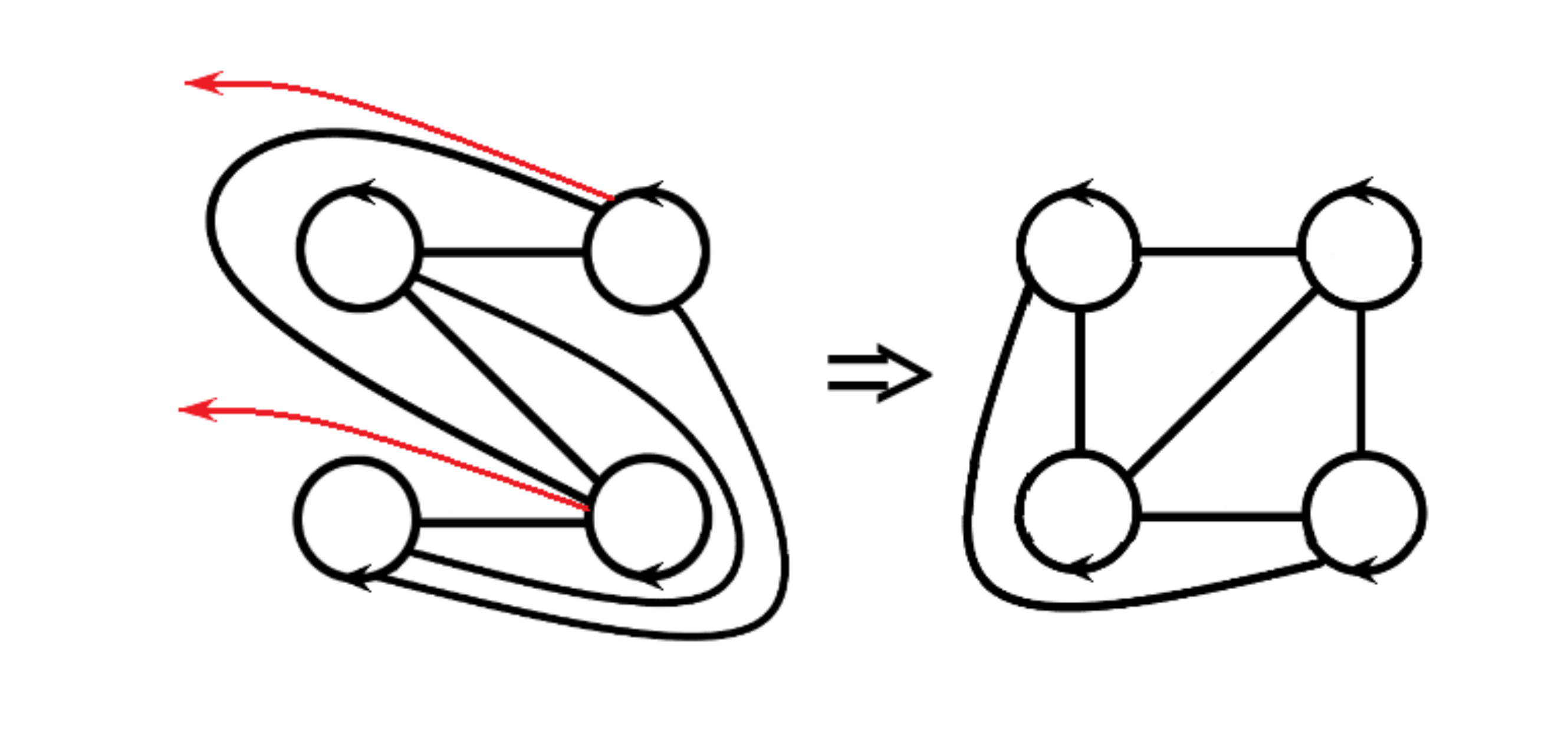}
\caption{}
\label{fig:heeg_4_mosse}
\end{subfigure}
\caption{Dealing with the fourth case.}
\end{figure}
This ends the proof. \\
\qed
\end{dimostrazione}
In Fig.\ref{fig:4610111} it is depicted the rich Heegaard diagram obtained from the \(6\)-tuple \(f = (4, 6, 10, 1, 1, 1)\); the manifold \(M_f\) is the Poincaré sphere (see \cite{bandieri2010census}). \\

\begin{figure}[h!]
\labellist
\small\hair 2pt
\pinlabel \(C_1^-\) at 100 125
\pinlabel \(C_0^-\) at 430 125
\pinlabel \(C_1^+\) at 100 410
\pinlabel \(C_0^+\) at 430 410
\endlabellist
    \centering
    \includegraphics[width = 0.4\textwidth]{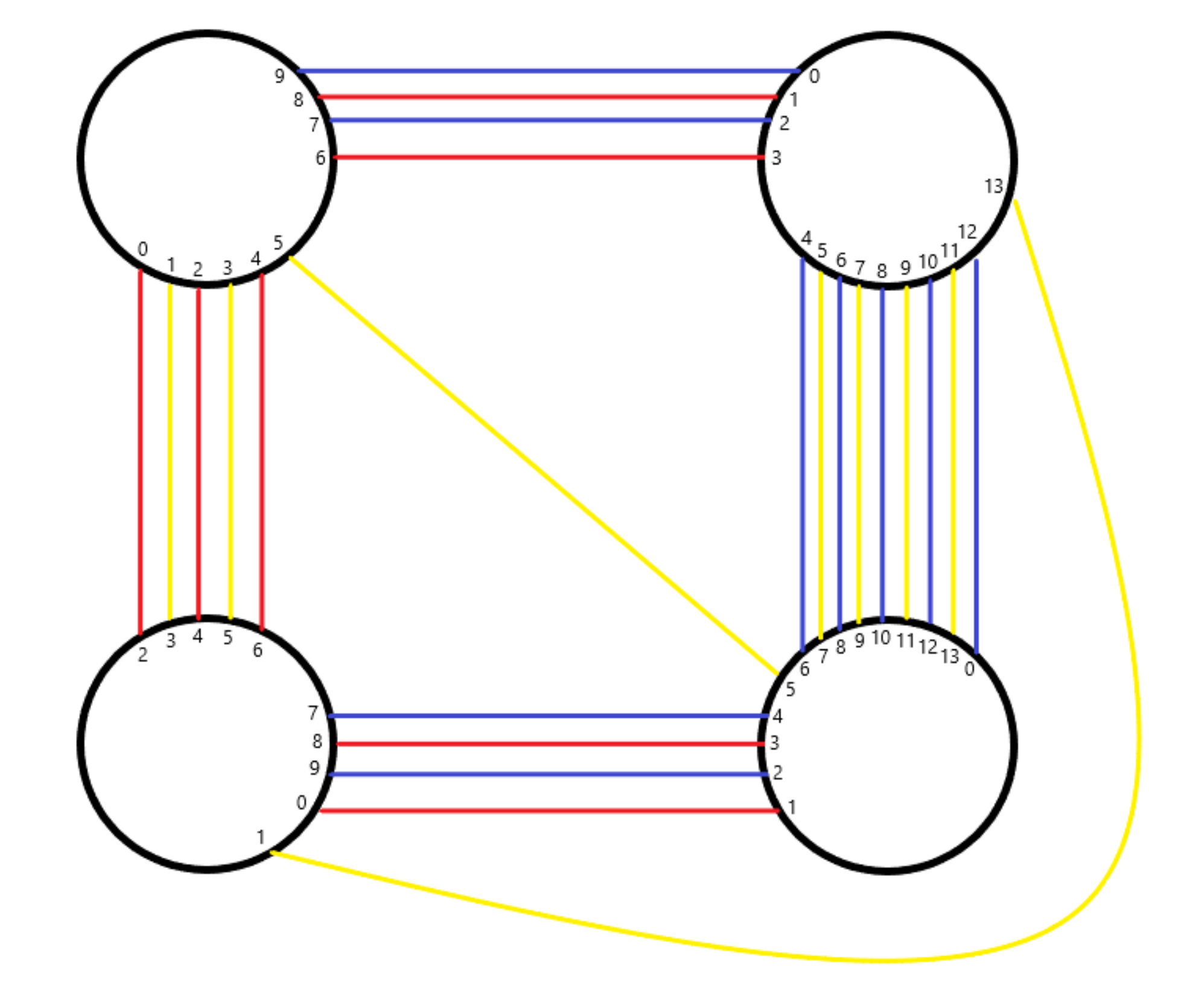}
    \caption{The rich Heegaard diagram obtained from \((4, 6, 10, 1, 1, 1)\).}
    \label{fig:4610111}
\end{figure}
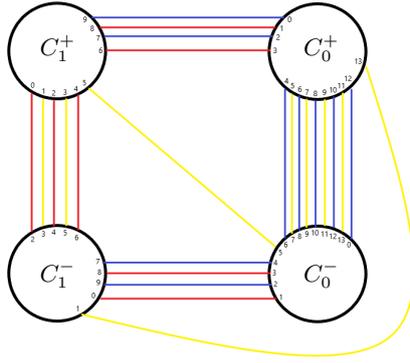
\begin{osservazione}
The normalization of the diagrams that we apply using homeomorphisms of the sphere in all cases different from the first one are not essential but are done in order to obtain a normalized diagram in the sense of \cite[Chapter 5]{fomenko2013algorithmic}. In this article, the authors describe three possible classes of genus 2 open Heegaard diagrams and prove that, up to equivalence, each genus 2 open Heegaard diagram belongs to one of them. The equivalence they consider between diagrams corresponds to a classification of the Heegaard splittings up to homeomorphism. Moreover, they prove that if one is interested only in manifolds (and not in their splittings), the class III is not necessary. Indeed, the diagrams that we obtain in the previous result belong to the classes I and II. \end{osservazione} 

\section{Main result}
\label{main}
In this section, starting from the 6-tuple \(f\), we construct an algorithm that allows to describe, in terms of elements in \(B_{2, 2n}\), the plat slide moves in \(M_f\) associated to the Heegaard diagram obtained from \(R_f\). \\

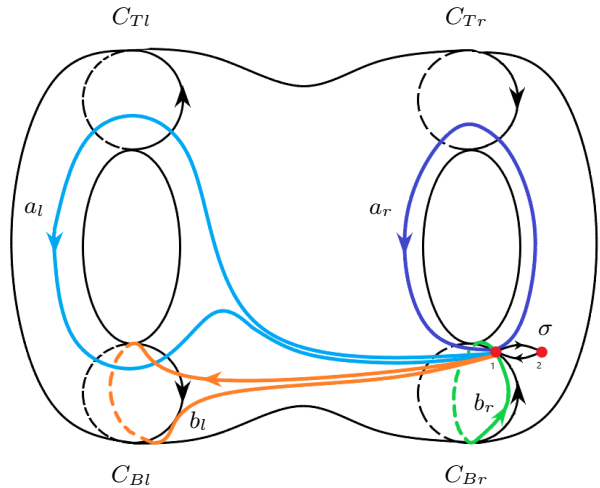
\begin{figure}[h!]
\labellist
\small\hair 2pt
\pinlabel \(C_{Bl}\) at 280 80
\pinlabel \(C_{Br}\) at 555 80
\pinlabel \(C_{Tl}\) at 280 460
\pinlabel \(C_{Tr}\) at 555 460
\pinlabel \(a_l\) at 200 300
\pinlabel \(a_r\) at 485 300
\pinlabel \(b_l\) at 335 125
\pinlabel \(b_r\) at 570 140
\pinlabel \(\sigma\) at 620 200
\endlabellist
    \centering
    \includegraphics[width = 0.9\textwidth]{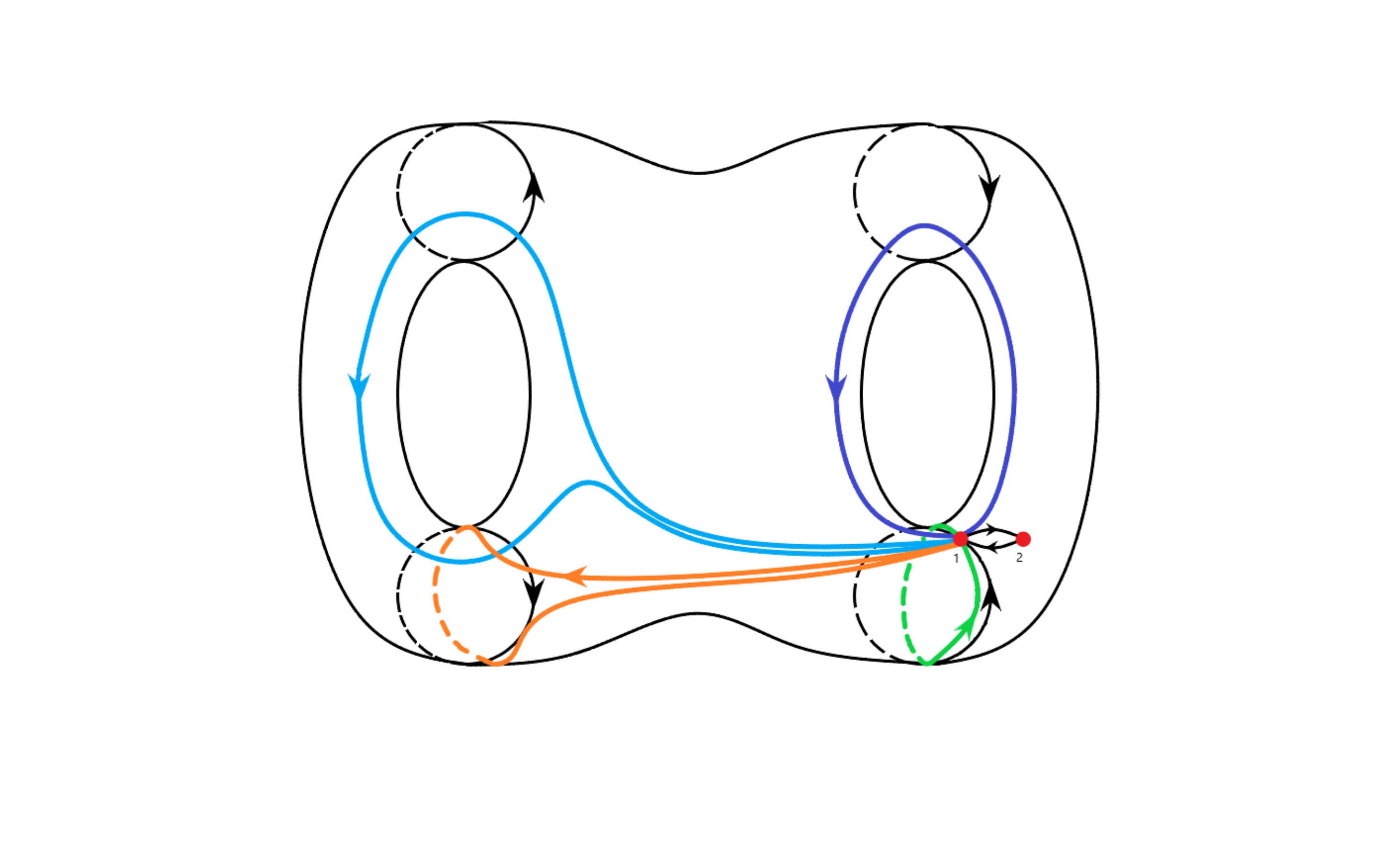}
\caption{The genus 2 surface.}
\label{fig:word_0}
\end{figure}

First of all, we need to fix some notation on the genus 2 surface. Referring to Fig.\ref{fig:word_0}, we denote with $b_l,b_r,a_l,a_r,\sigma$ the standard generators of $B_{2,2}$, where "\(l\)" and "\(r\)" stand for left and right. Moreover, we identify the boundary circles of 4-holed sphere containing $R_f$ with  the meridians $C_{Tr},C_{Tl}, C_{Br}, C_{Bl}$, where "$T$" stands for top and "$B$"  for bottom, so that $C_{Tl}$ is identified with $C_1^+$ in case \ref{prop:dati_1}, \ref{prop:dati_2}, \ref{prop:dati_4} of Proposition \ref{prop:dati} (corresponding to pictures \ref{fig:heeg_1}, \ref{fig:heeg_2_giusto}, \ref{fig:heeg_4_giusto} respectively) and with $C_0^-$ in case \ref{prop:dati_3} (corresponding to picture \ref{fig:heeg_3_giusto}). The four meridians are oriented as depicted in the figure in cases  \ref{prop:dati_1}, \ref{prop:dati_2}, \ref{prop:dati_3}, while we take the opposite orientation for all of them in case \ref{prop:dati_4}.

\begin{figure}
\begin{subfigure}{\textwidth}
    \centering
    \includegraphics[width = 0.4\textwidth]{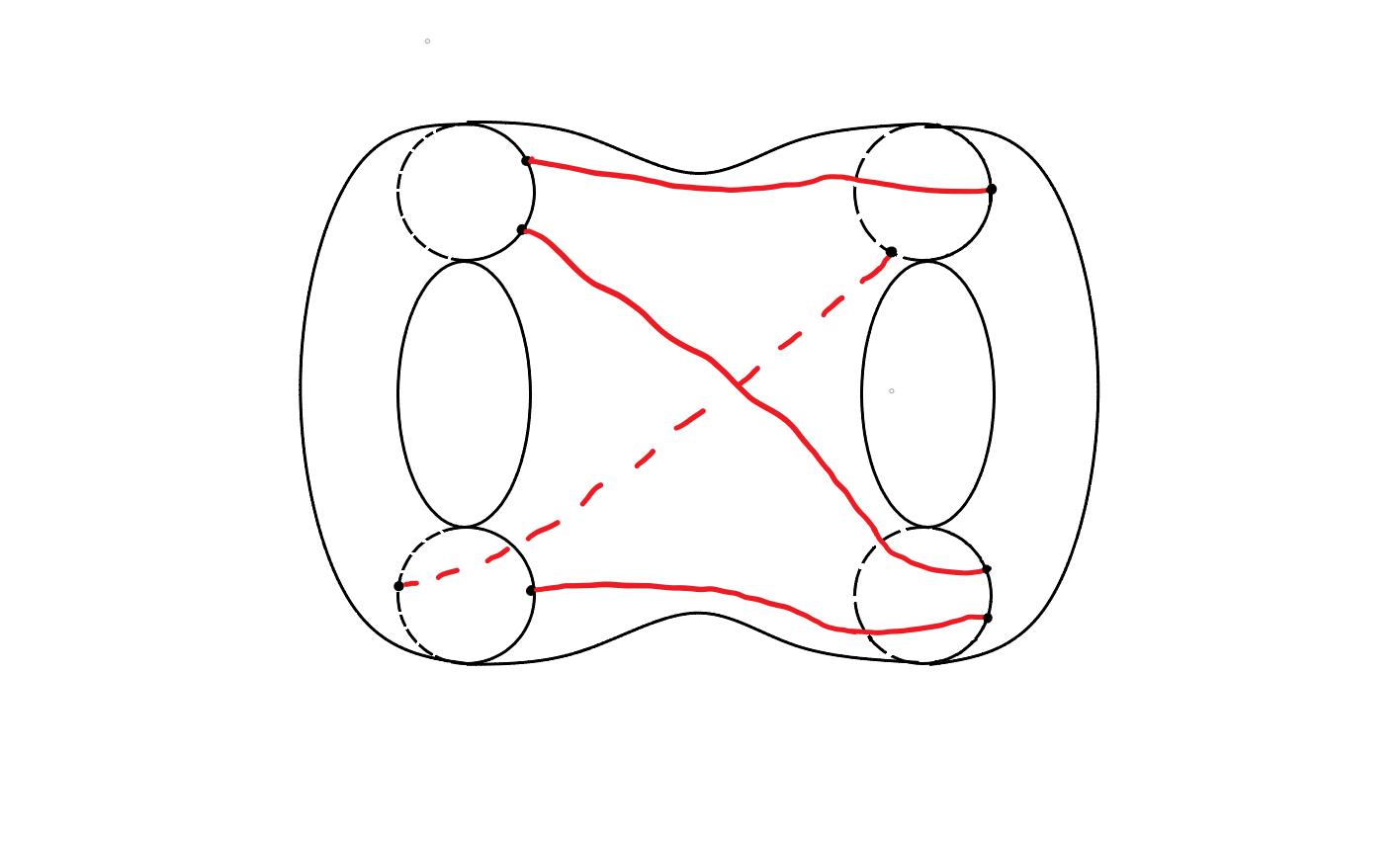}
    \caption{\((3,3,3,2,2,2)\) only red arcs.}
    \label{fig:333222_only}
\end{subfigure}

\begin{subfigure}{\textwidth}
    \centering
    \includegraphics[width = 0.6\textwidth]{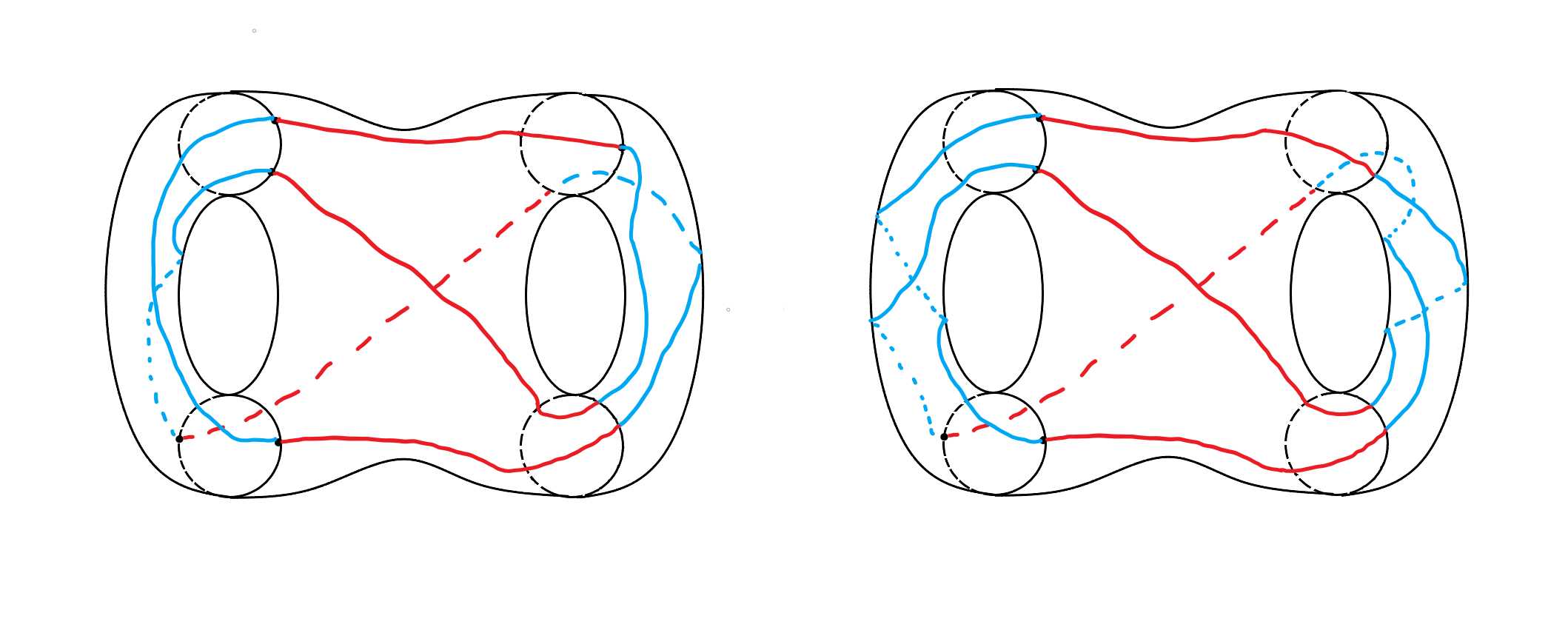}
    \caption{Two ways to connect the arcs of the \((3,3,3,2,2,2)\).}
    \label{fig:333222_1}
\end{subfigure}
\caption{The \((3,3,3,2,2,2)\) case.}
\end{figure}

\begin{figure}
\begin{subfigure}{.33\textwidth}
  \centering
  \includegraphics[width=\textwidth]{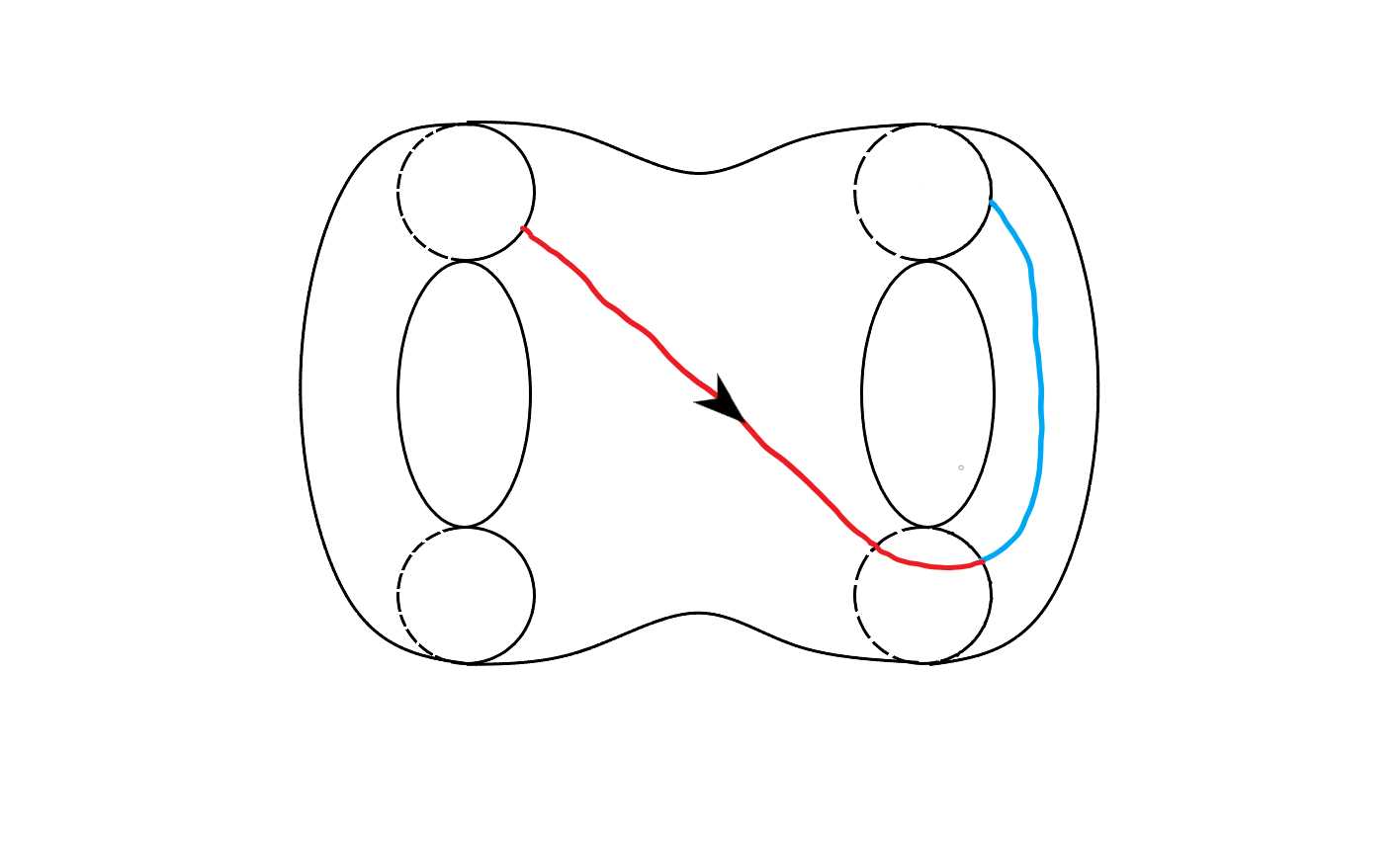}
  \caption{}
\end{subfigure}
\begin{subfigure}{.33\textwidth}
  \centering
  \includegraphics[width=\textwidth]{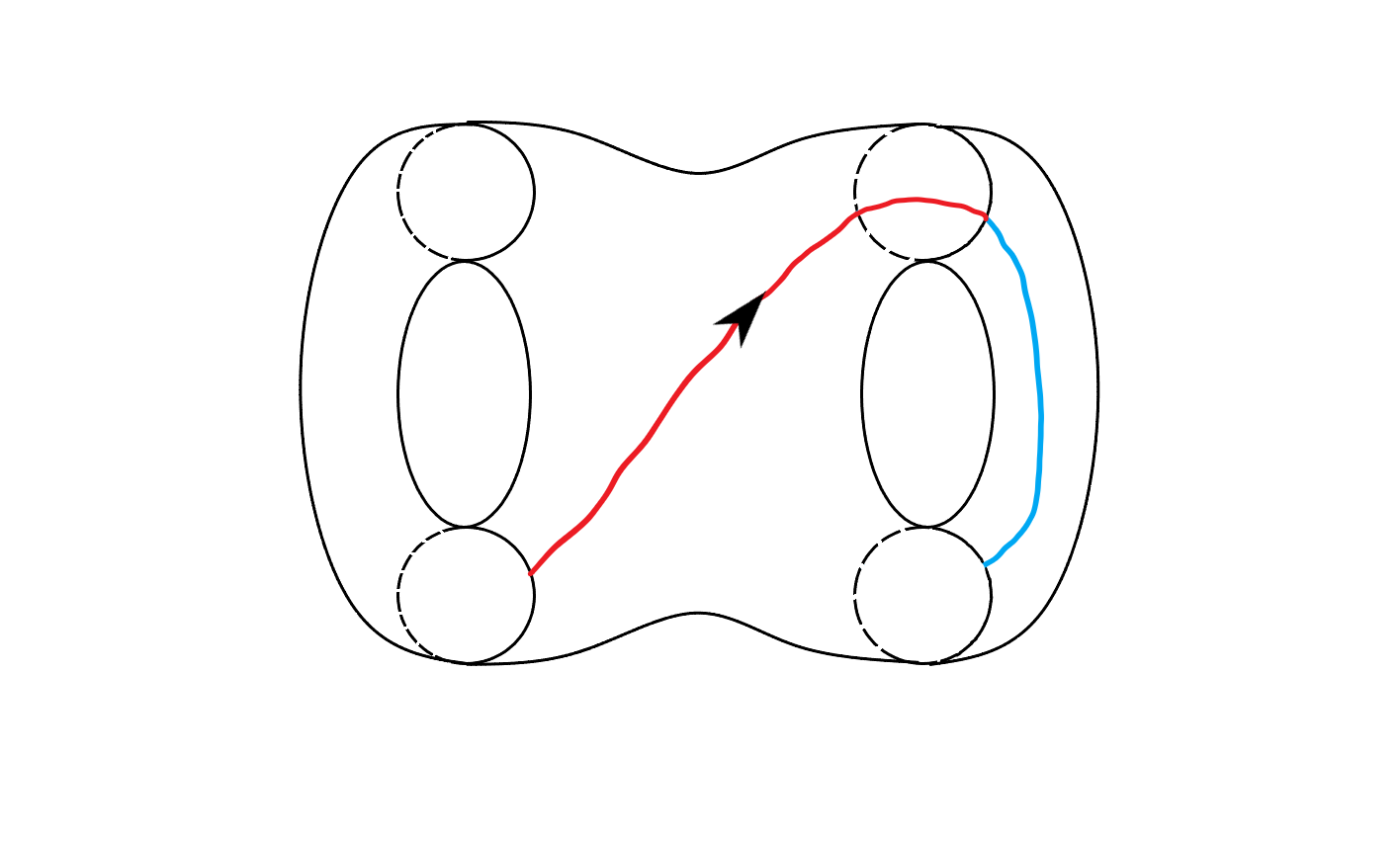}
  \caption{}
\end{subfigure}
\begin{subfigure}{.33\textwidth}
  \centering
  \includegraphics[width=\textwidth]{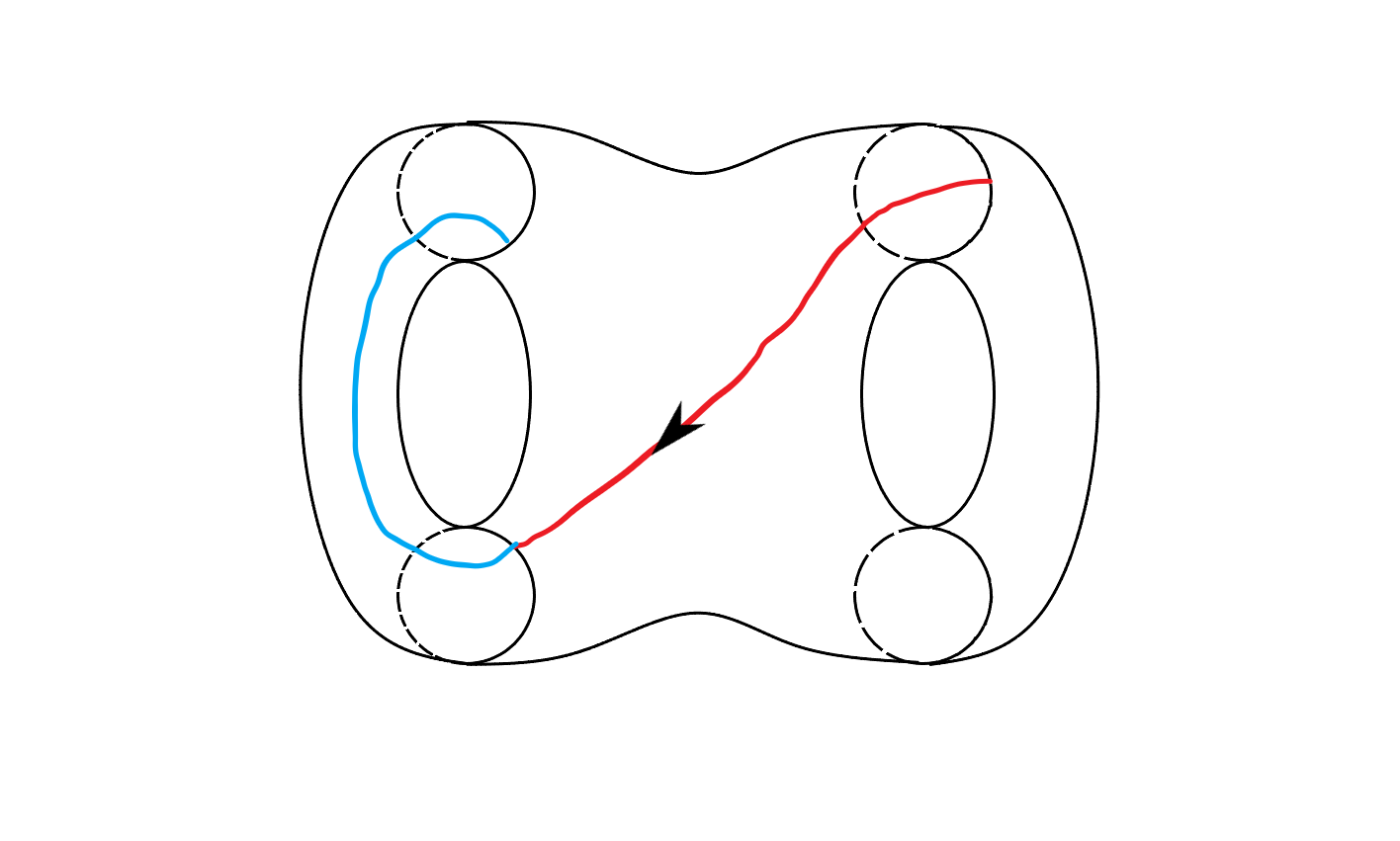}
  \caption{}
\end{subfigure}

\begin{subfigure}{.33\textwidth}
  \centering
  \includegraphics[width=\textwidth]{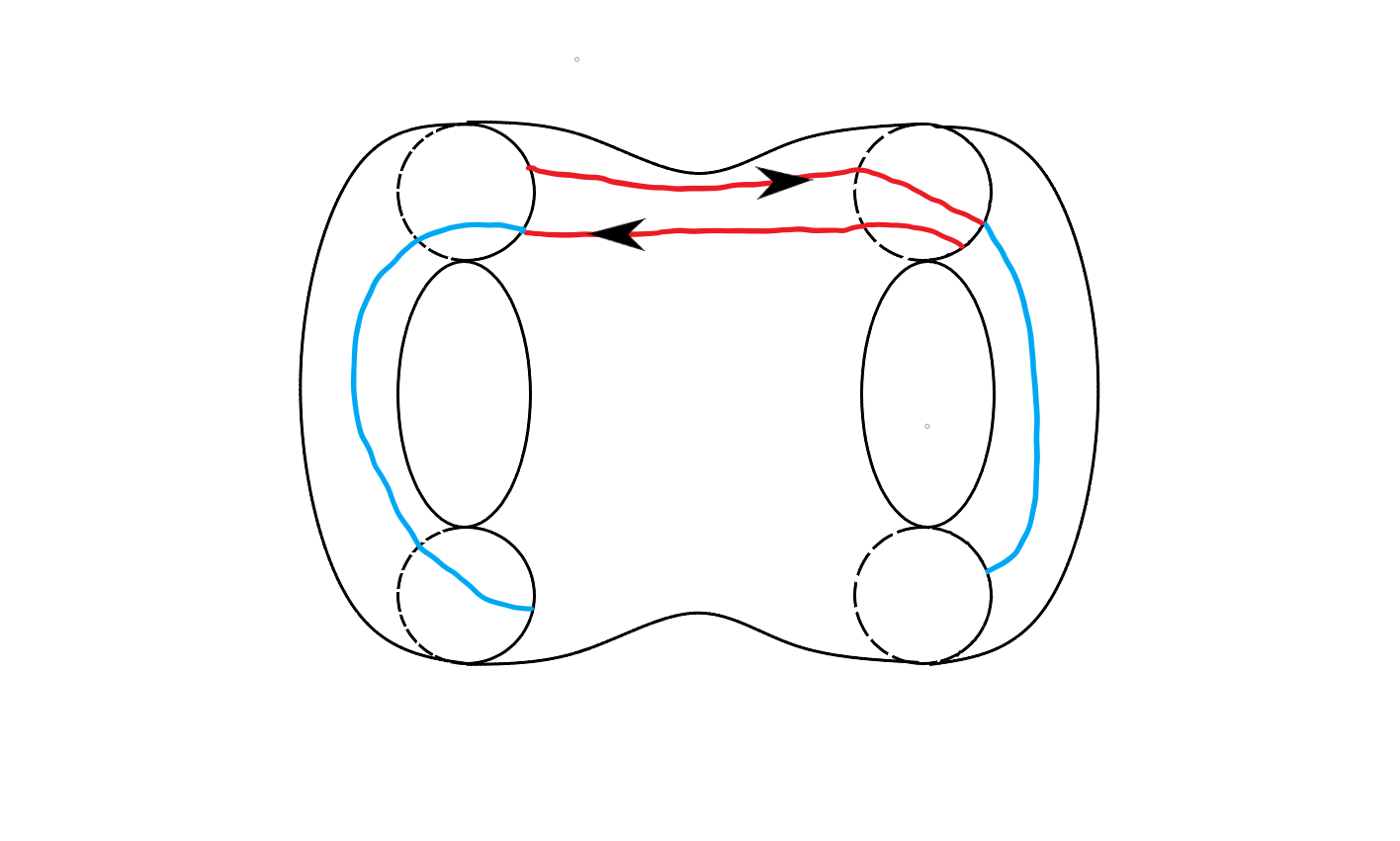}
  \caption{}
\end{subfigure}
\begin{subfigure}{.33\textwidth}
  \centering
  \includegraphics[width=\textwidth]{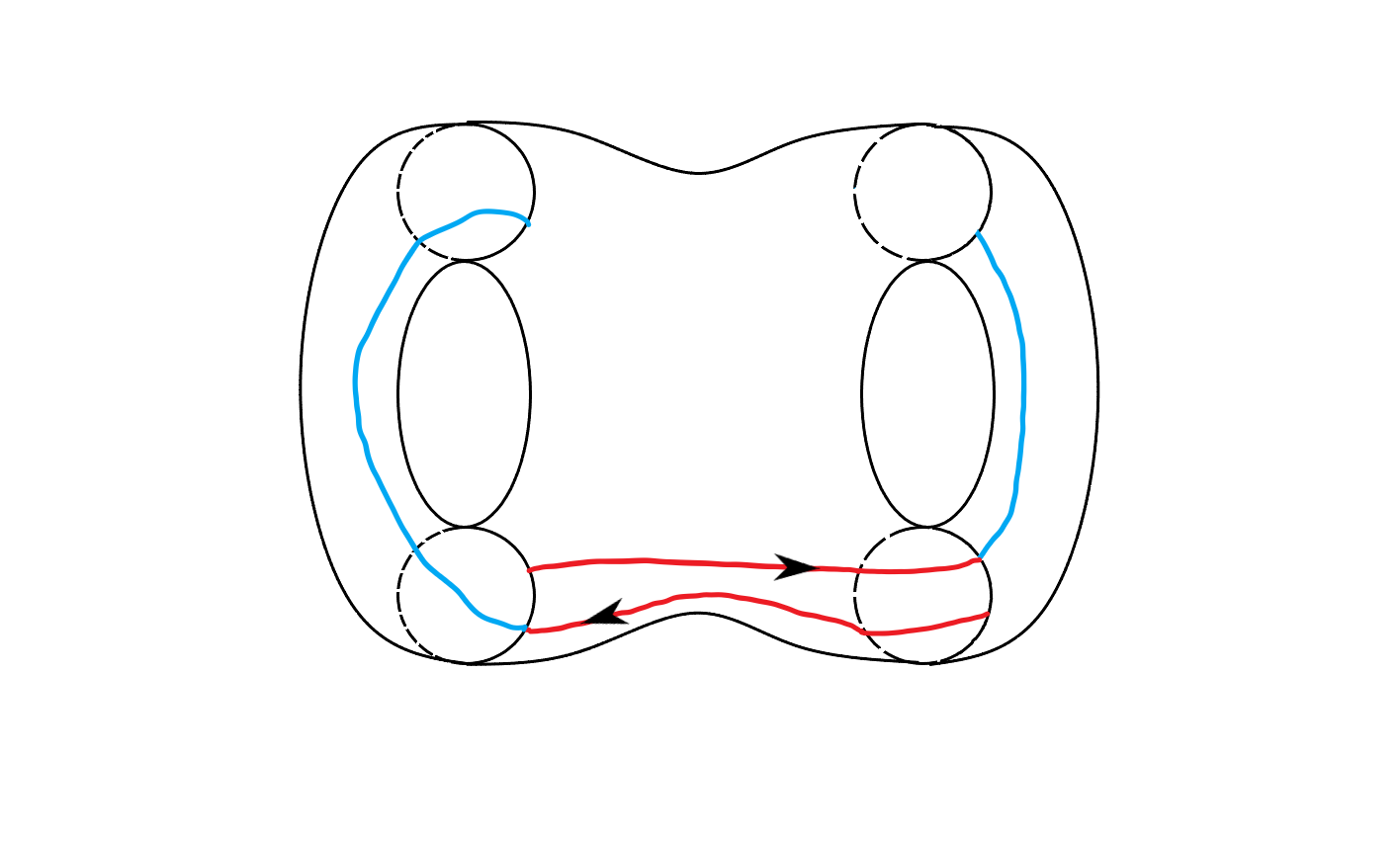}
  \caption{}
\end{subfigure}
\begin{subfigure}{.33\textwidth}
  \centering
  \includegraphics[width=\textwidth]{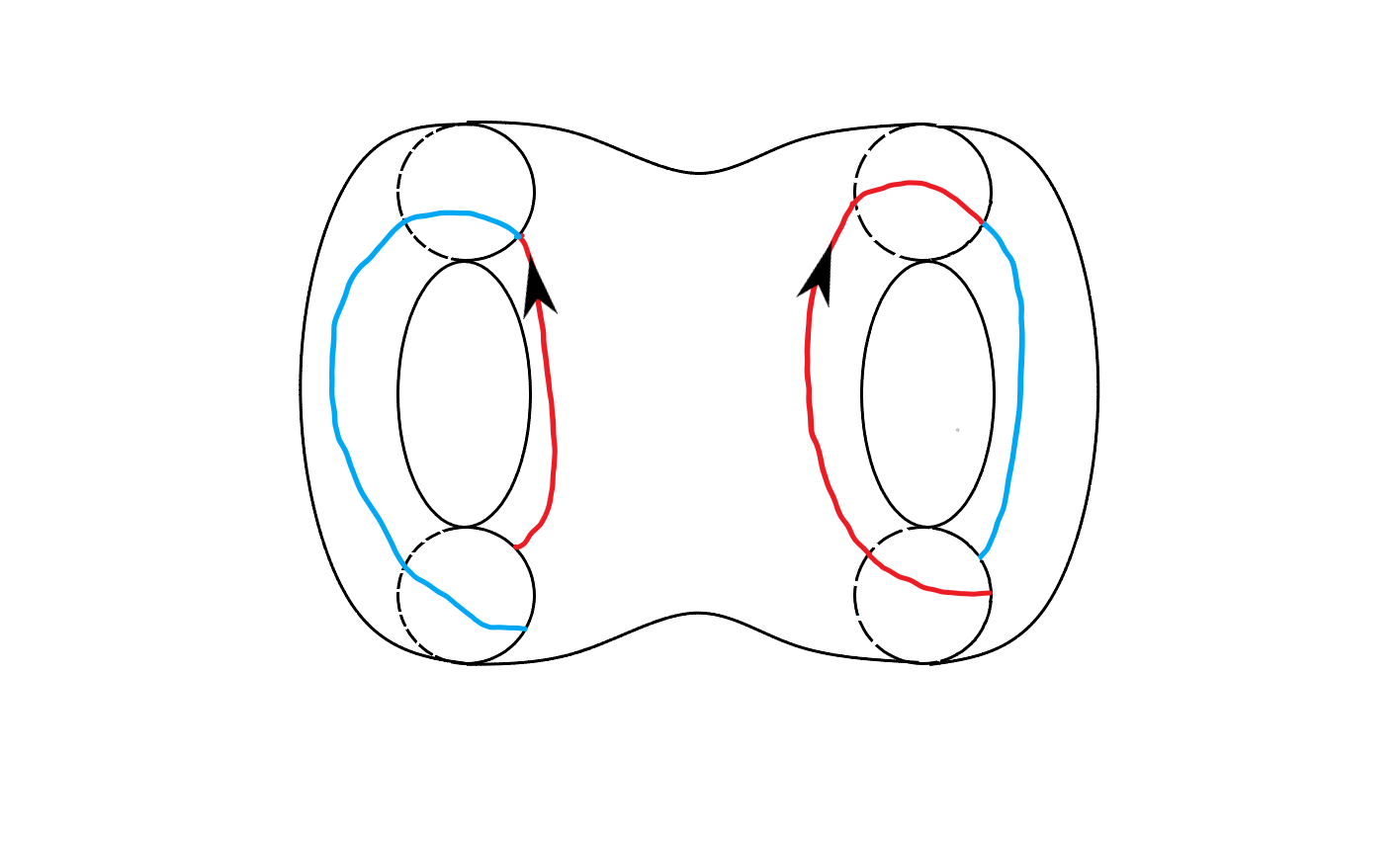}
  \caption{}
\end{subfigure}

\begin{subfigure}{.33\textwidth}
  \centering
  \includegraphics[width=\textwidth]{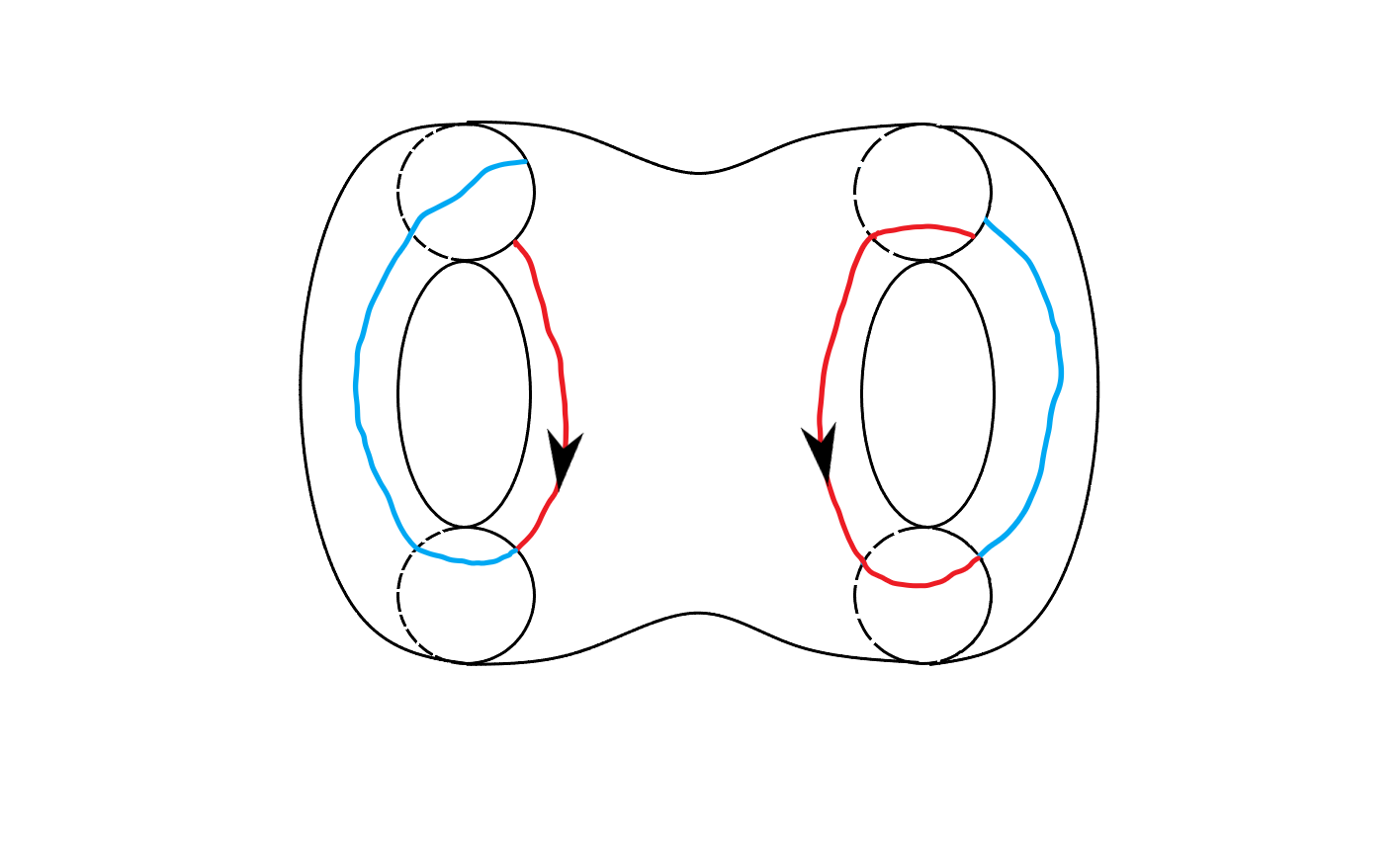}
  \caption{}
\end{subfigure}
\begin{subfigure}{.33\textwidth}
  \centering
  \includegraphics[width=\textwidth]{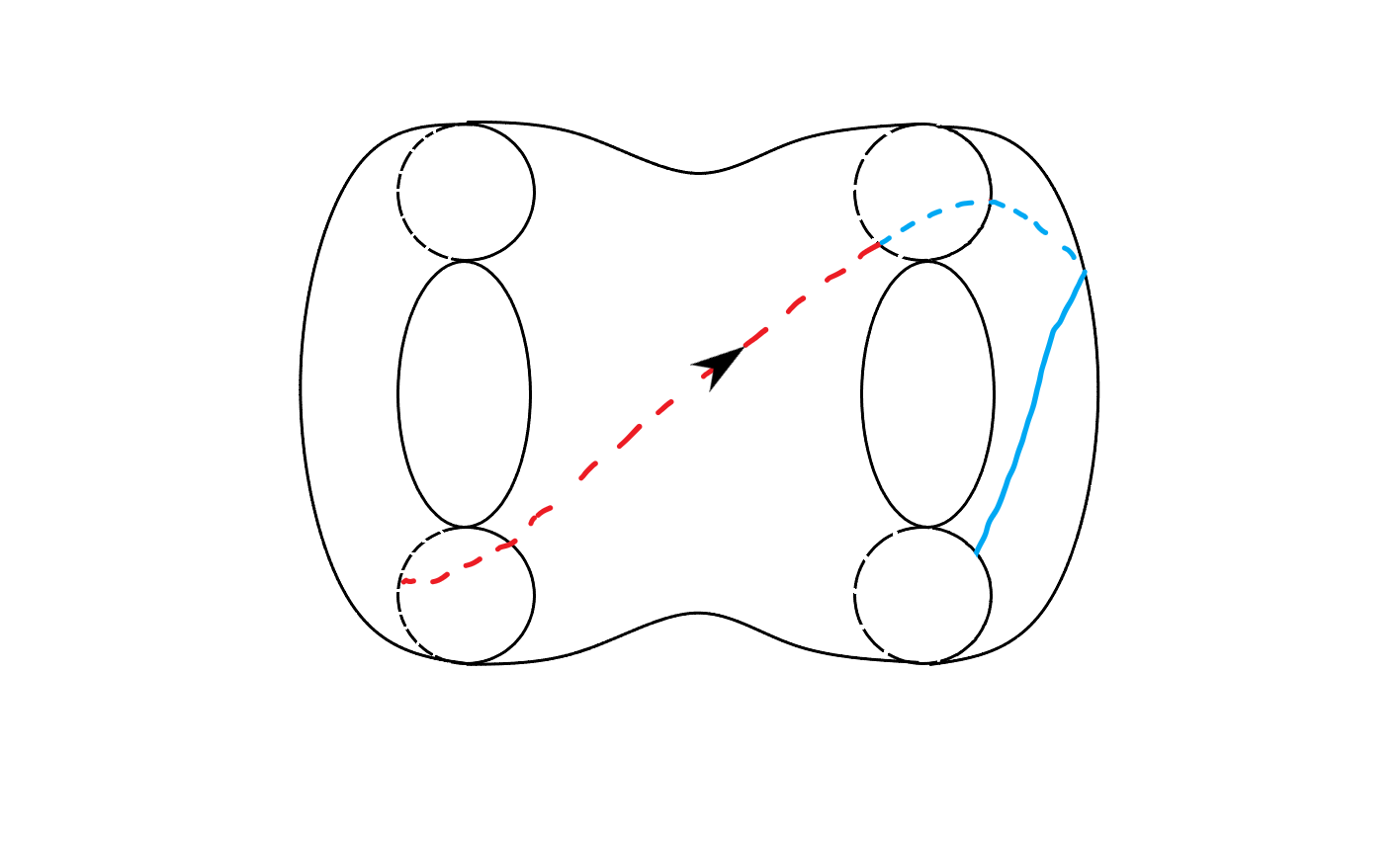}
  \caption{}
\end{subfigure}
\begin{subfigure}{.33\textwidth}
  \centering
  \includegraphics[width=\textwidth]{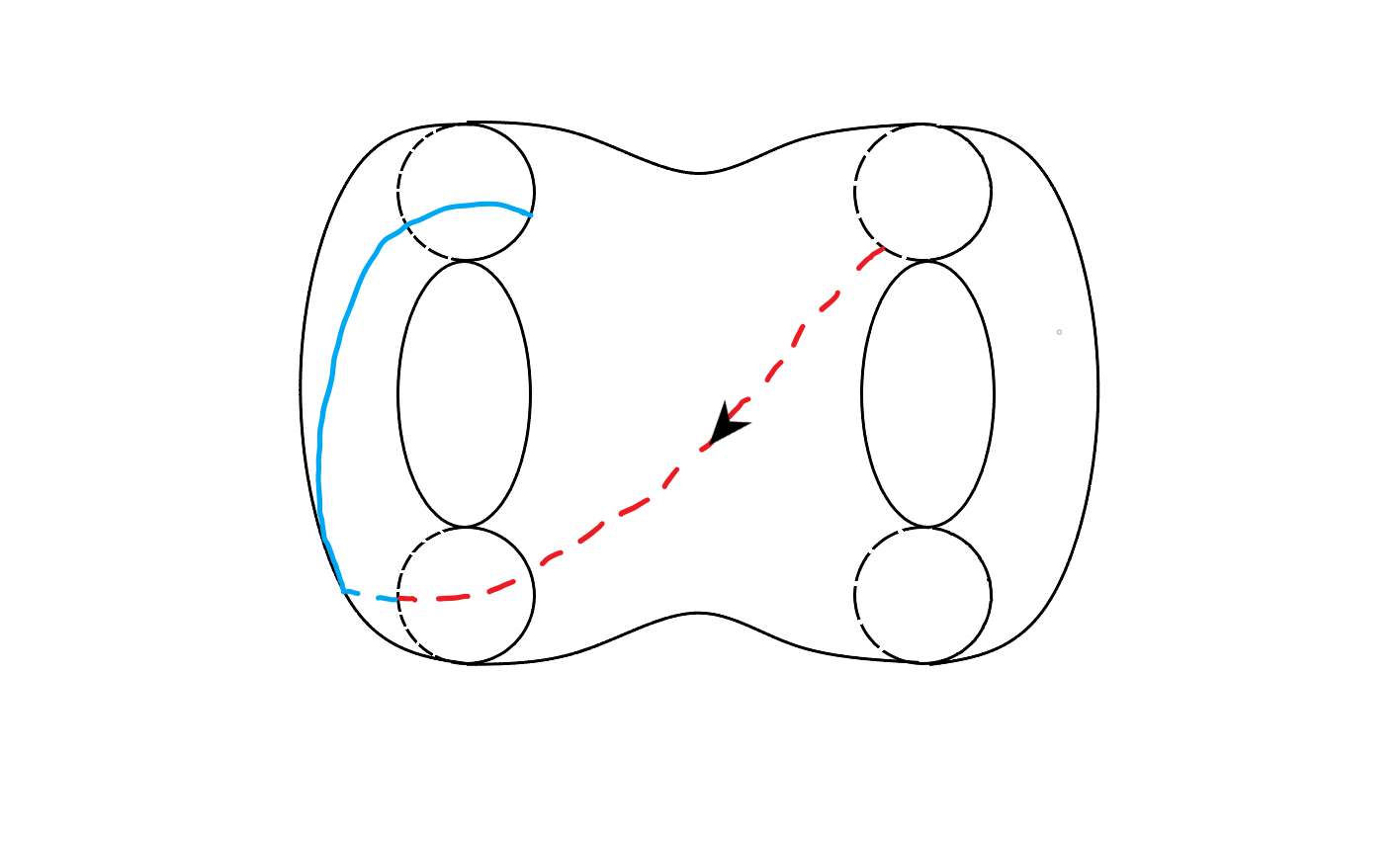}
  \caption{}
\end{subfigure}

\begin{subfigure}{.33\textwidth}
  \centering
  \includegraphics[width=\textwidth]{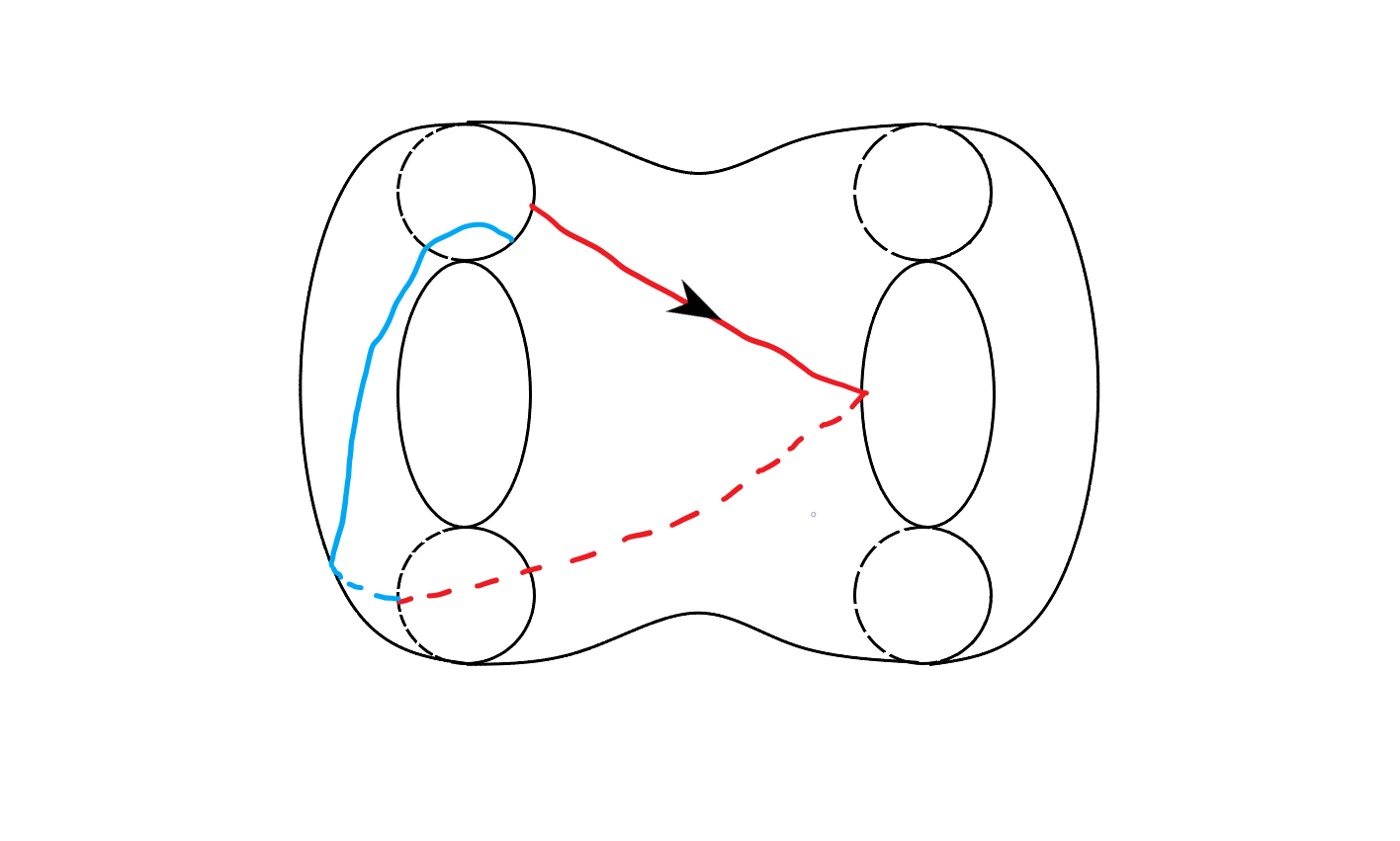}
  \caption{}
  \label{arco_c_dietro}
\end{subfigure}
\begin{subfigure}{.33\textwidth}
  \centering
  \includegraphics[width=\textwidth]{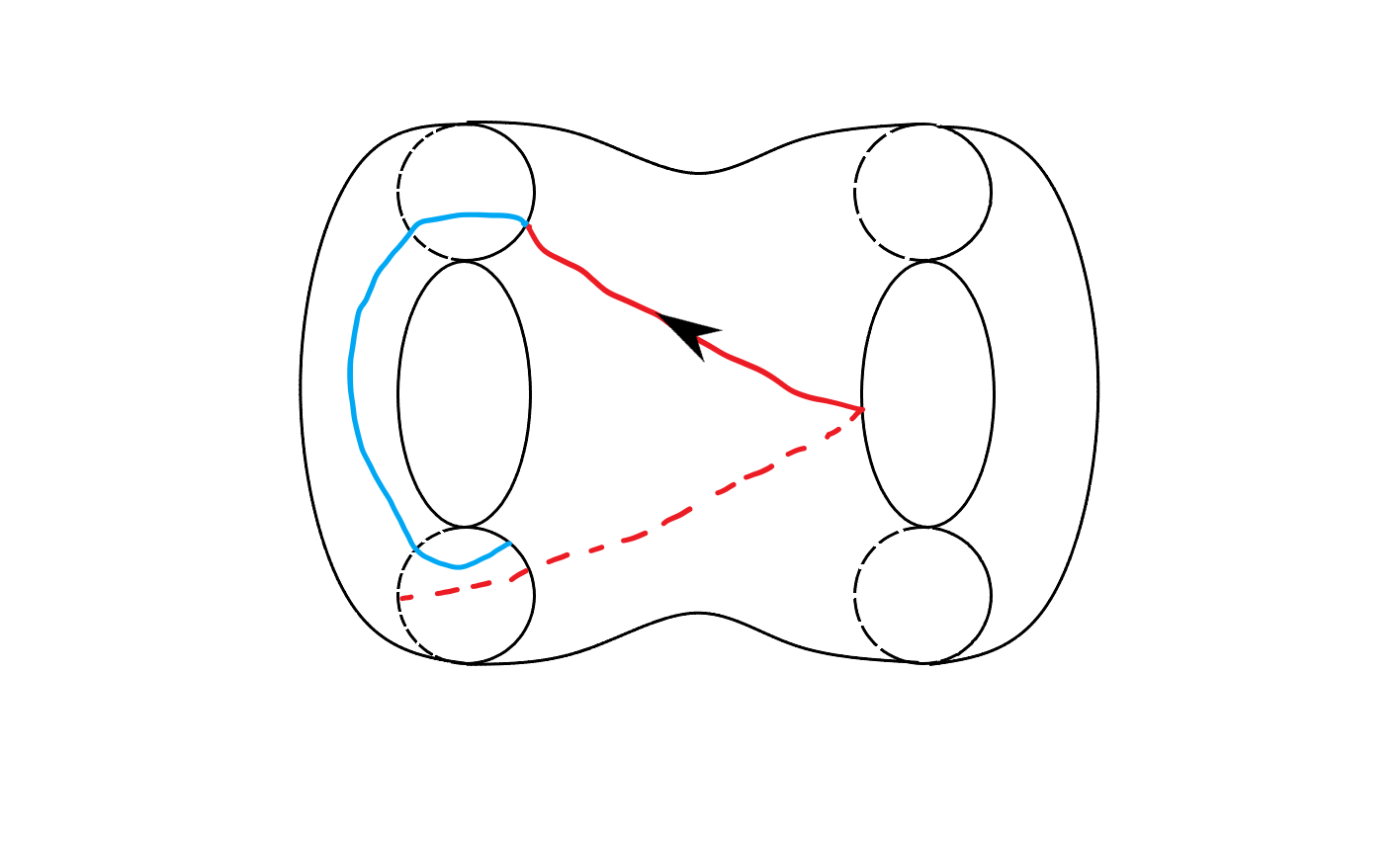}
  \caption{}
\end{subfigure}

\caption{The remaining cases without winding.}
\label{fig:fig_intro}
\end{figure}

We represent in red the arcs of $R_f$ on this surface (see Fig.\ref{fig:333222_only} for an example). Clearly, according to these identifications, some arcs (or part of), as well as some vertices, are contained in the "back part" of the surface (and so are dashed in figures) while some others are in the "front one" (and are not dashed in figures);  more precisely: (i) all the $b$ arcs connecting $C_{Tr}$ and $C_{Bl}$, and their endpoints, are in the back part and (ii) the $c$ arcs connecting $C_{Tl}$ and $C_{Bl}$ are as depicted in Fig.\ref{arco_c_dietro}, that is, their endpoints on $C_{Tl}$ are in the front while those on $C_{Bl}$ are in the back. All the remaining vertices and arcs are in the front part. 

On the surface we have also  arcs, represented in blue, arising from the identification of \(C_i^+\) with \(C_i^-\), \(i = 0, 1\). Clearly those arcs are contained in the two handles bounded by \(C_{Tl}\) and \(C_{Bl}\) on the left and \(C_{Tr}\) and \(C_{Br}\) on the right, and connect the couple of vertices with the same label in  the corresponding top  and bottom circles, eventually winding around the handle (see an example in Fig.\ref{fig:333222_1}). 

To avoid ambiguity we  describe precisely  how the connections are realized. 
Following their orientation, denote with \(X\) (resp. \(Y\))  the first (resp. last) vertex in the front part of both \(C_{Bl}\) and \(C_{Tr}\), as well as the corresponding vertices in \(C_{Tl}\) and \(C_{Br}\). Now denote with  \(L\): (i) the first (resp. last) vertex of \(C_{Tl}\) connected a bottom circle, and (ii) the last (resp. first) vertex of \(C_{Br}\) connected to a left circle in the first three cases (resp.  in case \ref{prop:dati_4}) of Proposition \ref{prop:dati}, as well as the corresponding vertices in \(C_{Bl}\) and \(C_{Tr}\). Moreover, given two vertices \(A, B\) on one circle, we denote with \(\overset{\frown}{AB}\) the oriented arc from \(A\) to \(B\) (endpoints included). 

Now we focus on the first three cases of Proposition \ref{prop:dati}. If \(L\in \overset{\frown}{XY}\),  connect: (i) each vertex in \(\overset{\frown}{LX}\subset C_{Ti}\), with the exception of \(X\), with the   corresponding vertex in $C_{Bi}$  winding once along the handle in the direction of $b_i$, with $i=l,r$, (ii) all the other vertices without winding. Otherwise, that is if \(L \notin \overset{\frown}{XY}\), connect: (i)  all vertices in \(\overset{\frown}{XY}\) without winding, (ii) the vertices in the internal part of \(\overset{\frown}{YL}\subset C_{Ti}\) winding once along the handle in the direction opposite to  $b_i$, (iii) the remaining vertices of  $C_{Ti}$ winding once in the direction of  $b_i$, with $i=l,r$. 
In case \ref{prop:dati_4},  we do the following changes with respect to the previous construction: if \(L\in \overset{\frown}{XY}\), the vertices labelled $L$ in  both circles are connected without winding and, in both cases, all the connections are done winding in the opposite directions. 

\begin{figure}[h!]
\begin{subfigure}{0.33\textwidth}
  \centering
  \includegraphics[width=\textwidth]{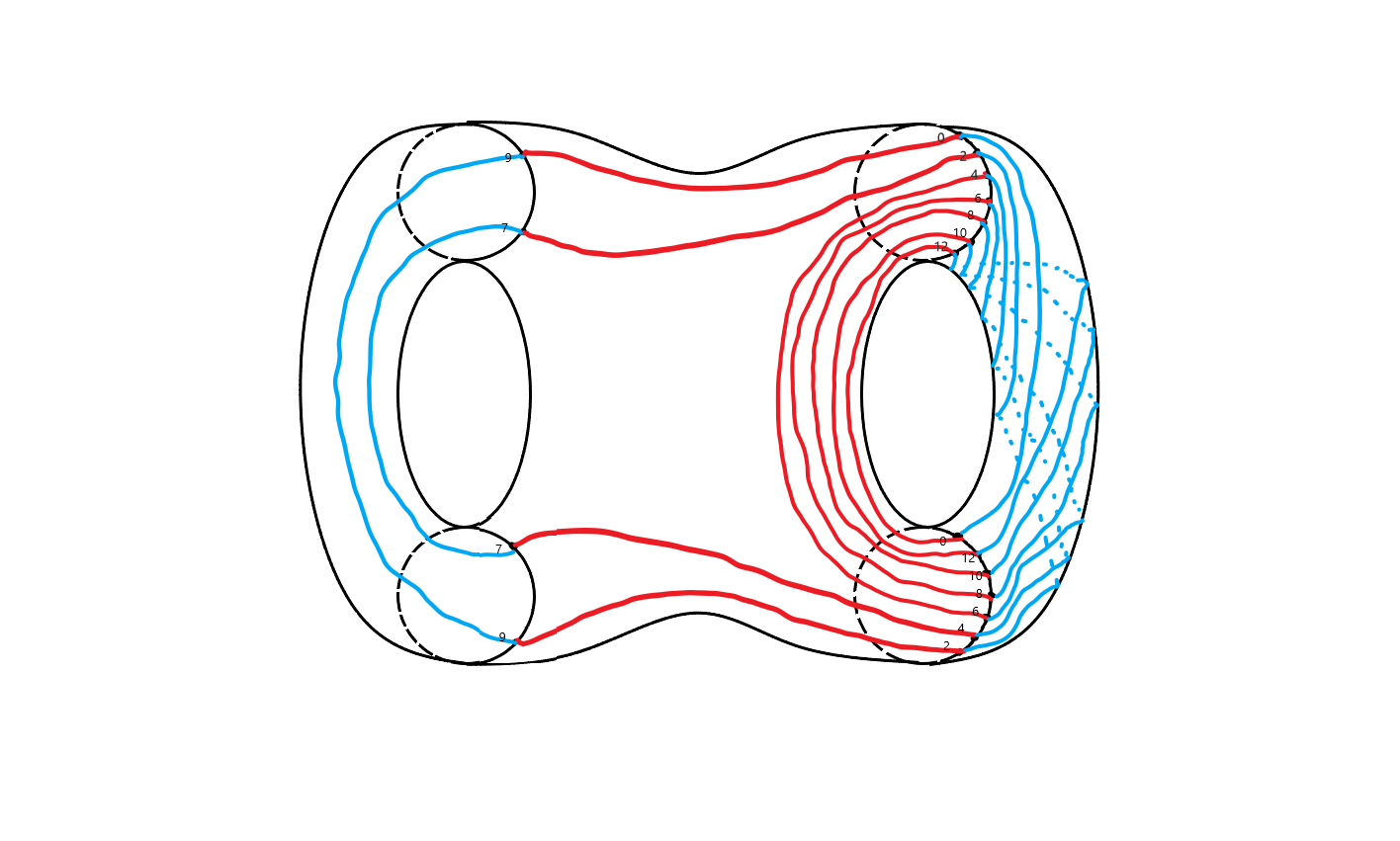}
\end{subfigure}
\begin{subfigure}{0.33\textwidth}
  \centering
  \includegraphics[width=\textwidth]{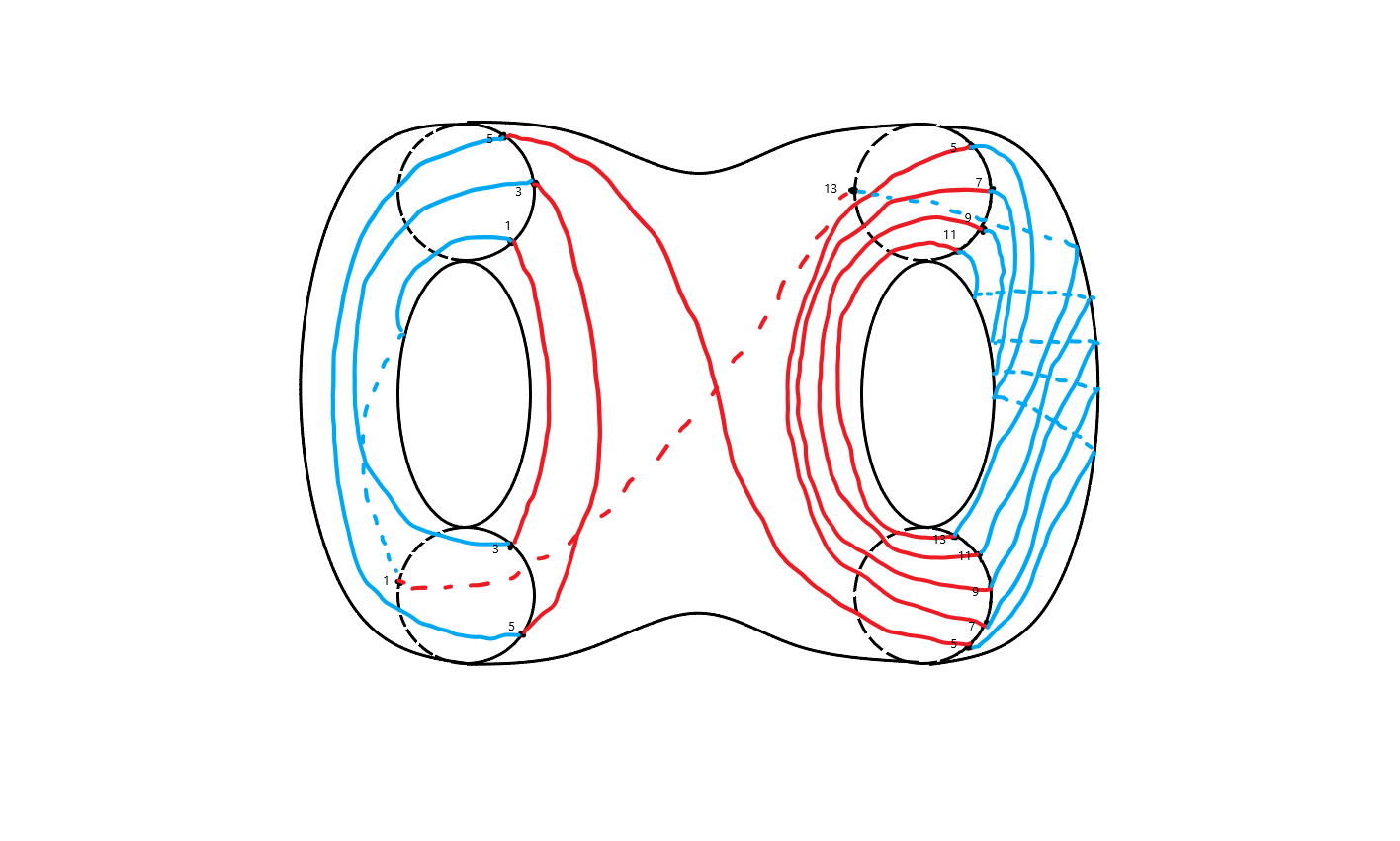}
\end{subfigure}
\begin{subfigure}{0.33\textwidth}
  \centering
  \includegraphics[width=\textwidth]{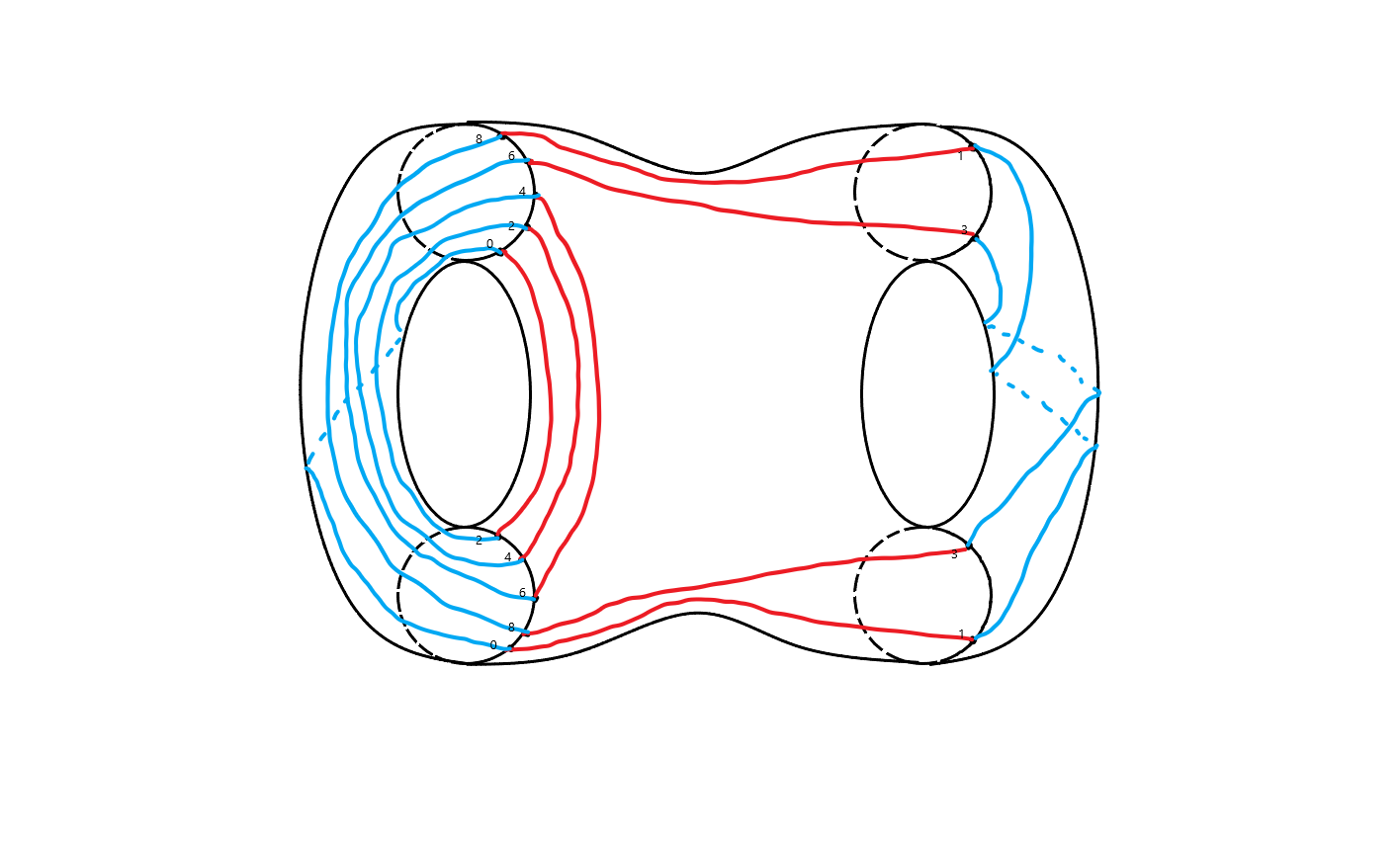}
\end{subfigure}
\caption{The curves \(x_f^0, x_f^1\) and \(x_f^2\) of \((4, 6, 10, 1, 1, 1)\).}
\label{fig:4610111_res}
\end{figure}

In this way we obtain three closed curves on the surface each corresponding to one of  the three colors  of $R_f$.  Note that along each curve the blue and red arcs alternate. We denote this curves with $x^0_f$, $x^1_f $, $x^3_f$. 
See the Fig.\ref{fig:4610111_res} for the three curves  corresponding to $R_f=(4,6,10,1,1,1)$.\\

As we have seen, any choice of a couple of curves among $x^0_f$, $x^1_f$, $x^2_f$ is meridian system for a Heegaard splitting of $M_f$ (being the other the standard one); so, in order to determine the $psl$-moves associated  to such a splitting it is enough to find words $\chi^i_f\in B_{2,2}$, such that $\widehat{\chi^i}_f=x_f^i$, for $i=0,1,2$, and the move takes the form: \[psl_i: \beta \rightarrow \chi_f^i \# \beta = \chi_f^i \beta\] 

\begin{figure}
\begin{subfigure}{0.33\textwidth}
  \centering
  \includegraphics[width=\textwidth]{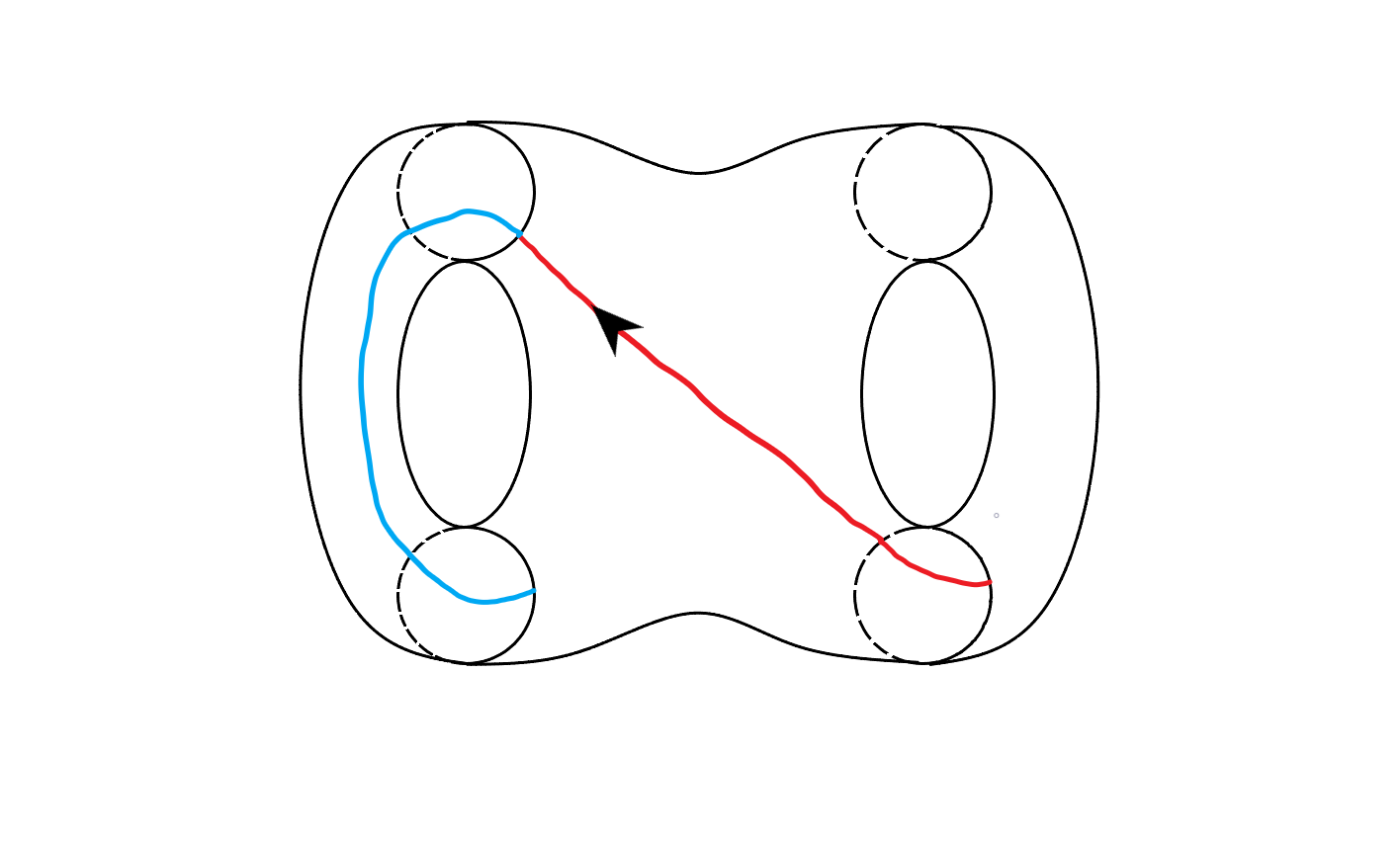}
  \caption{}
  \label{fig:fig_intro_1}
\end{subfigure}
\begin{subfigure}{0.33\textwidth}
  \centering
  \includegraphics[width=\textwidth]{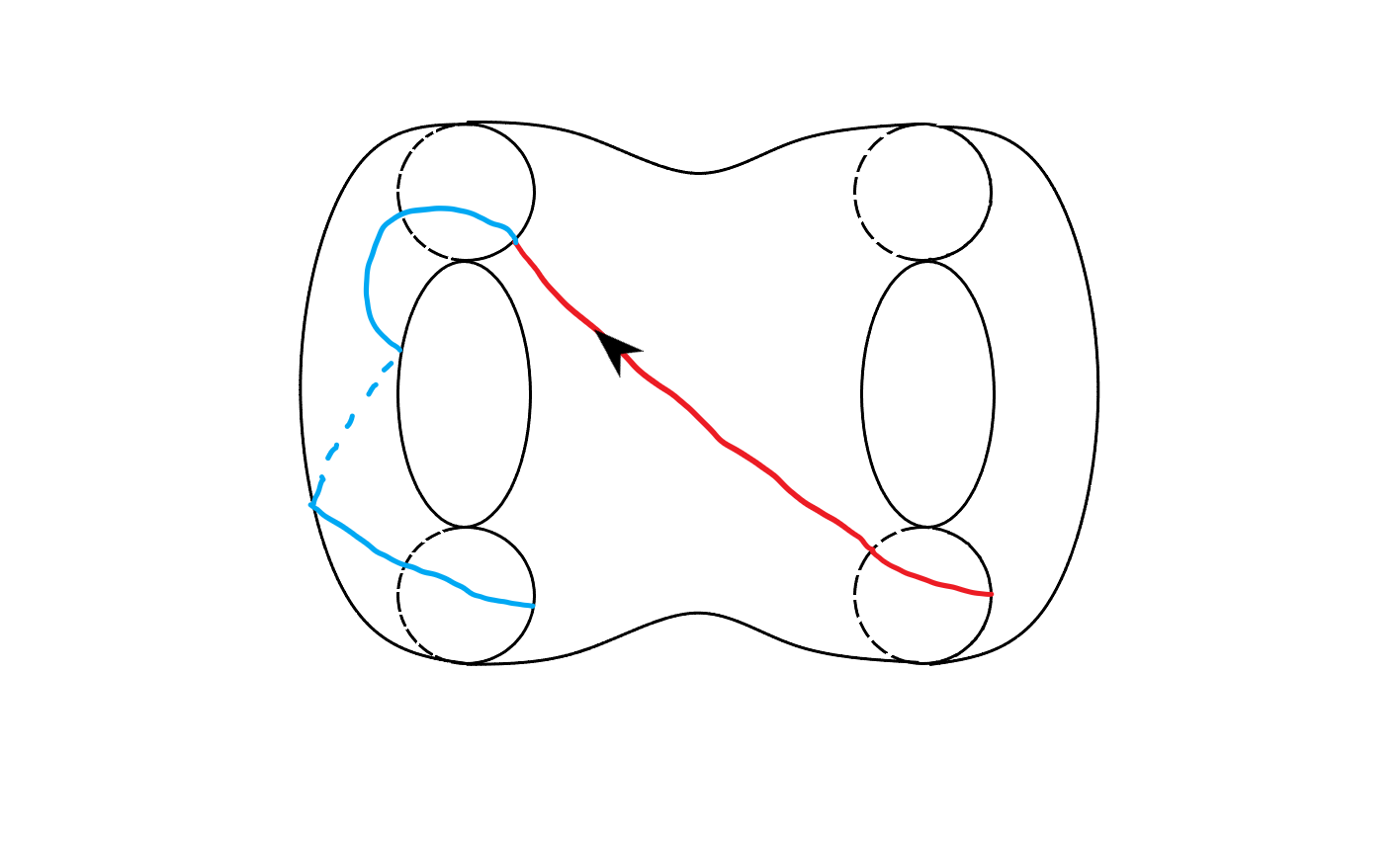}
  \caption{}
  \label{fig:fig_intro_1_a}
\end{subfigure}
\begin{subfigure}{0.33\textwidth}
  \centering
  \includegraphics[width=\textwidth]{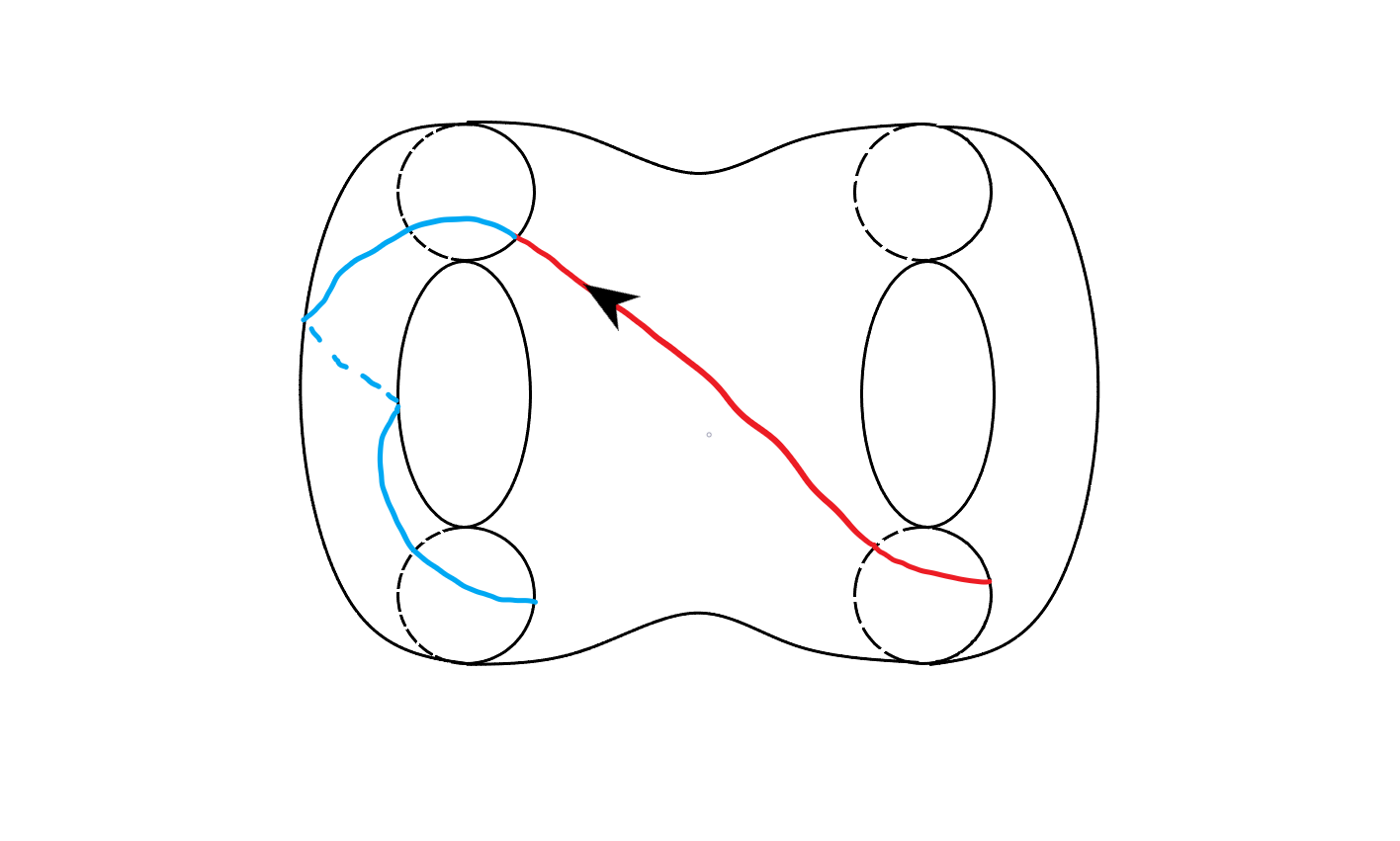}
  \caption{}
  \label{fig:fig_intro_1_b}
\end{subfigure}

\begin{subfigure}{0.33\textwidth}
  \centering
  \includegraphics[width=\textwidth]{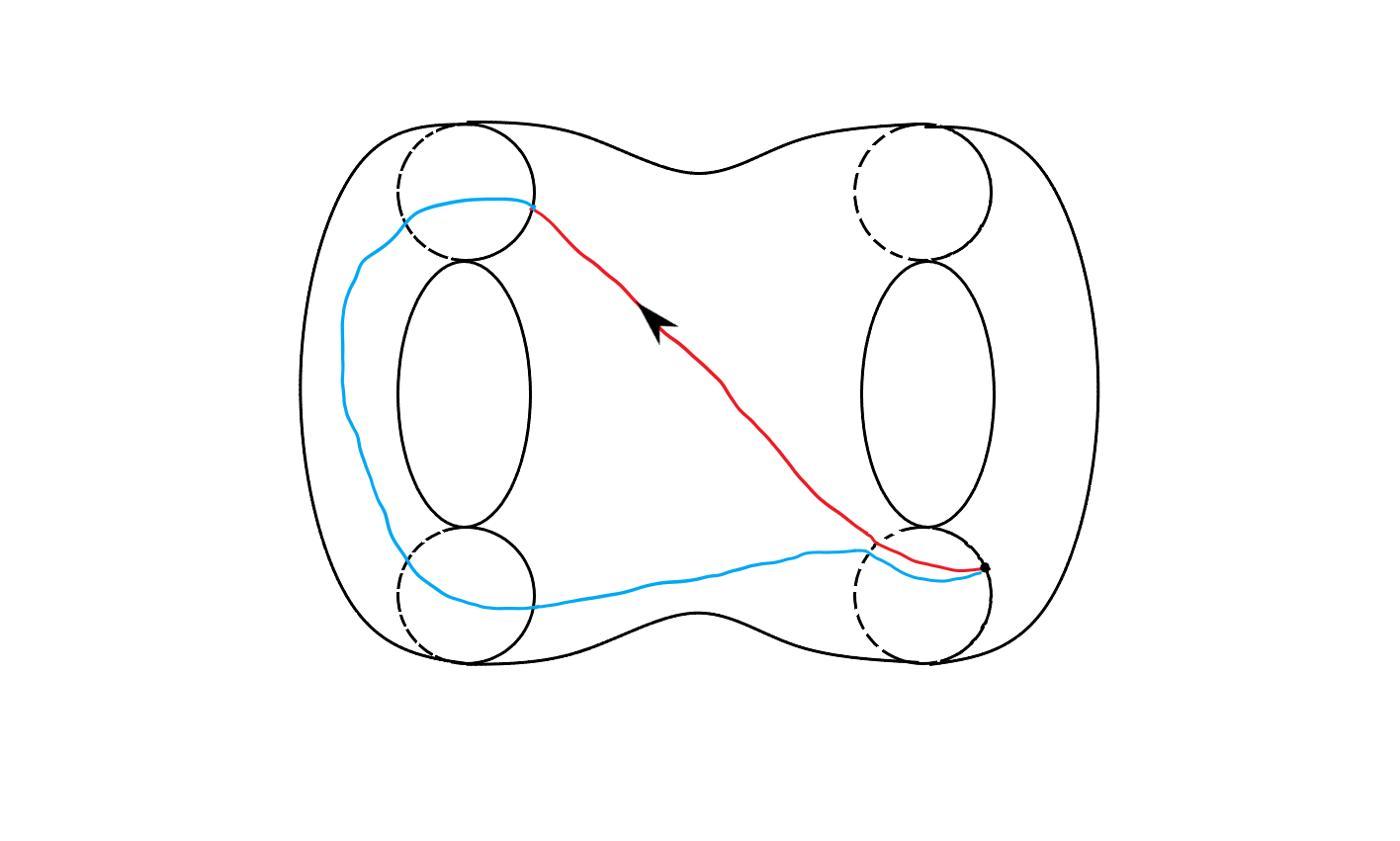}
  \caption{}
  \label{fig:fig_1}
\end{subfigure}
\begin{subfigure}{0.33\textwidth}
  \centering
  \includegraphics[width=\textwidth]{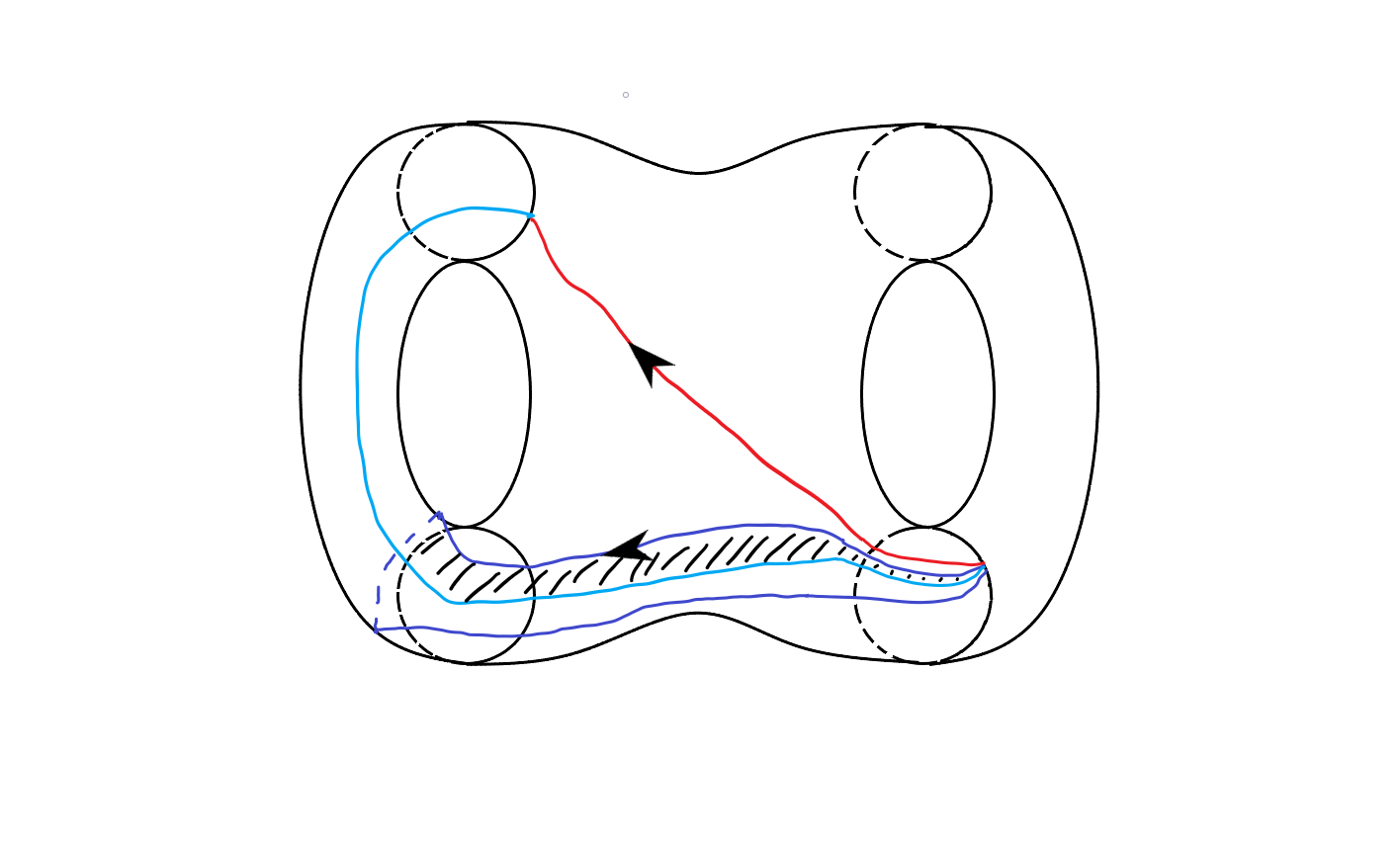}
  \caption{}
  \label{fig:fig_1_a}
\end{subfigure}
\begin{subfigure}{0.33\textwidth}
  \centering
  \includegraphics[width=\textwidth]{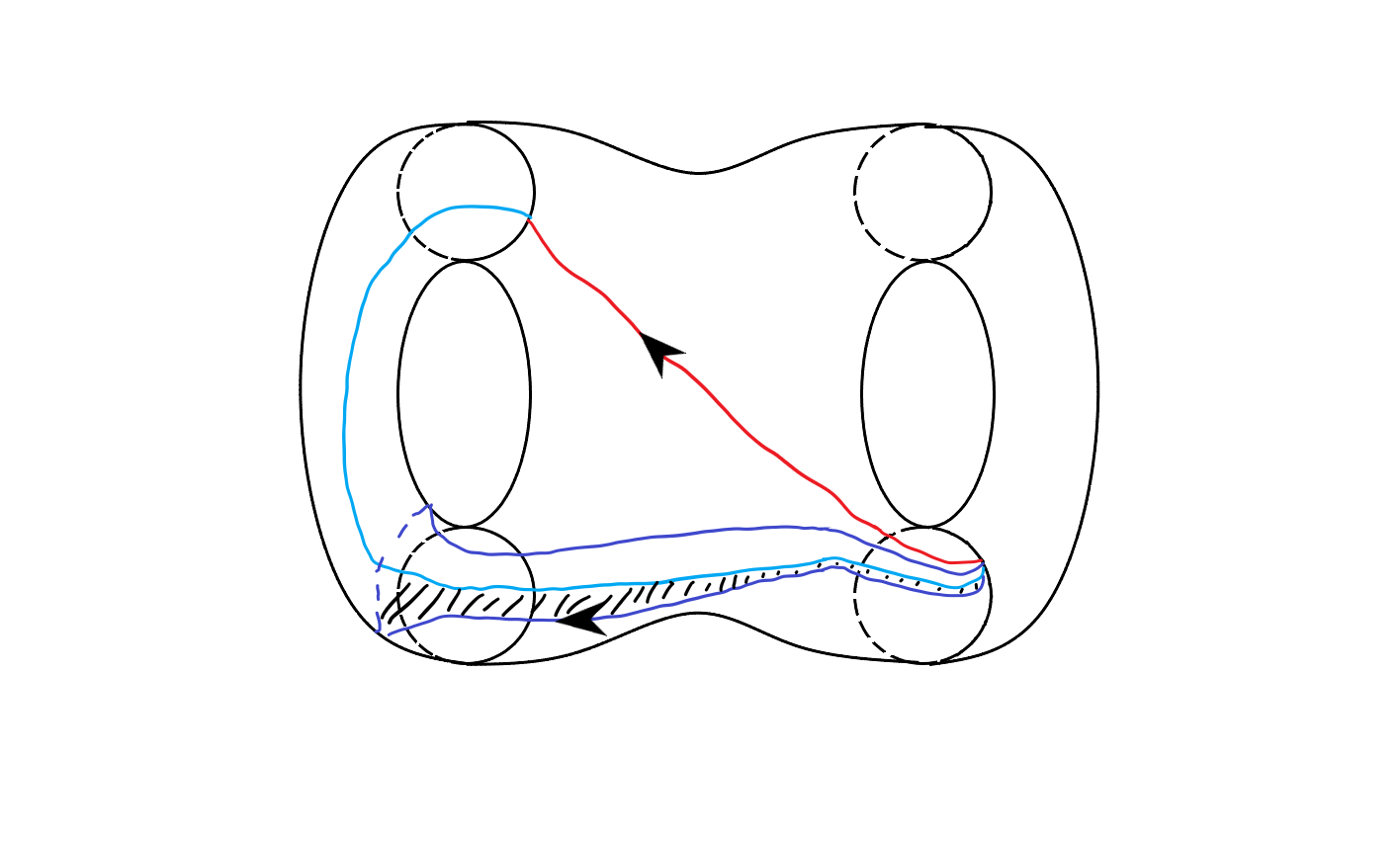}
  \caption{}
  \label{fig:fig_1_b}
\end{subfigure}
\caption{The first case with all the possible windings.}
\label{fig:fig_winding_intro}
\end{figure}

\begin{figure}
\begin{subfigure}{.33\textwidth}
  \centering
  \includegraphics[width=\textwidth]{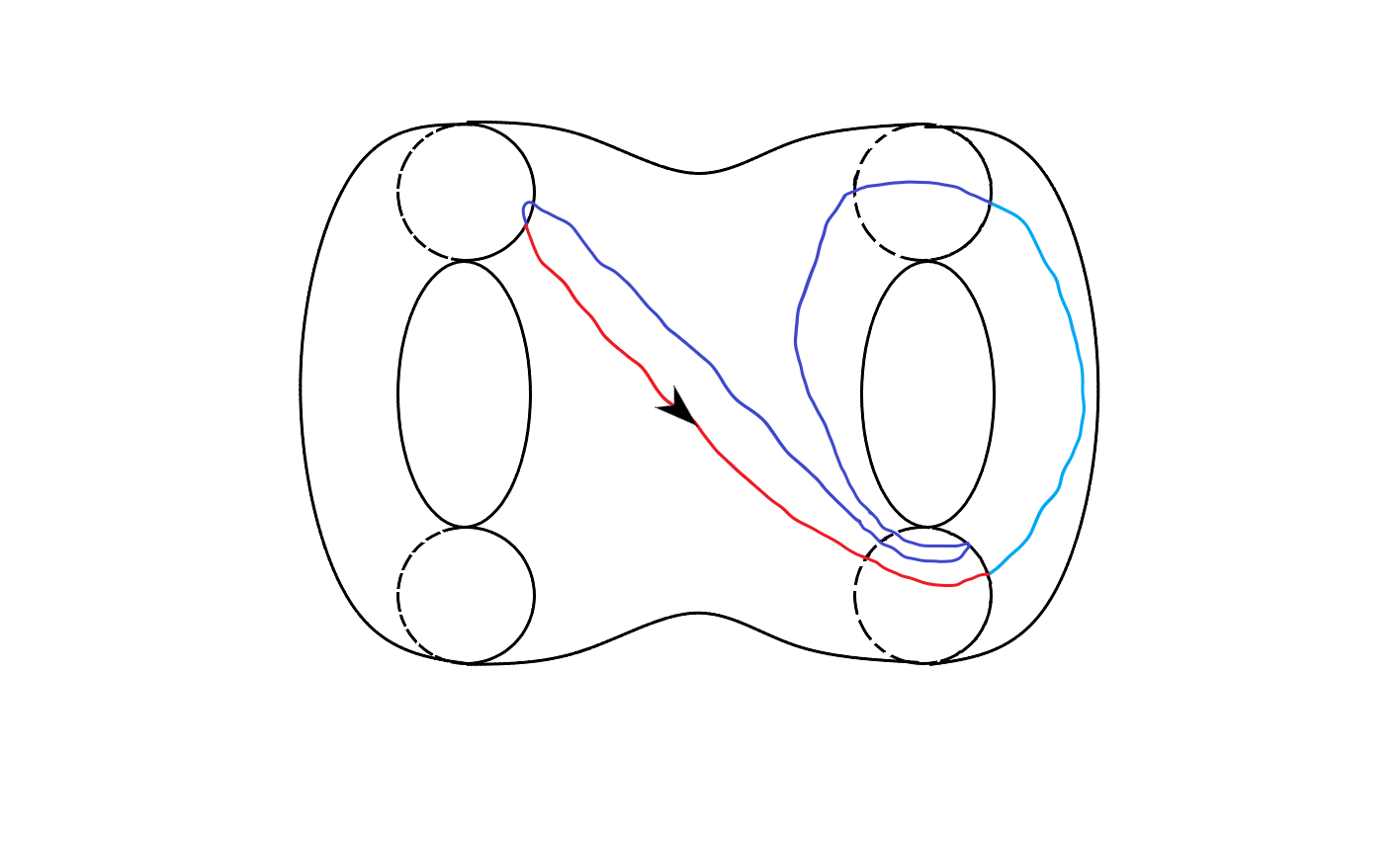}
  \caption{}
\end{subfigure}
\begin{subfigure}{.33\textwidth}
  \centering
  \includegraphics[width=\textwidth]{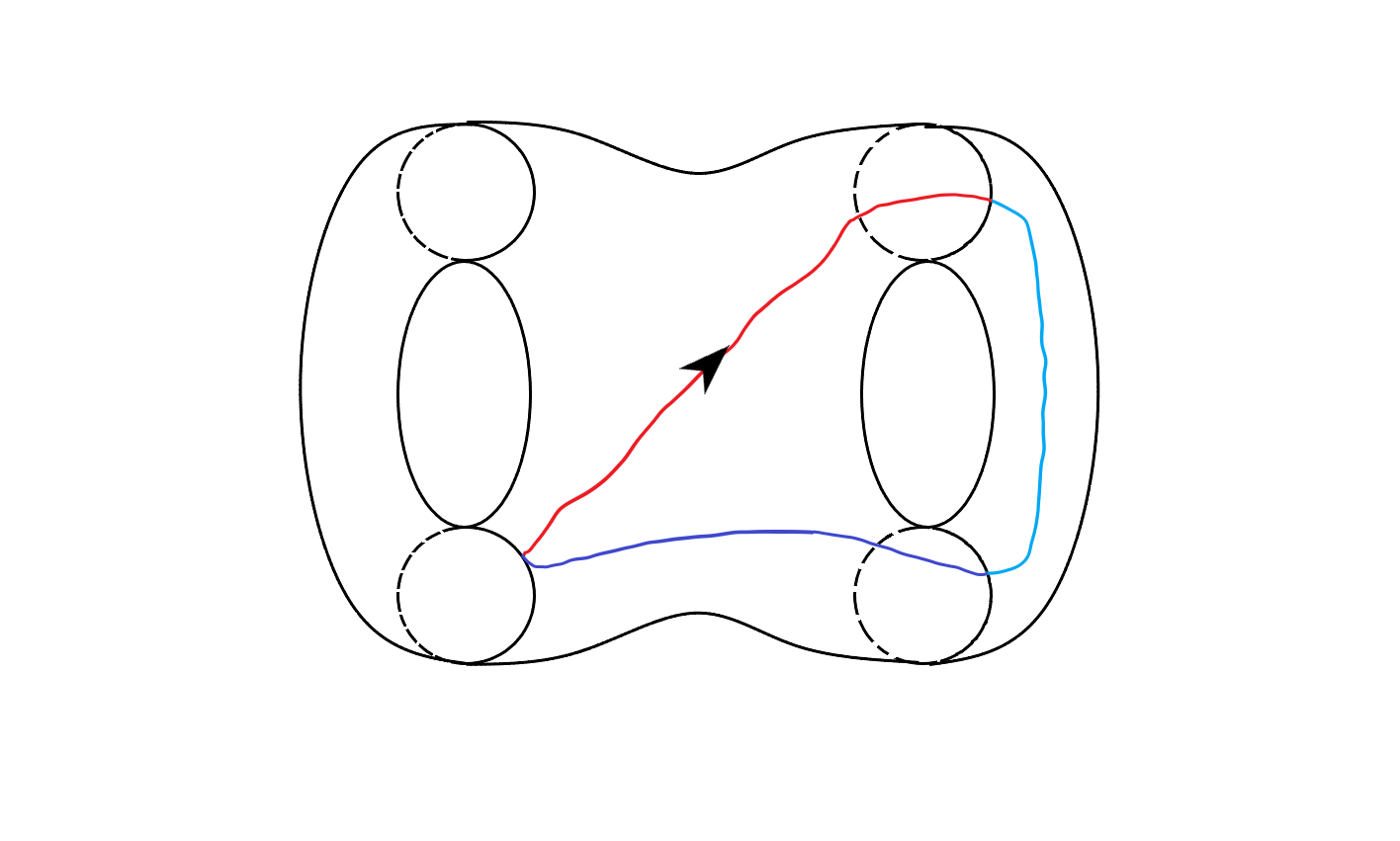}
  \caption{}
\end{subfigure}
\begin{subfigure}{.33\textwidth}
  \centering
  \includegraphics[width=\textwidth]{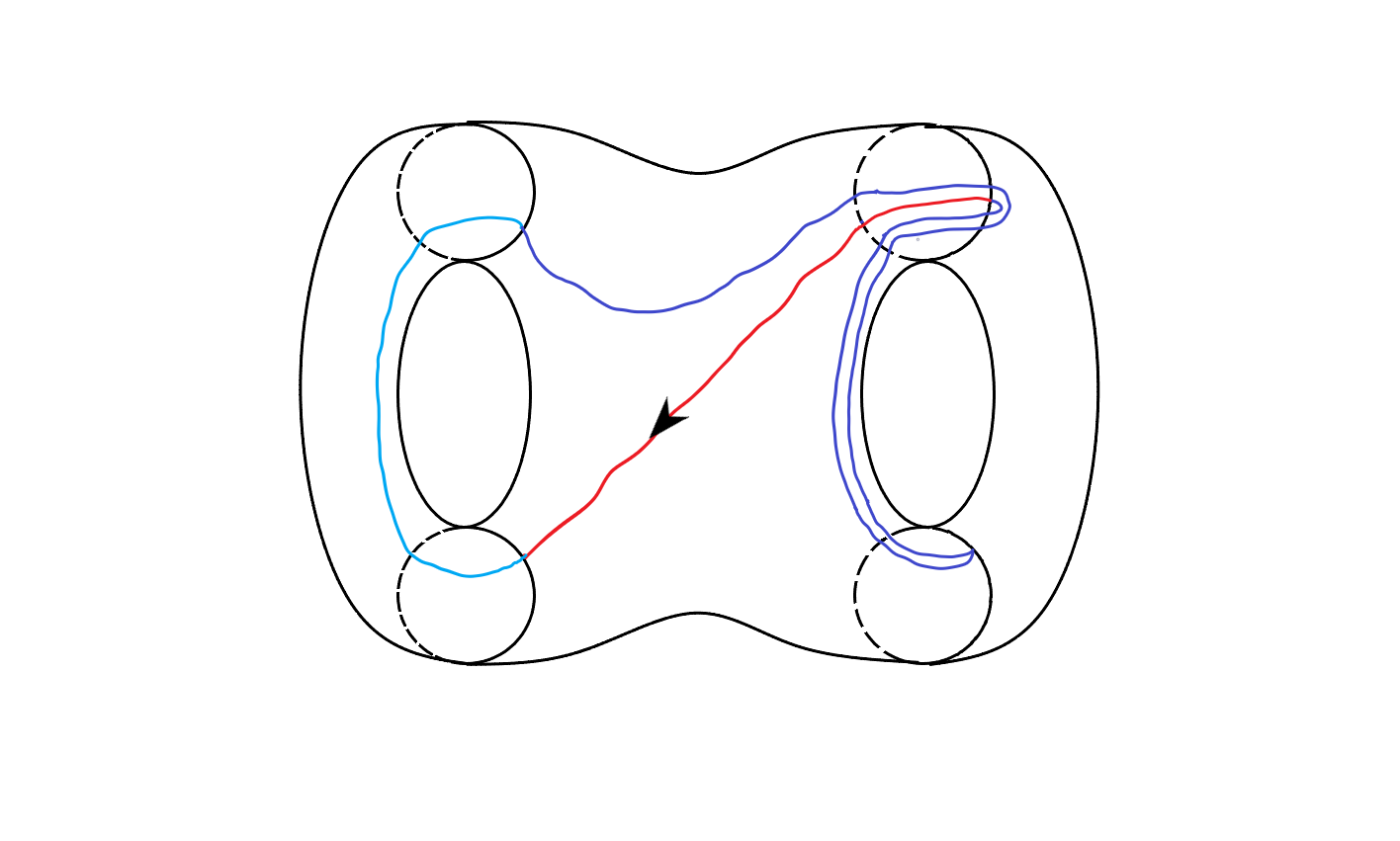}
  \caption{}
\end{subfigure}

\begin{subfigure}{.33\textwidth}
  \centering
  \includegraphics[width=\textwidth]{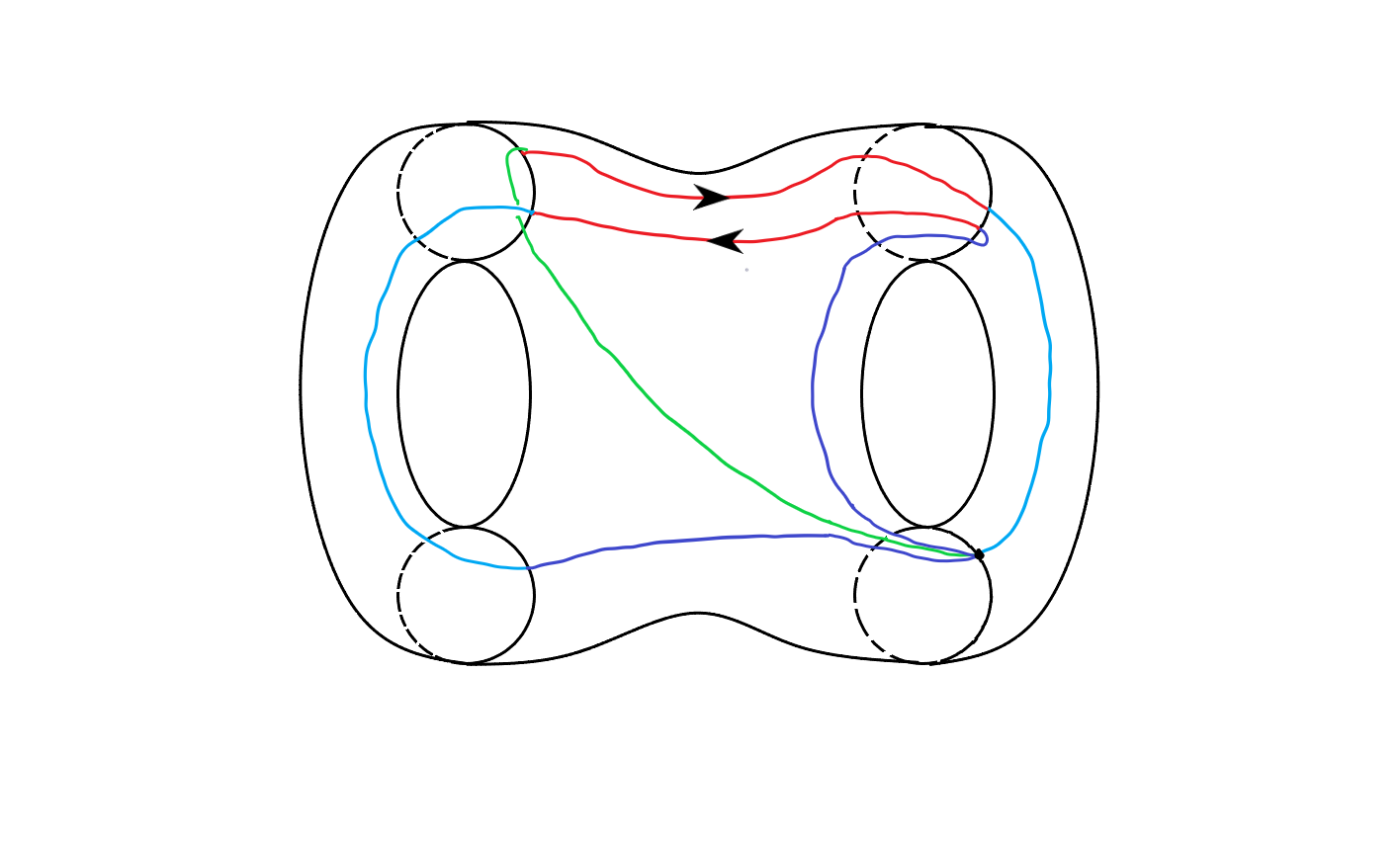}
  \caption{}
\end{subfigure}
\begin{subfigure}{.33\textwidth}
  \centering
  \includegraphics[width=\textwidth]{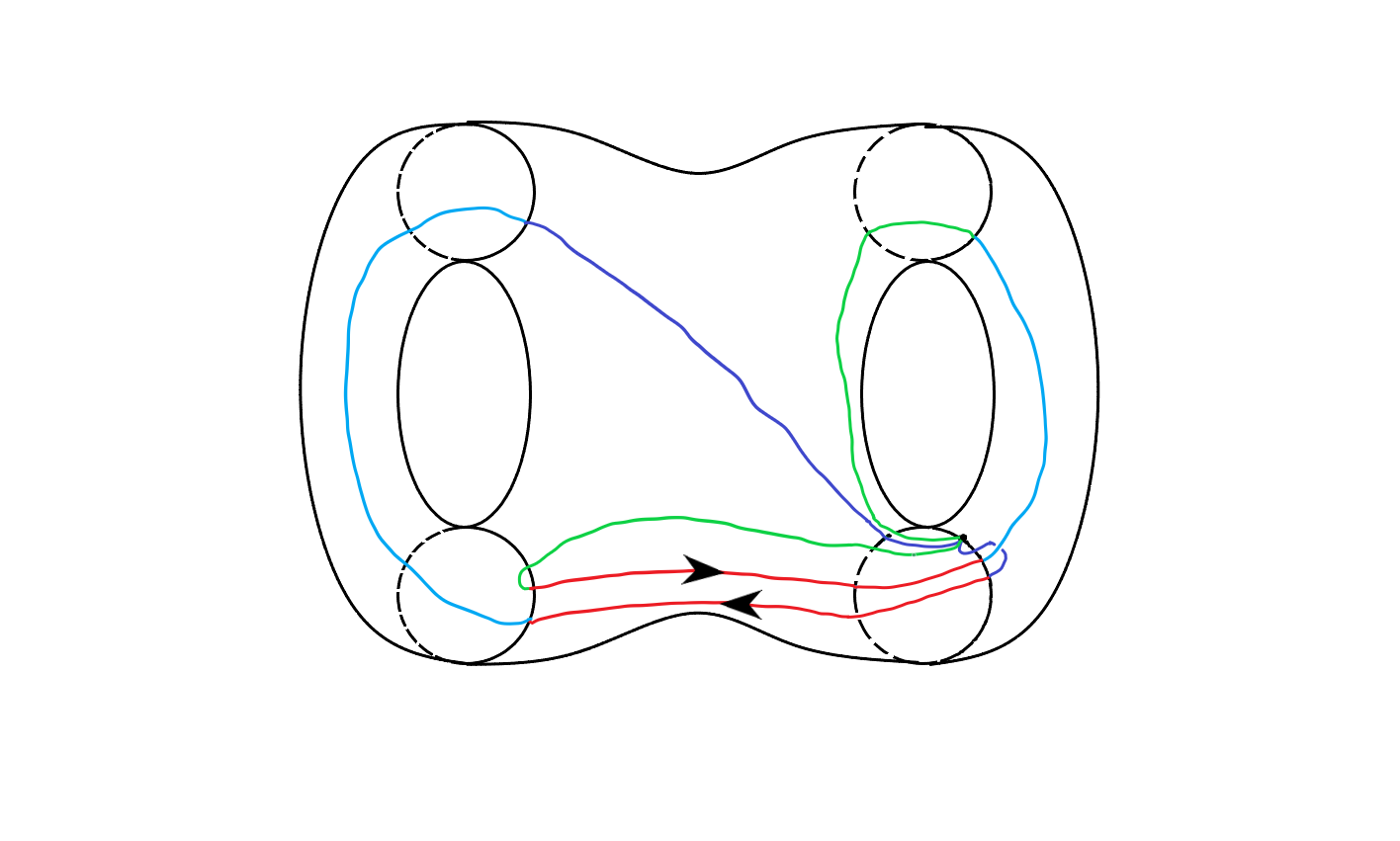}
  \caption{}
\end{subfigure}
\begin{subfigure}{.33\textwidth}
  \centering
  \includegraphics[width=\textwidth]{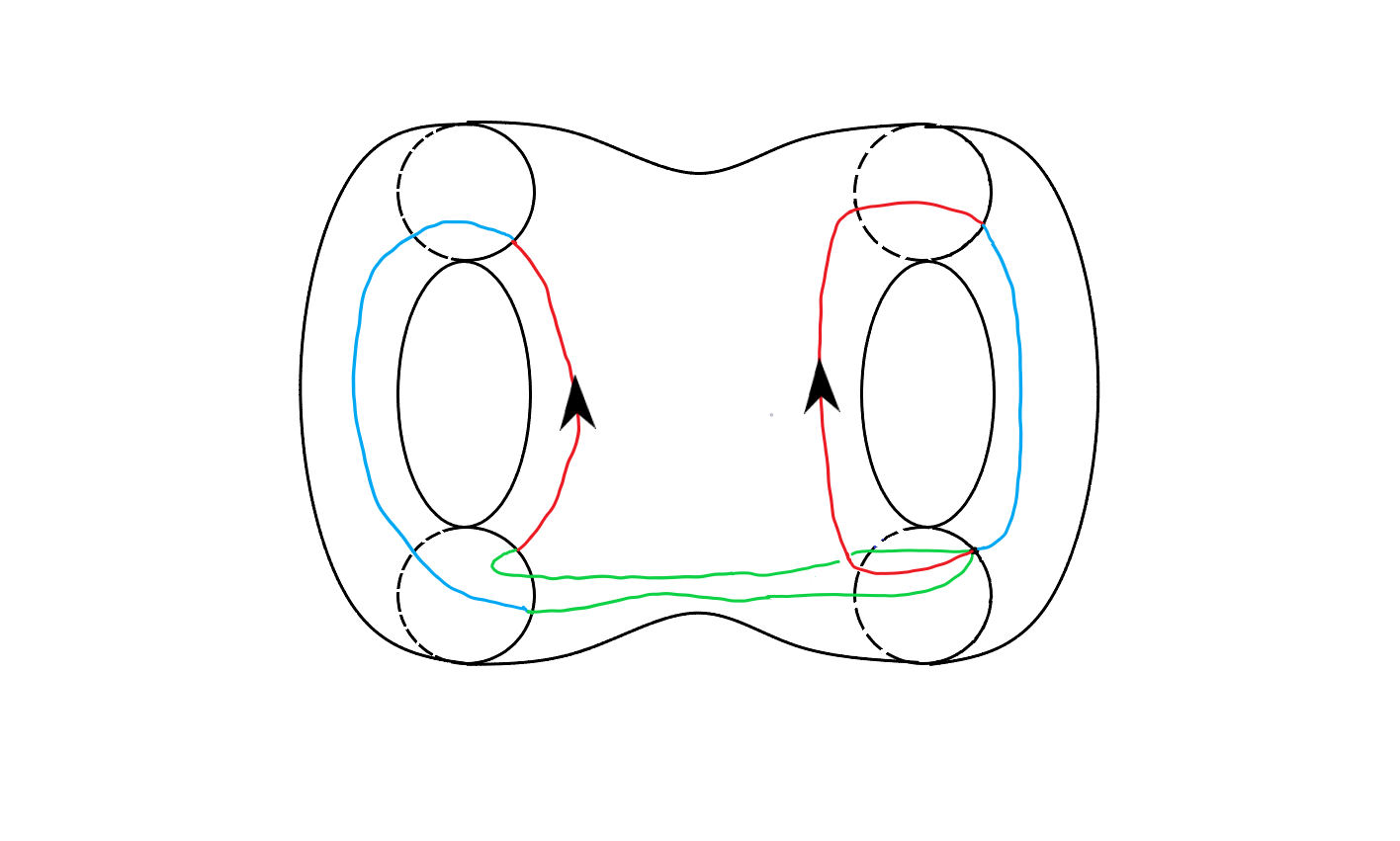}
  \caption{}
\end{subfigure}

\begin{subfigure}{.33\textwidth}
  \centering
  \includegraphics[width=\textwidth]{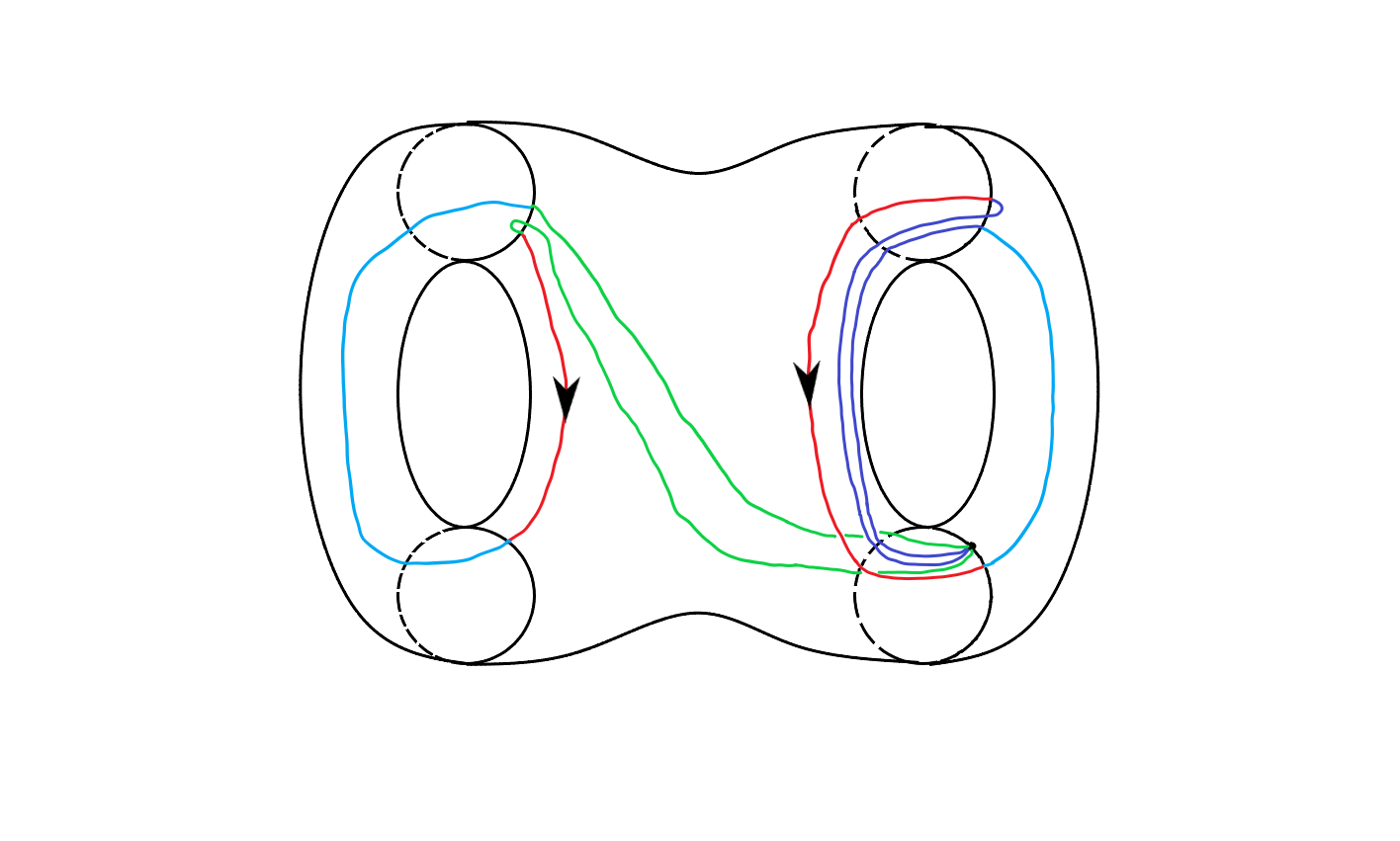}
  \caption{}
\end{subfigure}
\begin{subfigure}{.33\textwidth}
  \centering
  \includegraphics[width=\textwidth]{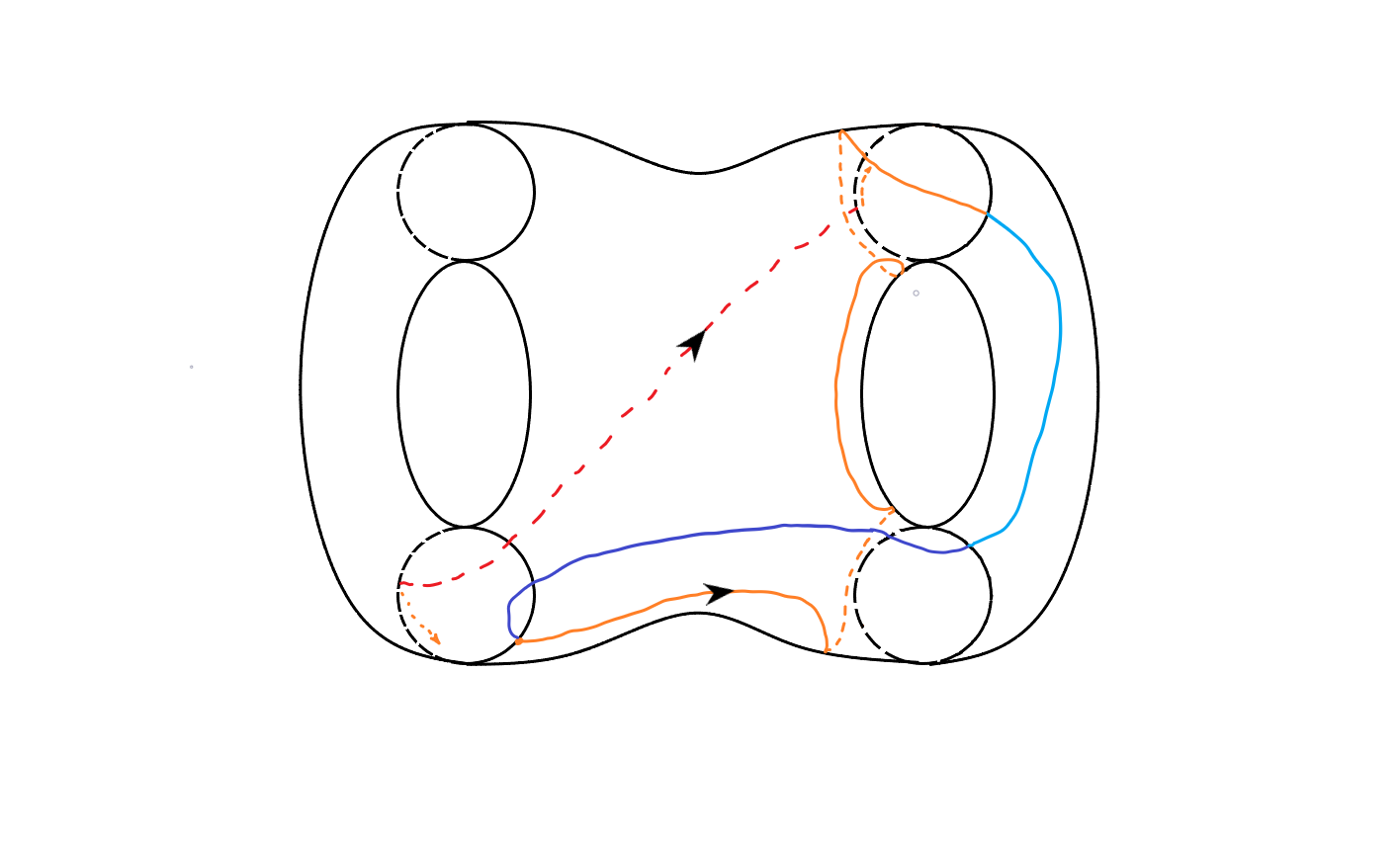}
  \caption{}
\end{subfigure}
\begin{subfigure}{.33\textwidth}
  \centering
  \includegraphics[width=\textwidth]{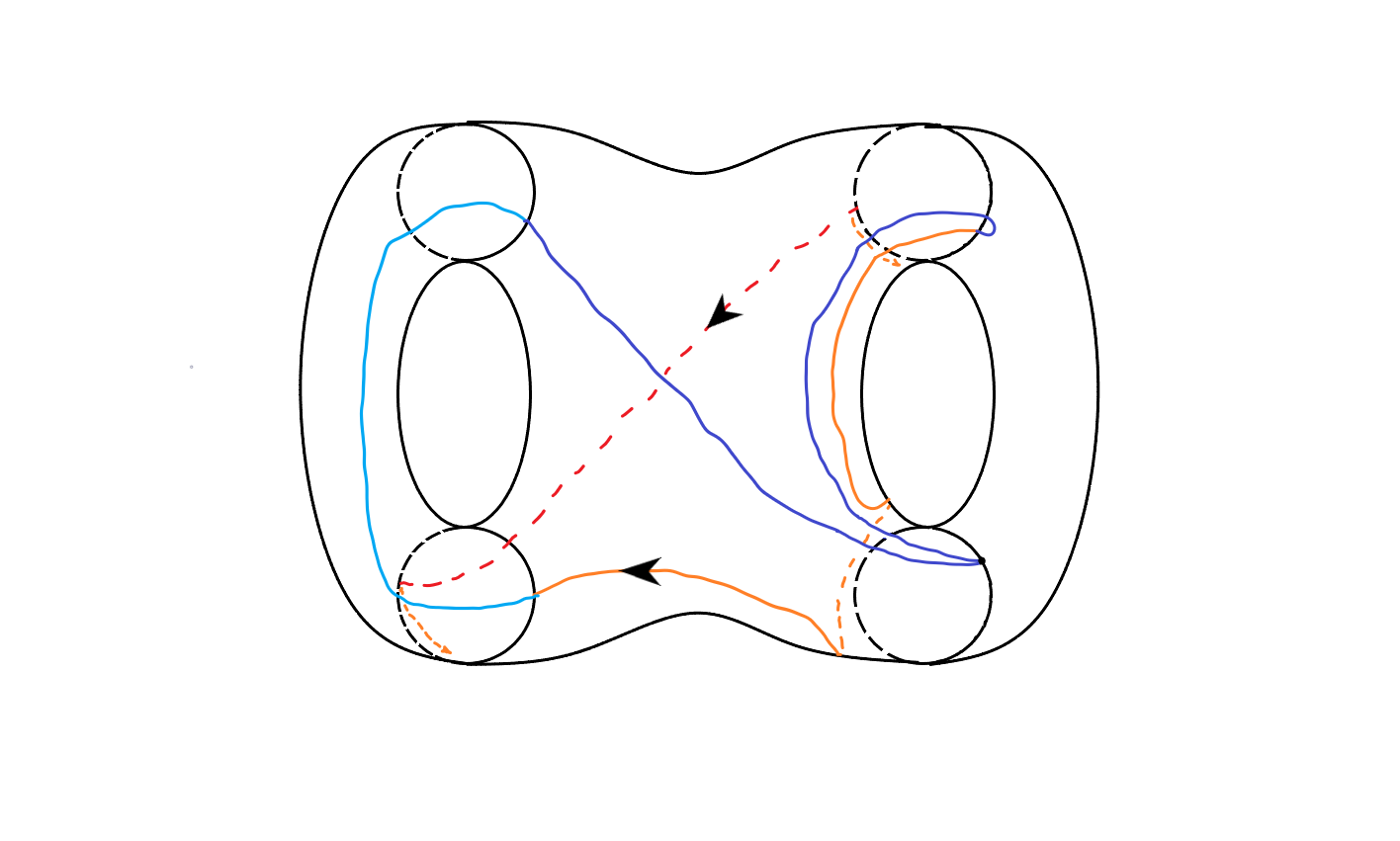}
  \caption{}
\end{subfigure}

\begin{subfigure}{.33\textwidth}
  \centering
  \includegraphics[width=\textwidth]{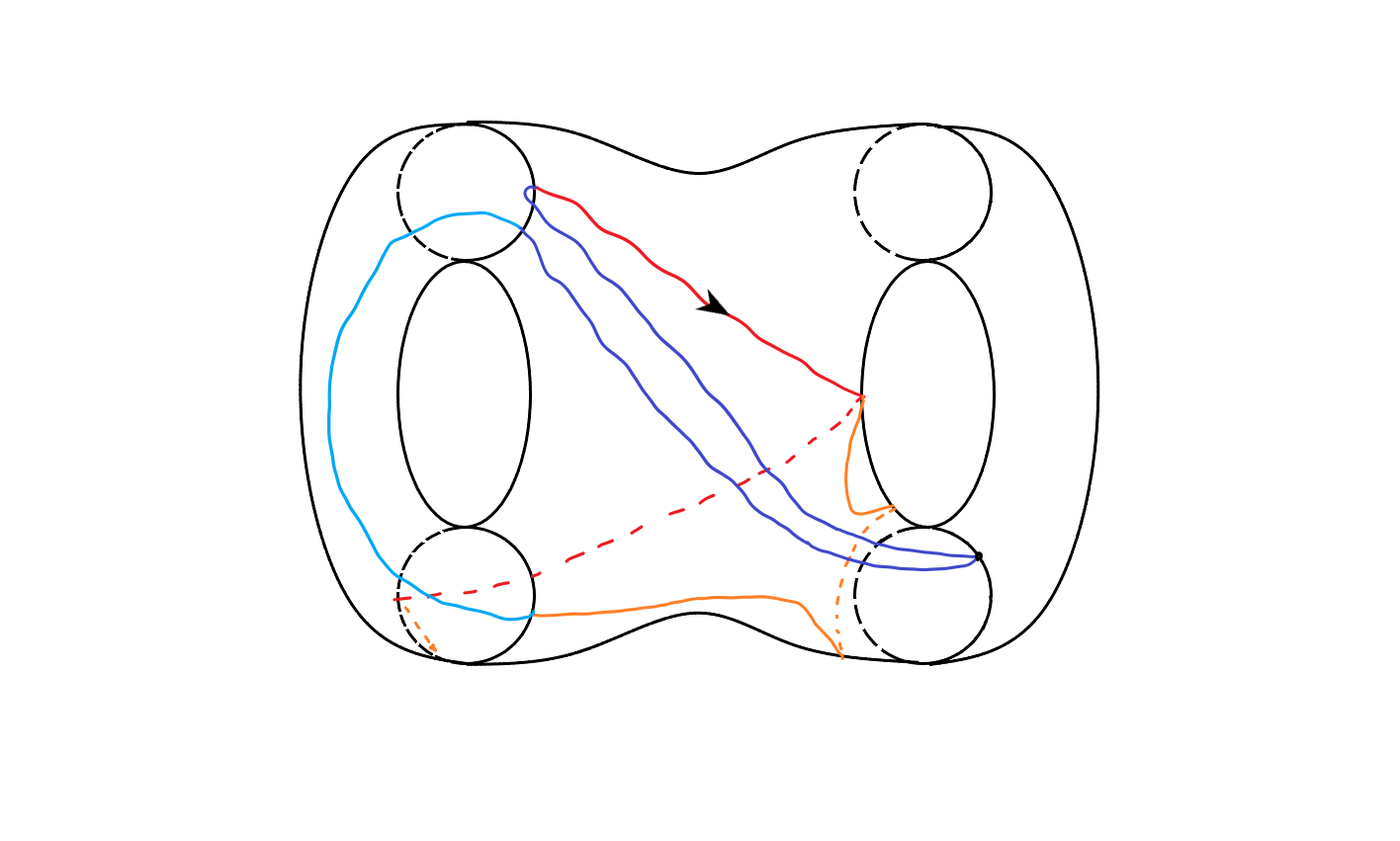}
  \caption{}
\end{subfigure}
\begin{subfigure}{.33\textwidth}
  \centering
  \includegraphics[width=\textwidth]{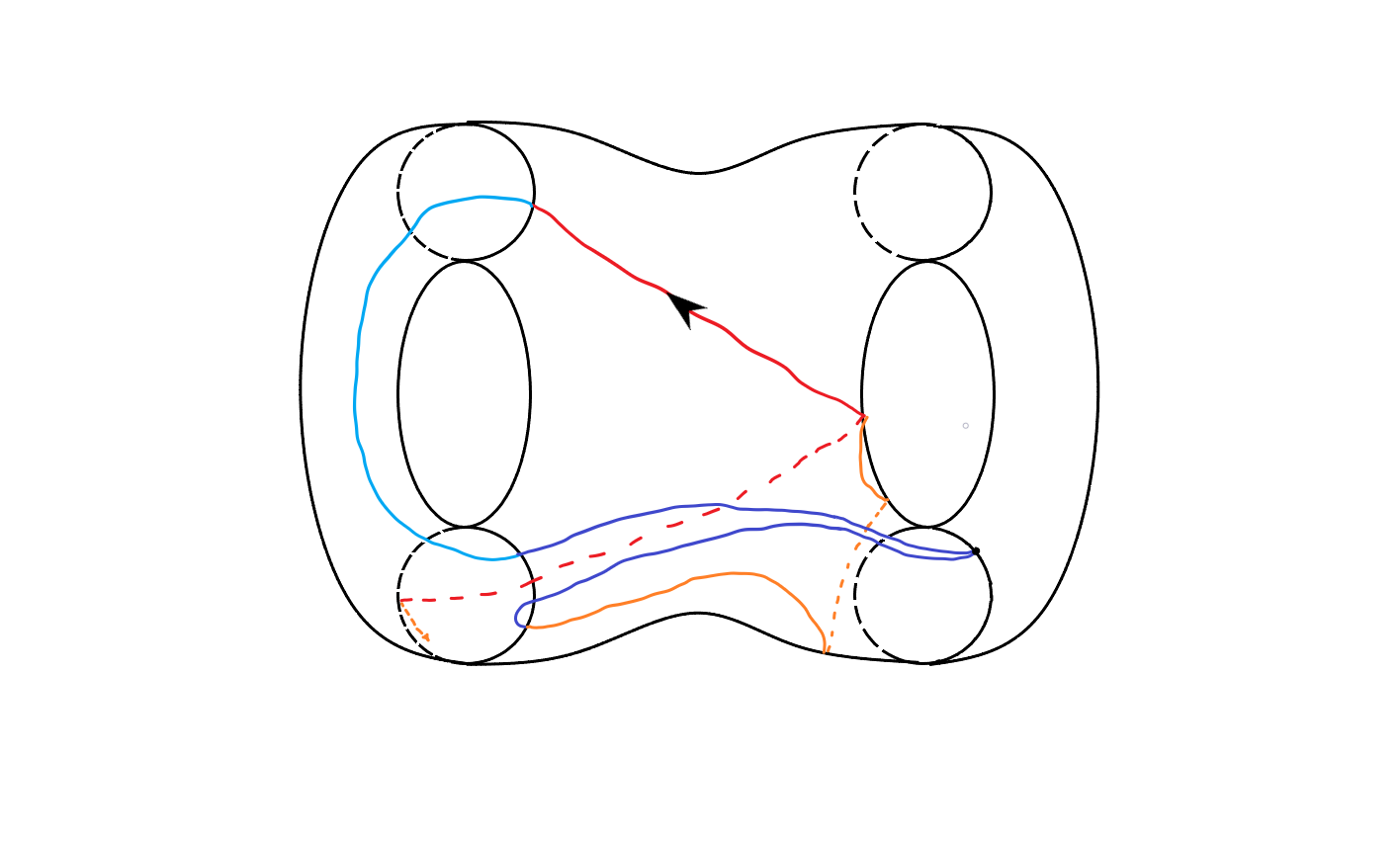}
  \caption{}
\end{subfigure}

\caption{All the possible closings of the previous cases. }
\label{fig:fig_closed}
\end{figure}

\begin{figure}

\begin{subfigure}{.5\textwidth}
  \centering
  \includegraphics[width=\textwidth]{mossa_intro_9.png}
  \caption{With a \(b_l\).}
\end{subfigure}
\begin{subfigure}{.5\textwidth}
  \centering
  \includegraphics[width=\textwidth]{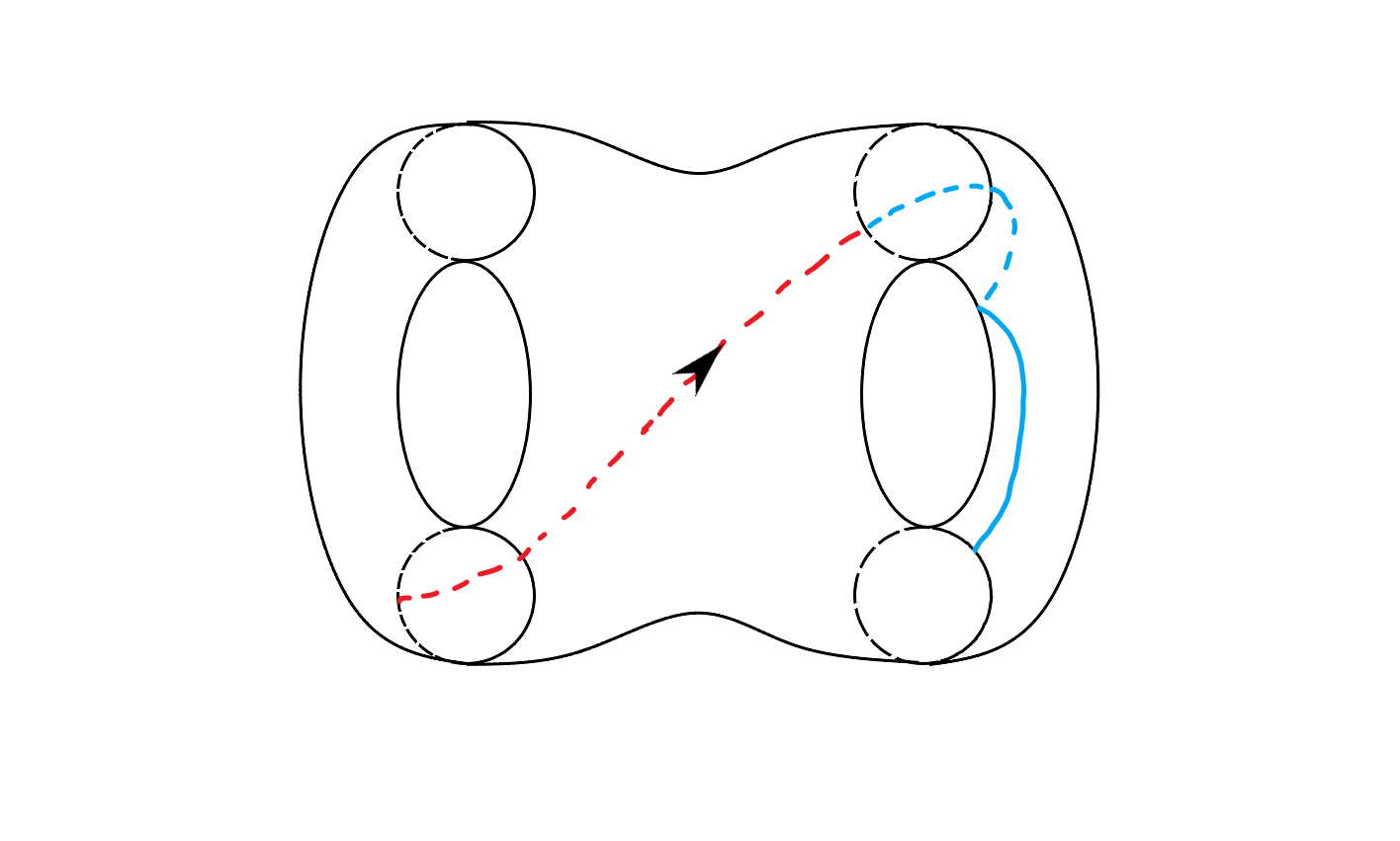}
  \caption{With a \(b_l^{-1}\).}
\end{subfigure}

\begin{subfigure}{.5\textwidth}
  \centering
  \includegraphics[width=\textwidth]{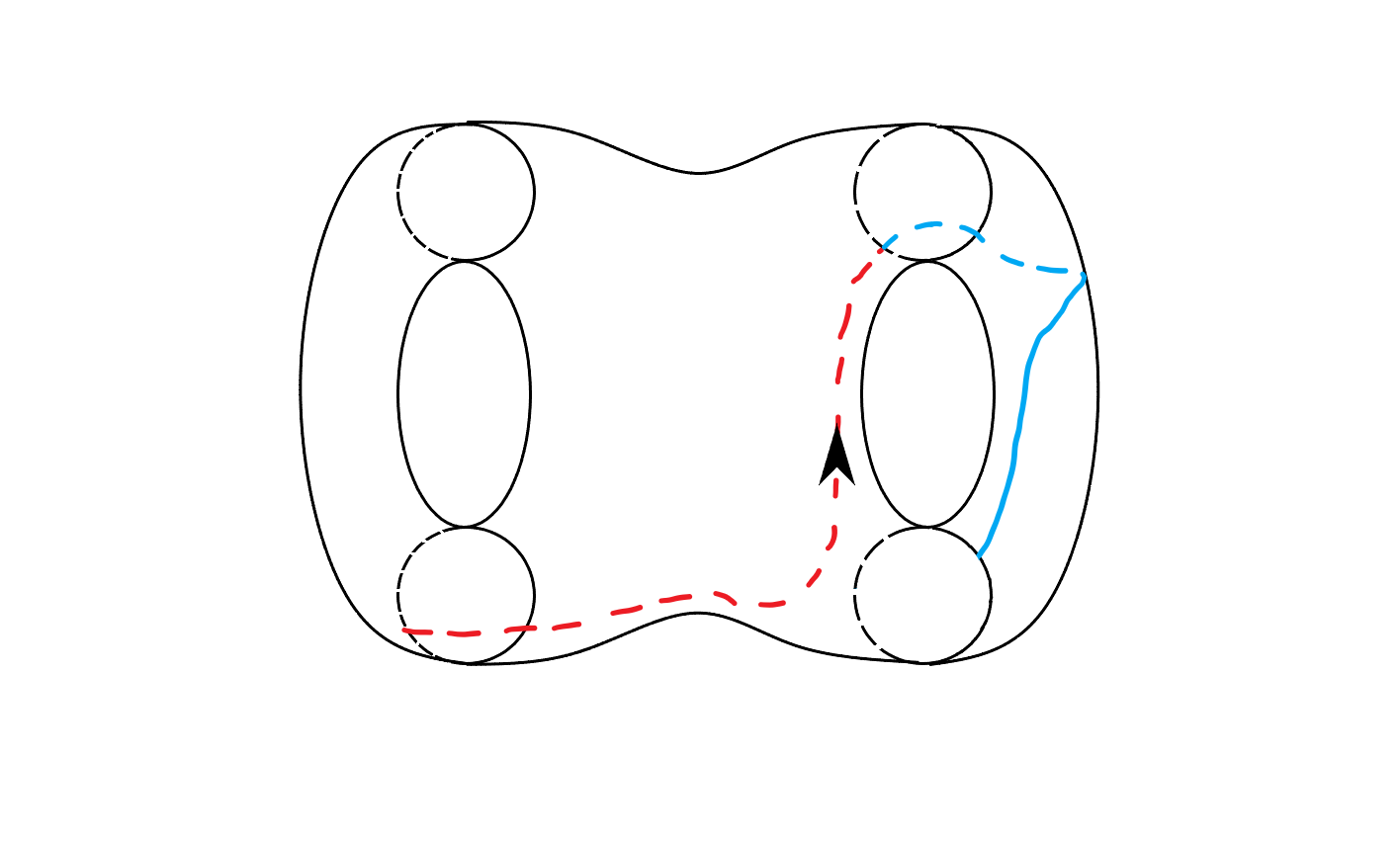}
  \caption{}
\end{subfigure}
\begin{subfigure}{.5\textwidth}
  \centering
  \includegraphics[width=\textwidth]{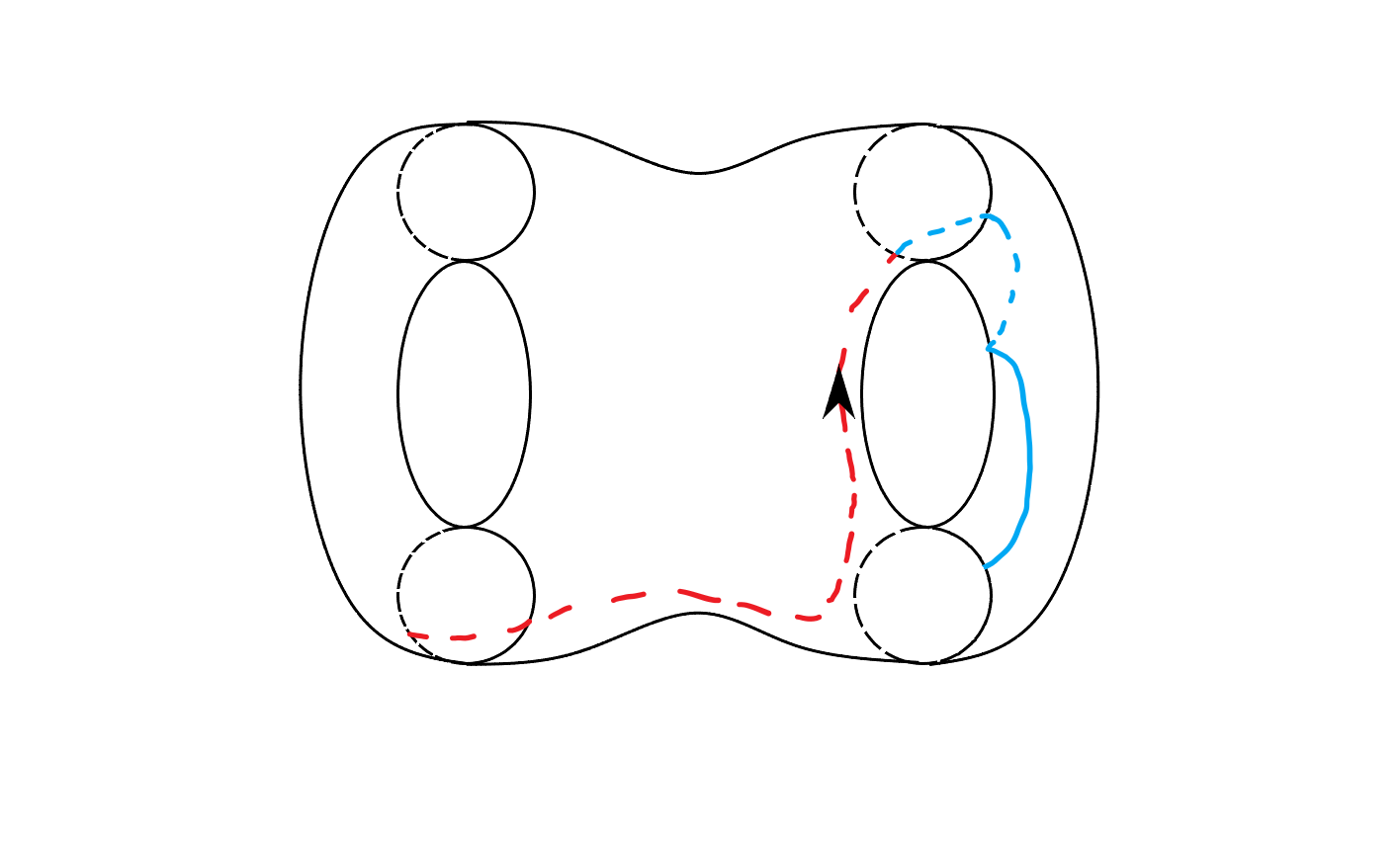}
  \caption{}
\end{subfigure}

\begin{subfigure}{.5\textwidth}
  \centering
  \includegraphics[width=\textwidth]{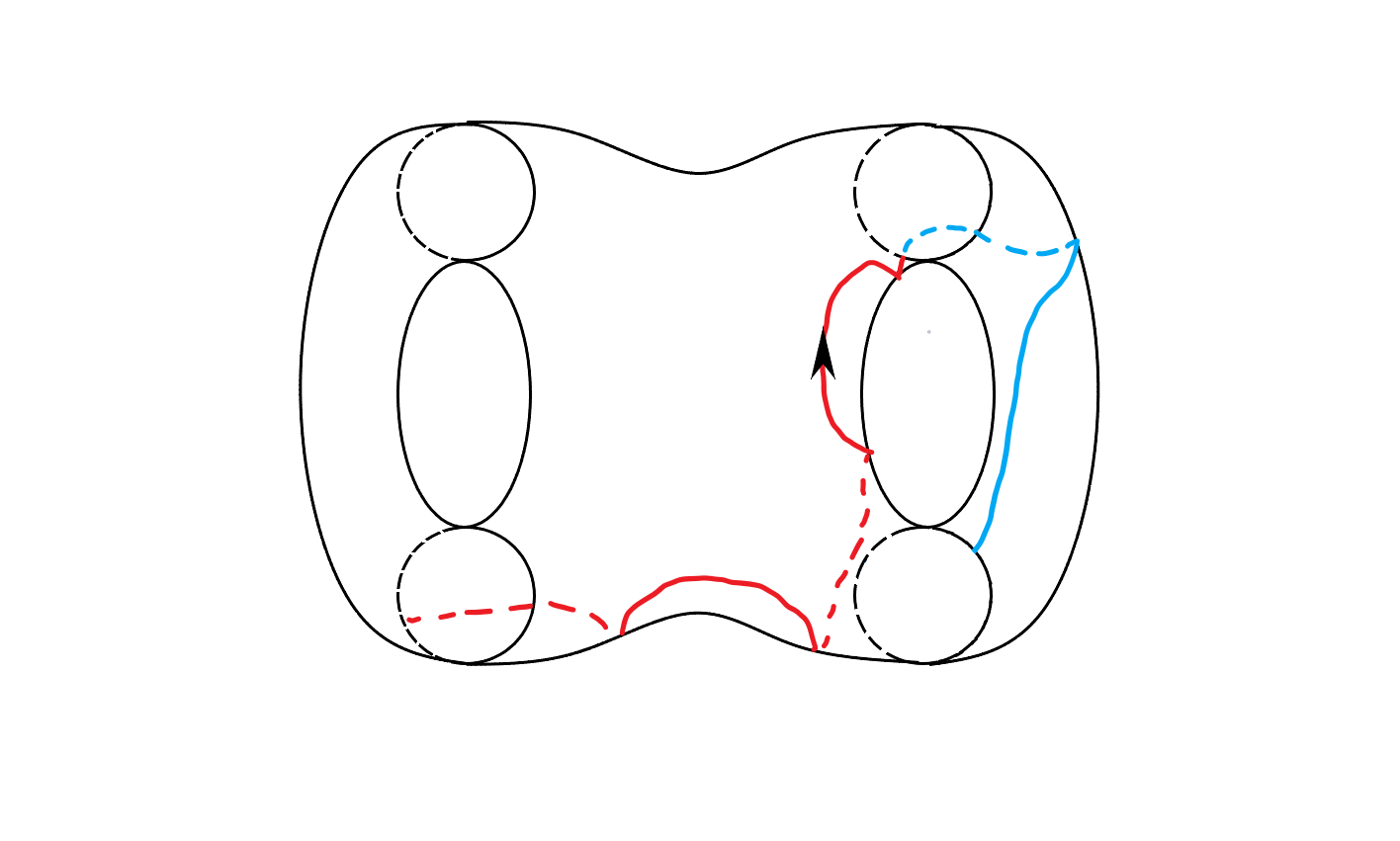}
  \caption{}
\end{subfigure}
\begin{subfigure}{.5\textwidth}
  \centering
  \includegraphics[width=\textwidth]{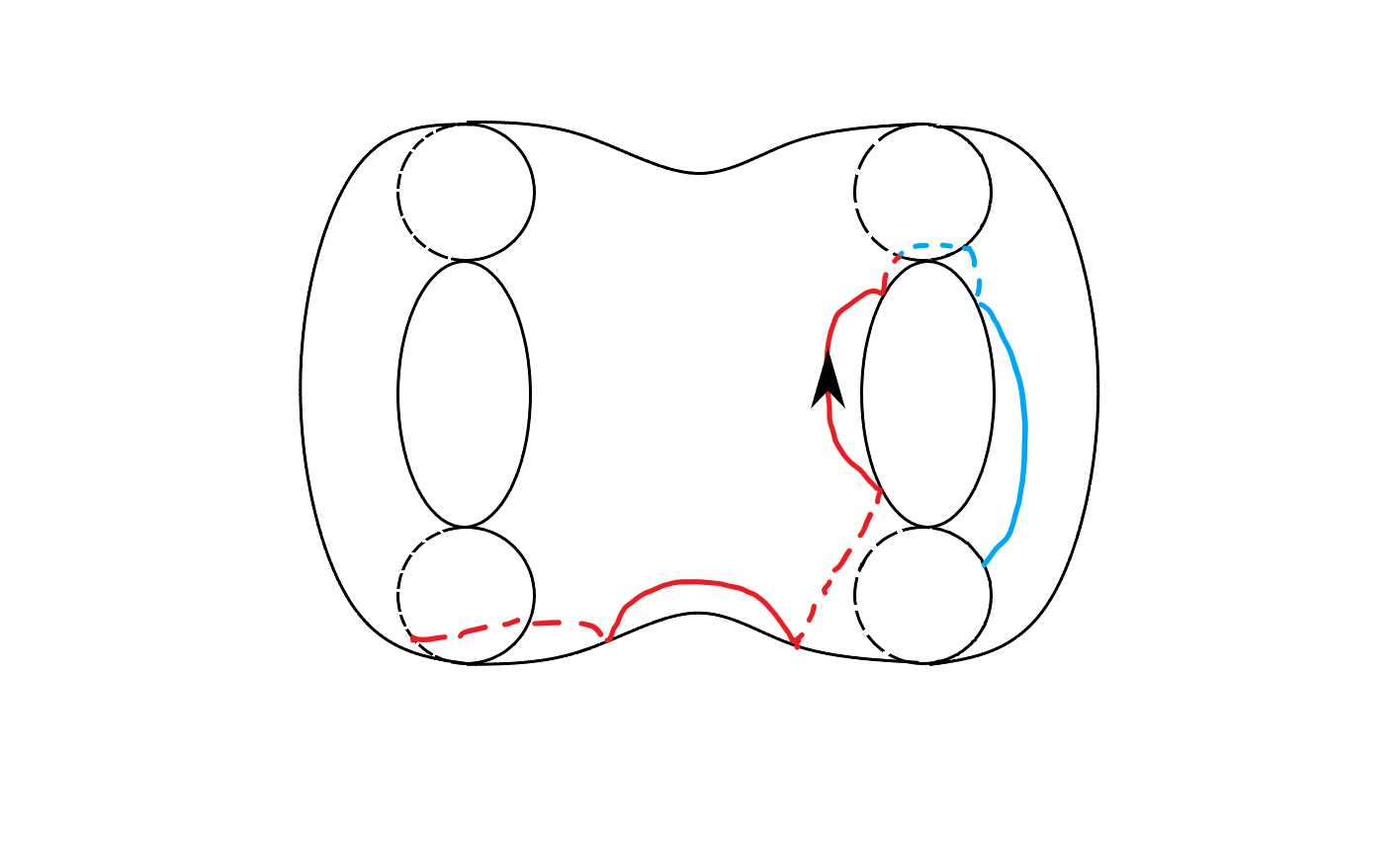}
  \caption{}
\end{subfigure}

\begin{subfigure}{.5\textwidth}
  \centering
  \includegraphics[width=\textwidth]{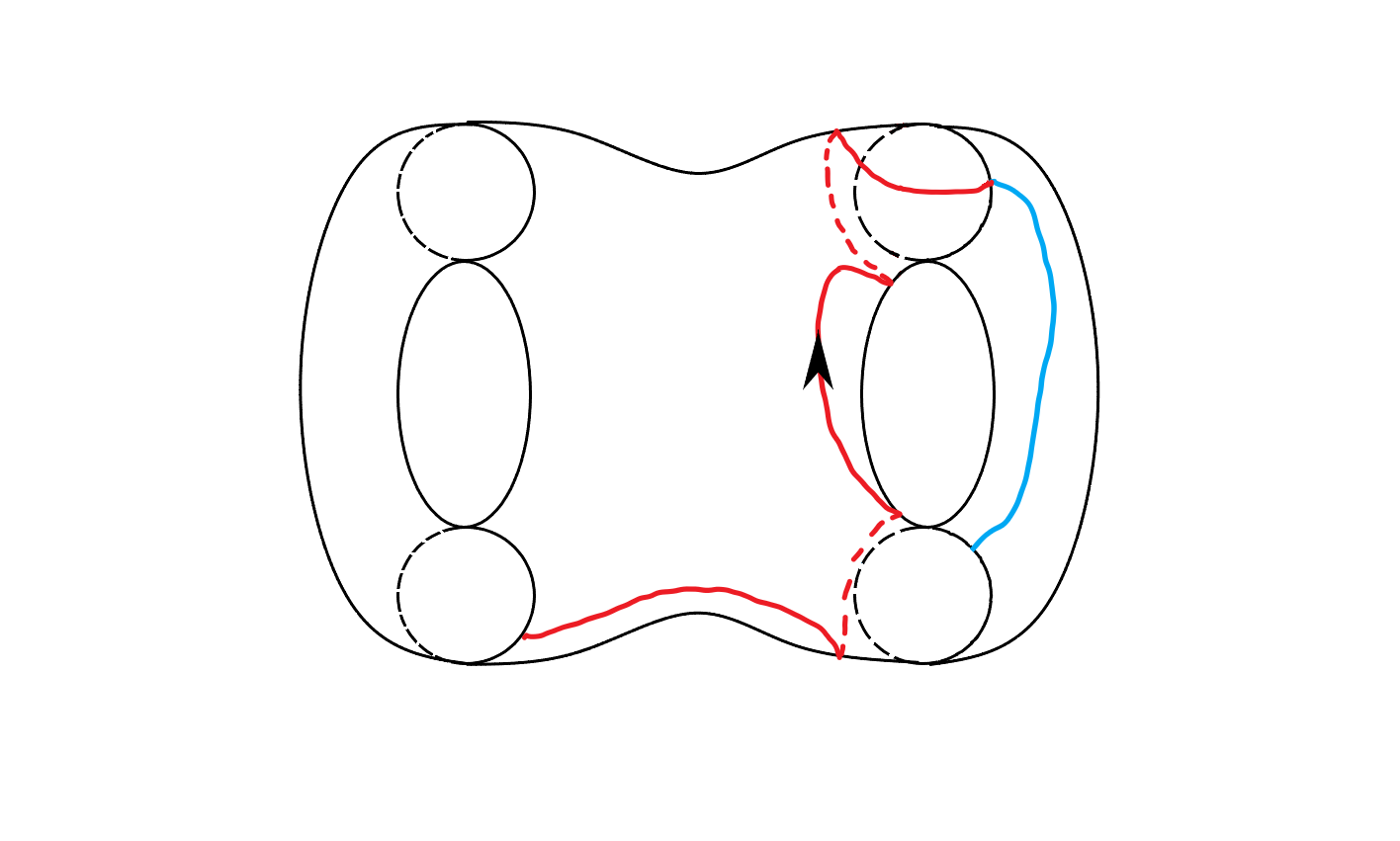}
  \caption{Ending situation: there is still the \(b_l\).}
\end{subfigure}
\begin{subfigure}{.5\textwidth}
  \centering
  \includegraphics[width=\textwidth]{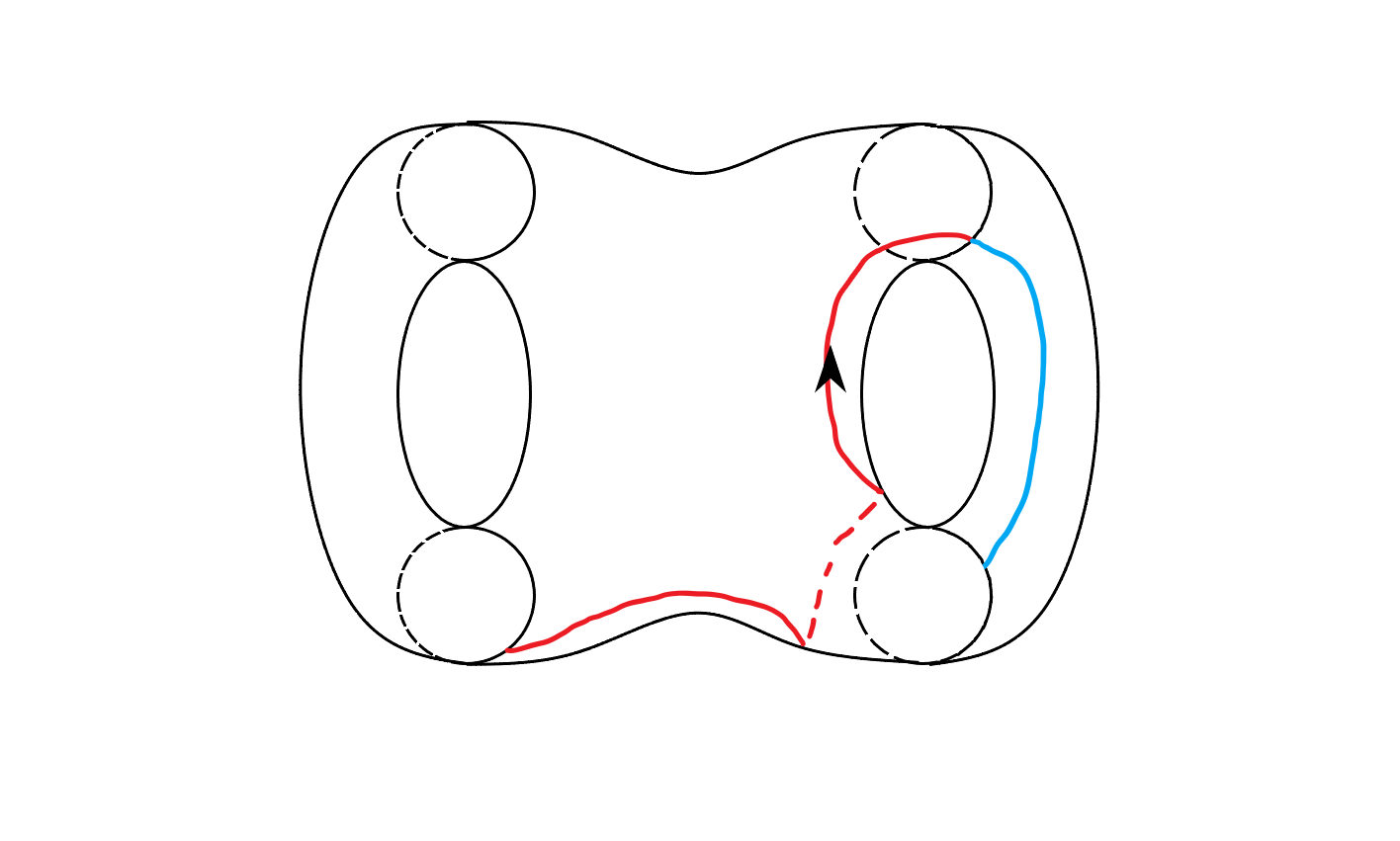}
  \caption{Ending situation: there is no \(b_l^{-1}\).}
\end{subfigure}

    \caption{The case (h).}
    \label{fig:case_11}
\end{figure}

\begin{teorema}
Let $f$ be an admissible $6$-tuple. There exists an algorithm depending only on $f$ that compute the $psl$-slide moves corresponding to the Heegaard splitting associated to the rich Heegaard diagram $R_f$. 
\end{teorema}
\begin{dimostrazione}
It is clear, using the previous construction, that we can determine, starting from \(f\) and in an algorithmic way, the curves \(x^0_f, x^1_f\) and \(x^2_f\). So to prove the statement we have only to describe, starting from \(x_f^i\), how to compute the words $\chi^i_f\in B_{2,2},\> i = 0, 1, 2$. Clearly the procedure is the same for each curve.\\
First of all we orient \(x_f^i\) as follows.  In case \ref{prop:dati_1}, \ref{prop:dati_2}, \ref{prop:dati_3} (resp.  \ref{prop:dati_4}) we orient the curve so that a red arc exit from the last (resp. first), vertices on $C_{Br}$. If the curve has no intersection with  $C_{Br}$ then we orient the curve so that a red arc exit from the last (resp. first), vertices on $C_{Bl}$. In fact at least two of three curves must intersect each \(C_i\) otherwise it would not be a proper Heegaard diagram \cite{gagliardi1981extending}. 
In Fig.\ref{fig:fig_intro_1} and Fig.\ref{fig:fig_intro} we depict all the unwound different cases emerging from considering, along $x^i_f$, a couple of consecutive red and blue arcs: clearly in each case a possible winding of the blue arc in both the directions may appear (see Fig.\ref{fig:fig_intro_1_a} and Fig.\ref{fig:fig_intro_1_b}). We associate to each of the previous cases a closed loop by connecting the arc to a based point in a standard way through the violet or green arc (see Fig.\ref{fig:fig_1}, \ref{fig:fig_1_a}, \ref{fig:fig_1_b} and Fig.\ref{fig:fig_closed}). Now we interpret each loop as an element of $B_{2,2}$ and call them elementary words. \\
We start with the first case, see Fig.\ref{fig:fig_intro_1} and Fig.\ref{fig:fig_1}: clearly the word associated is \(a_l\). The cases in Fig.\ref{fig:fig_intro_1_a} (resp. \ref{fig:fig_intro_1_b}) differ from the previous one for an outer (resp. inner) winding: we obtain \(a_l b_l\) (resp. \(a_l b_l^{-1}\)) with the isotopy realized along the dashed disk as in  Fig.\ref{fig:fig_1_a} (resp. Fig.\ref{fig:fig_1_b}). Note that the winding part affects the latter part of the word. \\
In general the same happens in the other cases, except for the last four in which the situation needs to be studied carefully. The words that we obtain, without winding, are the following (see Fig.\ref{fig:fig_intro} and Fig.\ref{fig:fig_closed}): 
\begin{enumerate}
    \item[(a)] \(a_r\) 
    \item[(b)] \(a_r^{-1}\) 
    \item[(c)] \(a_l^{-1}\) 
    \item[(d)] the green one is \(a_r^{-1}\), the violet one \(a_l\)
    \item[(e)] the green one is \(a_r\), the violet one \(a_l^{-1}\)
    \item[(f)] the green one is \(a_l\), the violet one \(a_r^{-1}\)
    \item[(g)] the green one is \(a_l^{-1}\), the violet one \(a_r\).
\end{enumerate}
Now we discuss the last four cases that are affected differently from the winding: 
\begin{enumerate}
    \item[(h)] In this case (see Fig.\ref{fig:case_11}) the word depends on the blue arc of the previous couple, so different cases arise due to this feature. More precisely, we have \(b_l b_r^{-1} a_r^{-1} b_r\) if the blue arc is of type \(b_r^{-1}\) and previous blue arc is of type \(b_l^{-1}\); otherwise we delete the first \(b_l\) or the last \(b_r\) if the corresponding requirement is not fulfilled
    \item[(i)] analogously to the previous case, we have \(b_r^{-1} b_r a_l^{-1} b_l^{-1}\) if the blue arc is of type \(b_l\) and the previous blue arc is of type \(b_r\); otherwise we delete the first \(b_r^{-1}\) or the last \(b_l^{-1}\) if the corresponding requirement is not fulfilled 
    \item[(j)] \(b_r a_l^{-1} b_l^{1}\) if the blue arc is of type \(b_l\); otherwise we delete the last \(b_l^{-1}\) 
    \item[(k)] \(b_l b_r^{-1} a_l\) if the previous blue arc is of type \(b_l^{-1}\); otherwise we delete the first \(b_l\).  
\end{enumerate}
Traveling along \(x_f^i\) and chaining the elementary words we obtaing the word representing \(\chi_f^i \in B_{2,2}\), \(i = 0, 1, 2\). 
\qed
\end{dimostrazione}

\section{Program and examples}
\label{examples}
In this last section, we describe the algorithm implemented in c++ which realizes the procedure described in the previous sections. The source code can be found on \href{https://github.com/Paolo-Cavicchioli/6_tuple_representation}{GitHub}. Moreover we report the results of the computation for a family of examples taken from \cite{bandieri2010census}. \\ 

The implementation is straightforward. After identifying among the different cases the one our \(6\)-tuple belongs in, all the variables are initialized using Proposition \ref{prop:dati}. Then, after computing all the connections and identifications, the software firstly keeps track of the words related to the single arcs then finds the three \(x_f^i\). The words \(\chi_f^i\) are then computed following each \(x_f^i\). If during the construction phase a different amount of \(x_f^i\) is found it means that \(f \notin \mathcal{F}\), so the program exits with an error.\\

We used the program to compute the plat-slide moves associated to the 6-tuples described in \cite{bandieri2010census} with associated orientable manifold. 
The results are collected in the following tables, where we used the following notation: 
\begin{itemize}
    \item \(\mathbb{S}^3/G\) is the quotient space of \(\mathbb{S}^3\) by the action of the group \(G\); the involved groups are: 
    \begin{itemize}
        \item \(Q_{4n} =\> < x, y\> |\> x^2 = (xy)^2 = y^n >\)
        \item \(D_{2^k(2n+1)} =\> < x, y\> |\> x^{2^k} = 1,\> y^{2n+1} = 1,\> xyx^{-1} = y^{-1}>\), with \(k \geq 3, n \geq 1\)
        \item \(P_{24} =\> < x, y\> |\> x^2 = (xy)^3 = y^3,\> x^4 = 1>\)
        \item \(P_{48} =\> < x, y\> |\> x^2 = (xy)^3 = y^4,\> x^4 = 1>\) 
        \item \(P_{120} =\> < x, y\> |\> x^2 = (xy)^3 = y^5,\> x^4 = 1>\) 
        \item \(P'_{3^k 8} =\> < x, y, z\> |\> x^2 = (xy)^2 = y^2,\> zxz^{-1} = y,\> zyz^{-1} = xy,\> z^{3^k} = 1>\), with \(k \geq 2\) 
        \item \(\mathbb{Z}_n\), the cyclic group over \(n\) elements, 
    \end{itemize}
    or direct products of these groups. 
    \item \(\mathbb{E}^3/G\) is the quotient space of \(\mathbb{E}^3\) by the action of the group \(G\), which is described using the International Tables for Crystallography \cite{hahn1983international}. 
    \item \((S, (\alpha_1, \beta_1), \dots, (\alpha_n, \beta_n))\) is the Seifert  fibered space whose orbit space is the surface \(S\) and having \(n\) exceptional fibers (with non-normalized parameters). 
    \item \(TB(A)\), with \(A \in GL(2; \mathbb{Z})\), is the torus bundle over \(\mathbb{S}^1\) with monodromy induced by \(A\).
    \item \(H_1 \cup_A H_2\) is the graph manifold obtained by gluing two Seifert manifolds \(H_1, H_2\), with \(\partial H_i = T, \> i = 1, 2\), along their boundary tori by means of the attaching map associated to matrix \(A\).
    \item \(Q_i(p,q)\) denotes the closed manifold obtained as the \((p,q)\) Dehn filling of the compact manifold \(Q_i\), whose interior is one of the 11 hyperbolic manifolds of finite volume with a single cusp and complexity at most 3 (see \cite{matveev2007algorithmic}). 
\end{itemize}

\begin{table}[]

\caption{}
\end{table}

\begin{table}[]
%
\caption{}
\end{table}

\begin{table}[]
%
\caption{}
\end{table}

\begin{table}[]
%
                                                           \\ \hline
\end{tabular}
\caption{}
\end{table}

\newpage

\printbibliography
\end{document}